\documentclass[final,3p,times]{elsarticle}

\usepackage{amssymb}
\usepackage{algorithm}
\usepackage{amsmath,amssymb,amsthm,mathtools}
\usepackage{newtxtext,newtxmath}
\usepackage{microtype}

\usepackage{amsthm}
\theoremstyle{plain}
\newtheorem{theorem}{Theorem}[section]
\newtheorem{proposition}{Proposition}[section]

\newtheorem{assumption}{Assumption}[section]
\newtheorem{lemma}[theorem]{Lemma}
\theoremstyle{definition}
\newtheorem{definition}[theorem]{Definition}

\theoremstyle{remark}
\newtheorem{rem}[theorem]{Remark}  
\newtheorem{Example}[theorem]{Example}  

\usepackage{caption}
\usepackage{subcaption}
\usepackage[colorlinks=true,allcolors=blue!55!black]{hyperref}

\usepackage{tikz}
\usetikzlibrary{positioning,arrows.meta,calc,fit,backgrounds}
\usetikzlibrary{arrows.meta,calc,matrix,positioning,fit}
\usepackage{xcolor}
\usepackage{booktabs}
\usepackage{enumitem}
\usepackage{newtxtext,newtxmath}
\usepackage{booktabs,tabularx,array}
\usepackage{microtype}
\usepackage{algorithm}
\usepackage{float}
\usepackage{algpseudocode}
\definecolor{colorA}{RGB}{78,121,167}
\definecolor{colorB}{RGB}{242,142,43}
\definecolor{colorC}{RGB}{89,161,79}
\definecolor{colorD}{RGB}{117,107,177}
\definecolor{lightA}{RGB}{205,220,238}
\definecolor{lightB}{RGB}{250,220,187}
\definecolor{lightC}{RGB}{205,230,200}
\definecolor{lightgray}{RGB}{242,242,242}

\newcommand{\R}{\mathbb{R}}
\newcommand{\Y}{\mathcal{Y}}
\newcommand{\U}{\mathcal{U}}
\newcommand{\Pmap}{\mathcal{P}}
\newcommand{\Fmap}{\mathcal{F}}
\newcommand{\Res}{\mathcal{R}}
\newcommand{\diff}{\mathrm{D}}

\newcommand{\ip}[2]{\left\langle #1,#2\right\rangle}
\newcommand{\norm}[1]{\left\lVert #1\right\rVert}

\newcommand{\al}{\alpha}
\newcommand{\be}{\beta}

\newcommand{\supp}{\operatorname{supp}}
\newcommand{\nnz}{\operatorname{nnz}}

\newcommand{\Dmat}{\mathsf D}

\newcommand{\Eplus}{\mathcal E_G^+}
\newcommand{\Patch}{\mathcal P}

\journal{Computer Methods in Applied Mechanics and Engineering}

\begin{document}

\begin{frontmatter}



\title{Block–Sparse Gauss–Newton Optimization of Finite Basis Physics-Informed Neural Networks: Sparse Automatic Differentiation Approach}


\author[label1]{Denis Korolev}
\ead{denis.korolev@wias-berlin.de}

\author[{label2a,label2b}]{Adrian Hill}
\ead{hill@tu-berlin.de}

\author[label3]{Jeeeun Lee}
\ead{jinlee7@kaist.ac.kr}

\affiliation[label1]{
    organization={Weierstrass Institute for Applied Analysis and Stochastics},
    addressline={Anton-Wilhelm-Amo-Str. 39},
    city={Berlin},
    postcode={10117},
    country={Germany}
}

\affiliation[label2a]{
    organization={BIFOLD -- Berlin Institute for the Foundations of Learning and Data},
    addressline={Franklinstr. 28/29},
    city={Berlin},
    postcode={10587},
    country={Germany}
}

\affiliation[label2b]{
    organization={Machine Learning Group, Technical University of Berlin},
    addressline={Franklinstr. 28/29},
    city={Berlin},
    postcode={10587},
    country={Germany}
}

\affiliation[label3]{
    organization={KAIST InnoCORE PRISM-AI Center, Korea Advanced Institute of Science and Technology (KAIST)},
    addressline={291 Daehak-ro, Yuseong-gu},
    city={Daejeon},
    postcode={34141},
    country={Republic of Korea}
}

\begin{abstract}
In this work, we develop an efficient Gauss--Newton optimization approach for finite basis physics-informed neural networks (FBPINNs). The key ingredient is the overlapping domain decomposition architecture of FBPINNs, which localizes the influence of each subnetwork and induces block sparsity in both the residual Jacobian and the corresponding Gauss--Newton matrix. We characterize this structure and exploit it through sparse automatic differentiation, using graph coloring to partition structurally orthogonal Jacobian columns or rows. This reduces differentiation costs and memory requirements while preserving all structurally nonzero derivative information. The resulting compressed Jacobian representations allow us to assemble the Gauss--Newton system directly in block-sparse form, enabling the use of sparse direct solvers and preconditioned iterative methods. In this way, the otherwise expensive second-order optimization problem is reduced to scalable sparse linear algebra. We demonstrate the approach on multiscale PDEs, where accurate approximation requires rich neural network parametrizations and large numbers of collocation points. Numerical experiments show that our block-sparse Gauss--Newton method attains high accuracy and substantially accelerates optimization relative to first-order baselines while remaining computationally tractable for large problem sizes.
\end{abstract}




\end{frontmatter}



\section{Introduction}

Training physics-informed neural networks (PINNs) consists in parametrizing the solution of a partial differential equation (PDE) by a neural network and incorporating the resulting PDE residuals into a least-squares objective \cite{raissi2019physics}. This residual-based objective is further evaluated over a set of collocation points $\{x_{1},\ldots,x_{N_{\Omega}}\}$ sampled throughout the computational domain $\Omega \subset \mathbb{R}^{d}$. 
In this work, motivated by applications in computational mechanics, 
we focus on spatial dimensions $d \leq 3$, with an emphasis on multiscale PDEs, for which a large number $N_{\Omega}$ of collocation points may be required to resolve the microstructural features. The approximation is then obtained by minimizing the objective with respect to the network parameters $\theta \in \mathbb{R}^{P}$. First-order methods are computationally attractive for this optimization task but often struggle to attain high accuracy. A major reason is the severe ill-conditioning of the resulting PINN optimization problem, which first-order methods are generally not well suited to address \cite{de2023operator,rathore2024challenges,krishnapriyan2021characterizing,cao2025analysis}. This motivates the use of suitable preconditioning strategies, in particular curvature-aware methods such as Gauss--Newton \cite{jnini2026curvature,jnini2025gauss}, which can also be viewed as belonging to the family of natural gradient methods \cite{muller2023achieving,muller2024function}.
It has been shown that the Gauss--Newton approach can substantially outperform first-order methods such as Adam \cite{kingma2014adam}
on a range of benchmark problems; cf. \cite{jnini2026curvature, jnini2025gauss}. 

The use of curvature-aware methods, however, may come at substantial
computational cost. To make this precise, let
$\mathrm R(\theta)\in\mathbb{R}^{N_{\mathrm{eq}}N_{\Omega}}$
denote the discrete residual vector, where $N_{\mathrm{eq}}$ is the
number of scalar PDE residual equations evaluated at each collocation
point. The associated residual Jacobian
$J_{\mathrm R}(\theta)\in\mathbb{R}^{(N_{\mathrm{eq}}N_{\Omega})\times P}$
gives rise to the Gauss--Newton Gramian
\begin{align}
G(\theta):=J_{\mathrm R}(\theta)^\top J_{\mathrm R}(\theta)
\in\mathbb{R}^{P\times P}.
\end{align}
At iteration $k$, the regularized Gauss--Newton direction
$\delta\theta_k\in\mathbb{R}^P$ is determined by solving
\begin{align}\label{eq:GN_system}
\big(G(\theta_k)+\mu I_P\big)\delta\theta_k
=-J_{\mathrm R}(\theta_k)^\top\mathrm R(\theta_k),
\end{align}
where $\mu>0$ is a regularization parameter. The network parameters
are subsequently updated according to
\begin{align}
\theta_{k+1}=\theta_k+\eta_k\delta\theta_k,
\end{align}
with step size $\eta_k>0$. Thus, among other tasks, such as the potentially expensive inversion of the regularized Gramian
in \eqref{eq:GN_system}, the realization of the Gauss--Newton approach requires repeated construction of the residual Jacobian using automatic differentiation (AD). This is computationally demanding, since evaluation of the PDE residual already involves derivatives of the neural network with respect to its inputs, while construction
of the residual Jacobian requires further differentiation of these quantities with respect to the network parameters. The resulting nested differentiation can produce a large computational graph, leading to substantial overhead in terms of memory and computational performance. We remind the reader that forward-mode AD evaluates a matrix-free Jacobian--vector product (JVP)  $v \mapsto J_{\mathrm R}(\theta)v$ (see e.g., \cite{griewank2008evaluating}).
Evaluating such a JVP with the $i$-th standard basis vector materializes the $i$-th column of the Jacobian,
such that constructing the full Jacobian column-by-column requires $P$ JVP evaluations, one for each network parameter. 
Conversely, reverse-mode AD computes matrix-free vector--Jacobian products (VJPs) of the form $w \mapsto J_{\mathrm R}(\theta)^\top w$, so row-wise materialization requires $N_{\mathrm{eq}}N_{\Omega}$ VJP evaluations.
In addition to these derivative passes, storing the numerical entries of the dense Jacobian requires memory proportional to $N_{\mathrm{eq}}N_{\Omega}P$,
not to mention the substantial additional memory overhead associated with storing the computational graph and intermediate activations.
These costs can become prohibitive for large-scale neural PDE discretizations, as in the case of multiscale PDEs, where resolving fine-scale structures typically requires both large numbers of collocation points and increasingly expressive (and therefore large) network parametrizations.
Consequently, Jacobian construction and storage constitute major bottlenecks for scaling Gauss--Newton training to such problems.  


\begin{figure}[H]
    \centering
    \includegraphics[width=0.8\linewidth]{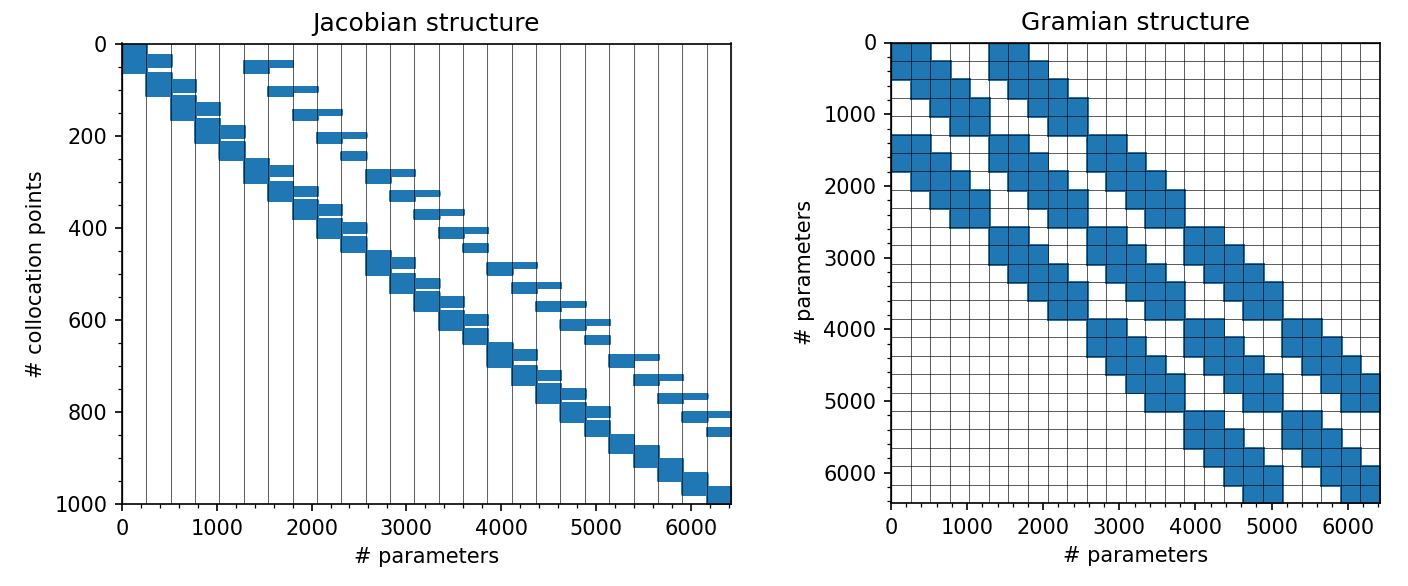}
    \caption{Block sparsity patterns of the FBPINN residual Jacobian (left)
and Gauss--Newton Gramian (right).
Each Jacobian column block
(highlighted by vertical lines)
corresponds to a local subnetwork, with nonzero entries restricted
to collocation points within its support. The Gramian inherits
a symmetric block sparsity pattern determined by interactions
between overlapping subnetworks.}
    \label{fig:structural_sparsity}
\end{figure}

%

Finite-basis physics-informed neural networks (FBPINNs) were introduced as an efficient domain decomposition neural network architecture, particularly for multiscale PDEs \cite{moseley2023finite}; see also \cite{anderson2026elm, botvinick2026ab, salvado2026multi, dolean2024multilevel, lee2026fbpinn, van2026local} for related variants and applications.  Instead of representing the solution by a single global network, FBPINNs combine multiple local subnetworks defined over an overlapping domain decomposition through compactly supported partition-of-unity windows. This localization reduces the complexity and effective frequency content that must be represented by each individual subnetwork, which then effectively helps mitigate spectral bias \cite{moseley2023finite, lee2026fbpinn}. In this work, we exploit a second consequence of the FBPINN construction: the compactly supported window functions localize each subnetwork, and hence the influence of its associated parameter block, to a small region of the computational domain. This induces block sparsity in the residual Jacobian and, consequently, in the Gauss--Newton Gramian, as shown in Fig.\ref{fig:structural_sparsity}, enabling sparse techniques for assembling and solving the resulting Gauss--Newton systems. 

In this work, we develop an efficient computational workflow for Gauss--Newton optimization of FBPINNs by carrying the resulting sparsity structure through AD and into the assembly and solution of the Gauss--Newton system.
The starting point is the \emph{structural orthogonality} of columns encoded in the block sparsity pattern of the residual Jacobian.
Jacobian columns are structurally orthogonal if they have no nonzero entries in a common row; analogously, rows are structurally orthogonal if they have no nonzero entries in a common column; cf. \cite{gebremedhin2005color}.
In FBPINNs, structural
orthogonality arises naturally from the domain decomposition: the Jacobian is organized into column blocks corresponding
to local subnetworks, each assigned to an individual subdomain.
Columns corresponding to nonoverlapping subdomains are
structurally orthogonal, while rows associated with different collocation points are structurally orthogonal whenever their
residuals do not depend on any common subnetwork.
Sparse AD exploits this property by materializing each group of structurally orthogonal columns or rows in a single JVP (or VJP) \cite{curtis1974estimation}.
Owing to structural orthogonality, their structural nonzeros occupy disjoint positions, allowing these contributions to be packed into a single column or row of a compressed Jacobian, from which the individual nonzero entries can be recovered.
Finding such groups is a classical graph-coloring problem determined by the Jacobian sparsity pattern \cite{gebremedhin2005color, griewank2008evaluating, coleman1983estimation}.
We adapt these techniques to FBPINNs by deriving an explicit optimal block coloring for forward-mode AD and employing greedy row coloring for reverse-mode AD. The resulting compressed Jacobians are then used to assemble the Gauss--Newton Gramian directly in block-sparse form,
avoiding the decompression of intermediate Jacobians.
This reduces differentiation and memory costs compared with conventional AD, while the resulting sparse systems can be treated using sparse direct solvers or preconditioned iterative methods.

\subsection{Contributions}

Our contributions are summarized as follows:
\begin{itemize}
    \item Starting from the continuous Gauss--Newton natural-gradient formulation, we derive its collocation discretization and characterize the locality-induced block sparsity of the Jacobian and Gramian, providing a basis for coloring-based sparse AD and the application of sparse numerical solvers.

    \item We adapt graph-coloring techniques to FBPINNs for sparse forward- and reverse-mode AD, using optimal periodic block coloring and greedy row coloring, respectively. The resulting compressed Jacobians enable direct block-sparse assembly of the Gauss--Newton Gramian without having to materialize a decompressed Jacobian.

    \item We discuss sparse direct and preconditioned iterative solvers and benchmark the method on three multiscale problems: Poisson with multiscale forcing, elliptic equations with oscillatory coefficients in a first-order system least-squares (FOSLS) formulation. Here, our approach enables Gauss--Newton optimization with up to a quarter of a million neural network parameters and half a million collocation points on a single GPU, while offering substantial potential for further GPU cluster parallelization.
\end{itemize}

\subsection{Related work}

We note that the block sparsity induced by the FBPINN ansatz has already been exploited in randomized variants of FBPINNs, where the hidden network parameters are randomly initialized and subsequently kept fixed; see, e.g., \cite{anderson2026elm, van2026local, shang2025overlapping}. In such settings, the resulting optimization problem is reduced to a structured least-squares problem whose system matrix is block-sparse, shifting the main computational challenge from nonlinear learning to numerical linear algebra. While computationally appealing, freezing the internal network parameters restricts the approximation capabilities and limits the applicability in settings where the representation itself or additional unknown quantities must be learned, as it is usual in applications of PINNs. However, keeping all network parameters trainable fundamentally changes the computing workflow, as the residual Jacobian  must be repeatedly constructed using automatic differentiation, and Gramian must be efficiently assembled.

Other approaches to efficient curvature-aware PINN optimization include
Kronecker-factored curvature approximations \cite{dangel2024kronecker},
randomized numerical linear algebra techniques for preconditioning
and efficient linear solves \cite{guzman2026improving,bioli2027accelerating},
and two-sided Jacobian sketching with SVD-based solves and adaptive
regularization \cite{webb2026optimisation}. Although these techniques substantially improve the efficiency of
Gauss--Newton optimization, they do not specifically exploit the localization-induced block sparsity of FBPINNs, which are particularly well suited to large-domain and multiscale PDEs. A different approach is pursued in \cite{korolev2026tensor}, where
coordinate separability of tensor-structured neural approximation classes is exploited to derive
exact compressed Gauss--Newton formulations. The present work follows the same structural principle, but exploits spatial locality in domain-decomposition-based neural classes and the resulting block sparsity, rather than coordinate separability and the
associated low-rank structure. Finally, local feature filtering, developed for domain-decomposed random feature methods \cite{van2026local} to improve conditioning, could in principle also be adapted to fully trainable FBPINNs, but this possibility is not investigated here. 

\subsection{Paper outline}
The remainder of the paper is organized as follows. Section~2 develops the block-sparse Gauss--Newton formulation and characterizes the sparsity structure of the Jacobian and Gramian. Section~3 presents coloring-based sparse AD and block-sparse Gramian assembly. Section~4 discusses sparse direct and iterative solvers, while Section~5 reports numerical experiments. The appendices provide details on the FBPINN construction and proofs of the coloring results.

\section{Block-Sparse Gauss--Newton FBPINN optimization}

We use the Spivak convention for derivatives. If $F:\R^m\to\R^q$ is differentiable, then
\begin{align}
\diff F(x):\R^m\to\R^q
\end{align}
is the derivative as a linear map, and $\diff F(x)[v]$ denotes its action on a direction $v\in\R^m$. If $J_F(x)$ is the matrix representation of this map in the standard Euclidean bases, then
\begin{align}
\diff F(x)[v]=J_F(x)v.
\end{align}
Throughout the paper, $0_n$ denotes the zero vector in $\R^n$, and $0_{m\times n}$ denotes the zero matrix in $\R^{m\times n}$. For a matrix $A$, $A_{i,:}$ and $A_{:,j}$ denote its $i$-th row and $j$-th column. For an index set $I$, $A_{i,I}$ denotes the restriction of row $i$ to the columns in $I$. We use the standard definition of the support of a function:
\begin{align*}
\supp(f):= \overline{\{x \, : \, f(x)\neq 0\}}.
\end{align*}
For  functions in $L^2(\Omega)$, support is understood in the essential sense. For simplicity, we present the analysis for a scalar PDE residual, i.e., $N_{\mathrm{eq}}=1$. The extension to systems of PDEs
is straightforward and is discussed in
Section~\ref{subsec:oscillatory-mixed}.

\subsection{Block-sparse Gauss--Newton geometry: from locality to sparsity}
Let $\Omega
=
\prod_{j=1}^{d}
(x_j^{\min},x_j^{\max})
\subset\R^d$ be a bounded domain with boundary $\partial\Omega$,
and let $\U$ be a suitable function space admitting boundary traces.
For a differential operator $\mathcal A:\U\to L^2(\Omega)$,
we consider the boundary value problem
\begin{equation}\label{Sec 1: PDE}
\mathcal A u=f
\quad\mathrm{in}\ L^2(\Omega),
\qquad
u|_{\partial\Omega}=g
\quad\mathrm{in}\ L^2(\partial\Omega).
\end{equation}
Here, we assume that $\mathcal A$ is linear.
For nonlinear operators, the linearization can be applied (see e.g. \cite{jnini2025gauss}), and
the methodology applies analogously.

We seek to approximate the solution of \eqref{Sec 1: PDE} by a suitably
parametrized neural network $u_\theta\in\U$ with trainable parameters
$\theta\in\R^P$. Throughout, we assume that
$u_\theta|_{\partial\Omega}=g$, i.e., that the neural networks
satisfy the boundary condition exactly; cf. \cite{sukumar2022exact}. Following the FBPINN approach \cite{moseley2023finite}, let
$\{\Omega_\al\}_{\al=1}^K$ be an overlapping cover of $\Omega$ with overlap ratio $\delta>1$. Furthermore, let
$\{\omega_\al\}_{\al=1}^K$ be window functions
forming a partition of unity subordinate to this cover:
\begin{align}\label{eq: partition of unity}
0\leq\omega_\al\leq1,
\qquad
\sum_{\al=1}^K\omega_\al=1
\quad\mathrm{on}\ \Omega,
\qquad
\supp(\omega_\al)\subseteq\Omega_\al.
\end{align}
We use normalized tensor-product windows
constructed from compactly supported squared raised-cosine profiles,
following \cite{dolean2024multilevel}. ~\ref{app1: finite basis networks}
defines the overlap ratio and gives the explicit construction. The FBPINN ansatz then takes the form
\begin{align}\label{eq:local-fbpinn-contributions}
u_\theta
=\sum_{\al=1}^K u_{\theta^\al}^\al,
\qquad
u_{\theta^\al}^\al(x)
:=\omega_\al(x)\,z_\al(x;\theta^\al),
\end{align}
where $z_\al(\,\cdot\,;\theta^\al):\Omega\to\R$ is the subnetwork assigned to
$\Omega_\al$, with trainable parameters $\theta^\al\in\R^n$.
Since $\supp(\omega_\al)\subseteq\Omega_\al$, each window confines the
contribution of its respective subnetwork to the corresponding subdomain. We assume here that the activation functions in the subnetworks and the window
functions are chosen such that $u_\theta$ is sufficiently smooth. 

The parameters are further collected block-wise as
\begin{align}\label{eq: parameter vector}
\theta=(\theta^1,\ldots,\theta^K),
\qquad
\theta^\al\in\R^n,
\qquad
P=Kn.
\end{align}
The associated \textit{parametrization} and nonlinear model class are given by
\begin{align}\label{eq: parametrization}
\Pmap:\R^P\to\U,
\qquad
\theta\mapsto\Pmap(\theta)=u_\theta,
\qquad
\mathcal{NN}_{\mathrm{FB}}
:=\operatorname{range}(\Pmap).
\end{align}
To approximate \eqref{Sec 1: PDE}, we aim to minimize a least-squares residual
over $\mathcal{NN}_{\mathrm{FB}}$. This requires the
\textit{residual operator}
\begin{align}
\Fmap:\U\to L^2(\Omega),
\qquad
\Fmap(u):=\mathcal A u-f,
\end{align}
and the parameter-to-residual map
\begin{align}\label{eq:continuous-residual-map}
\Res:=\Fmap\circ\Pmap:\R^P\to L^2(\Omega),
\qquad
\Res(\theta)=\mathcal A u_\theta-f.
\end{align}
The continuous FBPINN optimization problem then reads
\begin{align}\label{eq:continuous-pinn-loss}
\min_{\theta\in\R^P}\mathfrak{L}(\theta),
\qquad
\mathfrak{L}(\theta)
:=\frac12\norm{\Res(\theta)}_{L^2(\Omega)}^2
=\frac12\norm{\mathcal A u_\theta-f}_{L^2(\Omega)}^2.
\end{align}

We further note that the parametrization \eqref{eq: parametrization} is nonlinear, whereas
its differential $\diff\Pmap(\theta):\R^P\to\U$ is a linear map.
The parameter grouping \eqref{eq: parameter vector} induces the
block-column representation
\begin{align}\label{eq:fbpinn-parametrization-differential}
\diff\Pmap(\theta)
=\begin{bmatrix}
\diff_{\theta^1}\Pmap(\theta)&\cdots&
\diff_{\theta^K}\Pmap(\theta)
\end{bmatrix},
\end{align}
where each block is identified with the quasi-matrix of its
function-valued columns:
\begin{align}\label{eq:fbpinn-parametrization-blocks}
\diff_{\theta^\al}\Pmap(\theta)
:=\begin{bmatrix}
\partial_{\theta_1^\al}u_\theta&\cdots&
\partial_{\theta_n^\al}u_\theta
\end{bmatrix}
:\R^n\to\U.
\end{align}
By the chain rule, the residual push-forward is given by
\begin{align}\label{eq:continuous-residual-differential-general}
\diff\Res(\theta)
=\diff\Fmap(u_\theta)\circ\diff\Pmap(\theta)
=\mathcal A\,\diff\Pmap(\theta)
:\R^P\to\Y.
\end{align}
Defining the local residual blocks $\mathcal B_\al(\theta)
:=\mathcal A\,\diff_{\theta^\al}\Pmap(\theta)
:\R^n\to L^2(\Omega)$ yields the following block representation
\begin{align}\label{eq:continuous-residual-pushforward-blocks}
\diff\Res(\theta)
=\begin{bmatrix}
\mathcal B_1(\theta)&\cdots&\mathcal B_K(\theta)
\end{bmatrix}.
\end{align}
The Gauss--Newton Gramian is then given by
\begin{align}\label{eq:continuous-gramian}
\mathcal G(\theta)
:=\diff\Res(\theta)^*\diff\Res(\theta)
\in\R^{P\times P},
\end{align}
where $\diff\Res(\theta)^*$ denotes the Hilbert adjoint; cf. \cite{de2023operator}. The block
structure \eqref{eq:continuous-residual-pushforward-blocks} induces a
$K\times K$ block form with $n\times n$ block sizes,
\begin{align}\label{eq:continuous-gramian-blocks}
\mathcal G_{\al\be}(\theta)
=\mathcal B_\al(\theta)^*\mathcal B_\be(\theta), \qquad
\mathcal G_{\al\be}(\theta)_{ps}
=\ip{\mathcal A\,\partial_{\theta_p^\al} \, u_\theta}
{\mathcal A\,\partial_{\theta_s^\be} \, u_\theta}_{L^2(\Omega)},
\quad 1\leq p,s\leq n.
\end{align}
The loss gradient is
$\nabla\mathfrak{L}(\theta)=\diff\Res(\theta)^*\Res(\theta)$, and the
damped Gauss--Newton update is obtained from
\begin{align}\label{eq:continuous-natural-gn}
\bigl(\mathcal G(\theta_k)+\mu I_P\bigr)\delta\theta_k
=-\nabla\mathfrak{L}(\theta_k),
\qquad
\theta_{k+1}=\theta_k + \eta_k\delta\theta_k,
\end{align}
where $\mu>0$ is the damping parameter, $I_P \in \mathbb{R}^{P \times P}$ is the identity matrix,
and $\eta_k>0$ is a step size.

\begin{assumption}[Local differential operator]
\label{ass:local-differential-operator}
The operator $\mathcal A:\U\to L^2(\Omega)$ is a continuous linear
differential operator of finite order $m$, i.e.,
\begin{align}
\mathcal A v
=
\sum_{|\nu|\leq m}
a_\nu\,D^\nu v,
\qquad v\in\U,
\end{align}
where $a_\nu\in C(\overline{\Omega})$ for all multi-indices
$\nu$ with $|\nu|\leq m$.
\end{assumption}
This assumption is primarily used to establish that the differential operator
does not enlarge the support of a function.
\begin{lemma}[Locality]\label{lem:locality}
Let Assumption~\ref{ass:local-differential-operator} hold, let $v\in\U$, and
let $O\subseteq\Omega$ be open. Then
\begin{align}\label{locality_FBPINNs}
v=0 \quad \text{a.e.\ in } O
\qquad\Longrightarrow\qquad
\mathcal A v=0 \quad \text{a.e.\ in } O.
\end{align}
\end{lemma}

\begin{proof}
Fix $\nu$ with $|\nu|\le m$ and let $\varphi\in C_c^\infty(O)$. The extension of $\varphi$ by zero lies in $C_c^\infty(\Omega)$ and is
therefore admissible in the definition of the weak derivative $D^\nu v$ on
$\Omega$, so that
\begin{align}\label{eq:weak-def}
\int_\Omega (D^\nu v)\,\varphi\,dx
=(-1)^{|\nu|}\int_\Omega v\,D^\nu\varphi\,dx .
\end{align}
As $\supp\varphi\subseteq O$, the right-hand integrand vanishes outside $O$,
and since $v=0$ a.e.\ in $O$,
\begin{align}
\int_\Omega v\,D^\nu\varphi\,dx
=\int_O v\,D^\nu\varphi\,dx
=0 .
\end{align}
With \eqref{eq:weak-def} this gives $\int_O (D^\nu v)\varphi\,dx=0$ for every
$\varphi\in C_c^\infty(O)$, hence $D^\nu v=0$ a.e.\ in $O$ by the fundamental
lemma of the calculus of variations. Since $\nu$ was arbitrary and each
$a_\nu$ is bounded as a continuous function on the compact set
$\overline{\Omega}$, it follows that
$\mathcal A v=\sum_{|\nu|\le m}a_\nu D^\nu v=0$ a.e.\ in $O$.
\end{proof}
We have the following structural characterization of the continuous Gauss--Newton matrix.
\begin{proposition}\label{lem:fbpinn-support-locality}
Let Assumption~\ref{ass:local-differential-operator} hold. Then the following inclusions hold
\begin{align}\label{eq:fbpinn-support-locality}
\supp\!\left(
\mathcal A\,\partial_{\theta_p^\al}u_\theta
\right)
\subseteq
\supp\!\left(
\partial_{\theta_p^\al}u_\theta
\right)
\subseteq
\supp(\omega_\al).
\end{align}
for every $\al\in\{1,\ldots,K\}$ and $p=1,\ldots,n$. In particular, 
\begin{align}\label{eq:continuous-gramian-locality}
|\supp(\omega_\al)\cap\supp(\omega_\be)|=0
\quad\Longrightarrow\quad
\mathcal G_{\al\be}(\theta)=0_{n\times n},
\end{align}
where $|\cdot|$ denotes Lebesque measure.
\end{proposition}

\begin{proof}
Let
$x\notin \supp(\partial_{\theta_p^\al}u_\theta)$.
There exists an open neighborhood $O\subseteq\Omega$ of $x$ such that
$\partial_{\theta_p^\al}u_\theta=0$ in $O$. 
Assumption~\ref{ass:local-differential-operator} and Lemma \ref{lem:locality} imply that
$\mathcal A\,\partial_{\theta_p^\al}u_\theta=0$ in $O$, which proves the
first inclusion in \eqref{eq:fbpinn-support-locality}. Furthermore, \eqref{eq:local-fbpinn-contributions} gives
\begin{align}\label{eq:fbpinn-tangent-functions}
\partial_{\theta_p^\al}u_\theta
=
\omega_\al\,\partial_{\theta_p^\al}
z_\al(\,\cdot\,;\theta^\al).
\end{align}
Thus, $\omega_\al=0$ implies
$\partial_{\theta_p^\al}u_\theta=0$, and the second inclusion in \eqref{eq:fbpinn-support-locality} follows. Finally, the integrand in \eqref{eq:continuous-gramian-blocks} vanishes almost everywhere outside $\supp(\omega_\al)\cap\supp(\omega_\be)$, which proves \eqref{eq:continuous-gramian-locality}.
\end{proof}
\noindent
Therefore, the locality of the FBPINN ansatz, combined with that of the differential operator, implies block sparsity of the continuous Gauss–Newton matrix. The geometry of the overlapping cover then determines which blocks can be nonzero.

In what follows, we call $\Omega_\alpha$ \textit{active at} $x\in\Omega$ if
$\omega_\al(x)\neq0$. Accordingly, we define the active set
\begin{align}\label{eq:structural-active-set-continuous}
T(x)
:=
\{\alpha\in\{1,\ldots,K\}:\omega_\al(x)\neq0\}.
\end{align}
Furthermore, we define $p_\delta:=\lceil\delta\rceil$, where
$\lceil\cdot\rceil$ denotes the ceiling function.
Lemma~\ref{lem:overlap-multiplicity} then establishes the geometric bound
$|T(x)|\leq p_\delta^d$ for any $x\in\Omega$, which we use throughout.

\subsection{Sparsity structure of the discrete Jacobian and Gramian}
\label{sec:discrete-jacobian-gramian-sparsity}

In this section, we discretize the continuous Gauss--Newton formulation \eqref{eq:continuous-natural-gn} and
derive the corresponding quantities required for numerical optimization.
We show that the discretization preserves the block-sparse structure induced
by the FBPINN decomposition. Along the way, we characterize the sparsity
pattern of the residual Jacobian and the block sparsity pattern of the
Gauss--Newton Gramian. These structures will be used in the subsequent
sections for sparse automatic differentiation and efficient assembly of the
Gauss--Newton system. Exploiting these structures, however, requires
introducing some technical machinery.

\begin{figure}[H]
\centering
\includegraphics[width=1.0\linewidth]{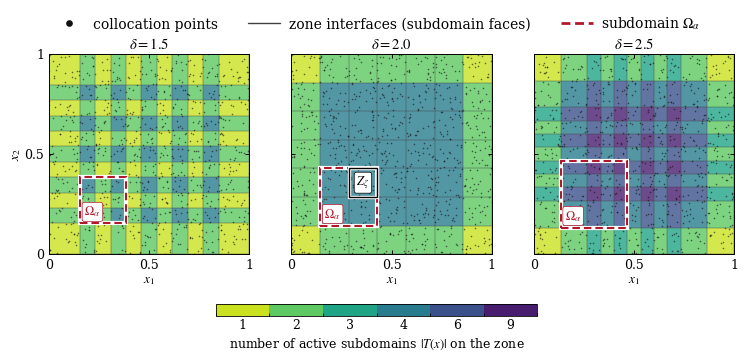}
\caption{Zonal sampling on a uniform $6\times6$ partition for three overlap
ratios $\delta$.}
\label{fig:zonal_sampling}
\end{figure}

For the discretization, we use a Monte--Carlo approach based on a
\textit{zonal sampling} strategy tailored to the overlapping geometry
of the FBPINN decomposition. The subdomain faces partition $\Omega$
into zones $\{Z_z\}_{z=1}^Z$, each with a fixed set of active subnetworks, see Fig. \ref{fig:zonal_sampling}.
Given a total budget of $N_\Omega$ collocation points, we distribute
the points proportionally to the zone volumes, rounding the allocations
to integers while preserving the total budget. Choosing $N_\Omega$
sufficiently large so that every zone is sampled, we draw $n_z$
points uniformly at random within each $Z_z$, where
\begin{align}\label{eq:zonal-allocation}
n_z \approx N_\Omega\frac{|Z_z|}{|\Omega|},
\qquad
\sum_{z=1}^{Z}n_z=N_\Omega.
\end{align}

\begin{rem}[Number of collocation points per subdomain]
\label{rem:zonal-rounding}
The proportional allocation also determines the number of points
in each subdomain. Indeed, $\Omega_\al$ is a union of zones, thus
\begin{align}\label{eq:zonal-subdomain-count}
\sum_{z:\,Z_z\subseteq\Omega_\al}n_z
=\frac{N_\Omega}{|\Omega|}
\sum_{z:\,Z_z\subseteq\Omega_\al}|Z_z|
+\mathcal O_{\delta,d}(1) =\frac{N_\Omega|\Omega_\al|}{|\Omega|}
+\mathcal O_{\delta,d}(1).
\end{align}
Here, the remainder accounts for the
difference between the integer counts $n_z$ and their proportional
quotas $N_\Omega|Z_z|/|\Omega|$.
\end{rem}
\noindent
Zonal sampling provides approximately uniform coverage across subdomains and, importantly, preserves the Jacobian sparsity pattern when collocation points are resampled within their respective zones.

We now introduce an important object that organizes much of the
bookkeeping. Based on
\eqref{eq:structural-active-set-continuous}, let
\begin{align}\label{eq: active set}
T_q
:= T(x_q^r)
= \{\al\in\{1,\ldots,K\}:\omega_\al(x_q^r)\neq0\}
\end{align}
be the \textit{active-subdomain set} of the $q$-th collocation point.
These sets are determined by evaluating the window functions at the
collocation points and can be computed before training.

Each point $x_q^r\in Z_z$ is assigned the quadrature weight
\begin{align}\label{eq:zonal-quadrature-weights}
w_q:=\frac{|Z_z|}{n_z},
\end{align}
so that each zone contributes to the numerical quadrature proportionally to its volume, and the normalization holds:
\begin{align}
\sum_{q=1}^{N_\Omega} w_q=
\sum_{z=1}^{Z}
\sum_{q:\,x_q^r\in Z_z}
\frac{|Z_z|}{n_z}
=
\sum_{z=1}^{Z}|Z_z|
=
|\Omega|.
\end{align}
We define the quadrature-weighted discrete residual by
\begin{align}\label{eq:normalized-discrete-residual}
\mathrm{R}(\theta)
:=
\begin{bmatrix}
\sqrt{w_1}\,\mathcal R(\theta)(x_1^r)\\
\vdots\\
\sqrt{w_{N_\Omega}}\,\mathcal R(\theta)(x_{N_\Omega}^r)
\end{bmatrix}
\in\R^{N_\Omega}.
\end{align}
The discrete FBPINN loss is therefore given by
\begin{align}\label{eq:normalized-discrete-loss}
 L(\theta)
:=\frac12\norm{\mathrm{R}(\theta)}_2^2
=\frac12\sum_{q=1}^{N_\Omega}
w_q\left|\mathcal A u_\theta(x_q^r)-f(x_q^r)\right|^2,
\end{align}
which is the zone-wise Monte--Carlo approximation of
\eqref{eq:continuous-pinn-loss}.

The differential
of the discrete residual is the linear map
\begin{align}\label{Sec 1.2: DDR}
\diff\mathrm{R}(\theta):\R^P\to\R^{N_\Omega},
\qquad
\diff\mathrm{R}(\theta)[v]
=
\begin{bmatrix}
\sqrt{w_1}\,\diff\mathcal R(\theta)[v](x_1^r)\\
\vdots\\
\sqrt{w_{N_\Omega}}\,\diff\mathcal R(\theta)[v](x_{N_\Omega}^r)
\end{bmatrix}.
\end{align}
The matrix representing $\diff\mathrm{R}(\theta)$ in the standard
Euclidean bases is the residual Jacobian
$J_{\mathrm R}(\theta)\in\R^{N_\Omega\times P}$. Computationally, obtaining this representation amounts to evaluating
$\diff\mathrm R(\theta)$ along suitable AD seed directions, which are discussed
in detail in Section~\ref{Sec: sparse AD}.
With the block ordering \eqref{eq: parameter vector}, the residual Jacobian
inherits the block-column form
\begin{align}\label{eq:block jacobian representation}
J_{\mathrm R}(\theta)
=
\begin{bmatrix}
\Phi_1 & \cdots & \Phi_K
\end{bmatrix},
\end{align}
where each $\Phi_\al\in\R^{N_\Omega\times n}$ is the collocated counterpart
of the corresponding local residual block in
\eqref{eq:continuous-residual-pushforward-blocks}:
\begin{align}\label{eq: Jacobian block}
\Phi_\al
:=
\begin{bmatrix}
\sqrt{w_1}\,
\mathcal A\partial_{\theta_1^\al}u_\theta(x_1^r)
&
\cdots
&
\sqrt{w_1}\,
\mathcal A\partial_{\theta_n^\al}u_\theta(x_1^r)
\\
\vdots
&
&
\vdots
\\
\sqrt{w_{N_\Omega}}\,
\mathcal A\partial_{\theta_1^\al}u_\theta(x_{N_\Omega}^r)
&
\cdots
&
\sqrt{w_{N_\Omega}}\,
\mathcal A\partial_{\theta_n^\al}u_\theta(x_{N_\Omega}^r)
\end{bmatrix}.
\end{align}
Respectively, the discrete Gauss--Newton matrix is given by
\begin{align}\label{eq:discrete-gn-objects}
G(\theta):=J_{\mathrm R}(\theta)^\top J_{\mathrm R}(\theta)
\in\R^{P\times P}.
\end{align}
Owing to \eqref{eq:block jacobian representation}, it consists of $K \times K$ blocks:
\begin{align}\label{eq:block-gramian-form}
G(\theta)
=
\begin{bmatrix}
\Phi_1^\top\Phi_1&\cdots&\Phi_1^\top\Phi_K\\
\vdots&\ddots&\vdots\\
\Phi_K^\top\Phi_1&\cdots&\Phi_K^\top\Phi_K
\end{bmatrix},
\end{align}
where each $n\times n$ block is the zone-wise Monte--Carlo approximation
of \eqref{eq:continuous-gramian-blocks}.

We proceed by determining the sparsity pattern of the Jacobian
\eqref{eq:block jacobian representation}. Throughout, a block or entry is called
\emph{structurally nonzero} if it is admitted by the sparsity pattern induced
by the window supports; its numerical value may still vanish for particular
parameter values. The sparsity pattern of $J_{\mathrm R}$ is defined as the set
of structurally nonzero row--column positions, denoted by
\[
\operatorname{Sp}(J_{\mathrm R})
\subseteq
\{1,\ldots,N_\Omega\}\times\{1,\ldots,P\}.
\]
Furthermore, let
\begin{align}\label{eq:jacobian-column-blocks}
\mathcal C_\al
:=
\{(\al-1)n+1,\ldots,\al n\}
\end{align}
be the columns of $J_{\mathrm R}$ occupied by the block $\Phi_\al$. For the $q$-th row of $J_{\mathrm R}$, define its block support by
\begin{align}
\operatorname{bsupp}\bigl((J_{\mathrm R})_{q,:}\bigr)
:=
\left\{
\al\in\{1,\ldots,K\}:
(\Phi_\al)_{q,:}
=
(J_{\mathrm R})_{q,\mathcal C_\al}
\text{ is structurally nonzero}
\right\}.
\end{align}
For collocation points obtained using the zonal sampling, Proposition~\ref{lem:fbpinn-support-locality} yields the following chain of implications
\begin{align}
\al\notin T_q
&\quad\Longrightarrow\quad
x_q^r\notin \supp(\omega_\al)
\quad\Longrightarrow\quad
x_q^r\notin
\supp\!\left(
\mathcal A\,\partial_{\theta_p^\al}u_\theta
\right)\quad\Longrightarrow\quad
\bigl(\mathcal A\,\partial_{\theta_p^\al}u_\theta\bigr)(x_q^r)=0
\end{align}
for each $p=1,\ldots,n$.  Invoking \eqref{eq: Jacobian block}, we conclude that
$(\Phi_\al)_{q,:}=0_{1\times n}$. Treating each active subnetwork as a dense
parameter block therefore gives
\begin{align}\label{eq: structural equivalence of active sets}
\operatorname{bsupp}\bigl((J_{\mathrm R})_{q,:}\bigr)
=
T_q,
\end{align}
which identifies the Jacobian blocks that are structurally present in row $q$. The Jacobian sparsity pattern is then given by
\begin{align}\label{eq:scalar-sparsity-proof}
\operatorname{Sp}(J_{\mathrm R})
=
\bigcup_{q=1}^{N_\Omega}
\bigcup_{\al\in T_q}
\{q\}\times\mathcal C_\al.
\end{align}
Observe that it is completely determined by the active-subdomain sets
$T_q$ and the block column size $n$.

We now shift our attention to the discrete Gramian matrix. Its sparsity is governed by co-occurrence of subdomains in the active
sets. For a subdomain $\al$ and a pair $(\al,\be)$, define the
\textit{active-row sets}
\begin{align}\label{eq: active-row set}
Q_\al:=\{q:\al\in T_q\},
\qquad
Q_{\al\be}:=Q_\al\cap Q_\be
=\{q:\al\in T_q\ \text{and}\ \be\in T_q\},
\end{align}
both of which are read off from $\operatorname{Sp}(J_{\mathrm R})$.
Writing $G_{\al\be}(\theta)=\Phi_\al^\top\Phi_\be$ row-wise and
using \eqref{eq: structural equivalence of active sets} gives
\begin{align}\label{eq:gramian-block-Q-note}
G_{\al\be}(\theta)
=\sum_{q\in Q_{\al\be}}
(\Phi_\al)_{q,:}^\top(\Phi_\be)_{q,:}.
\end{align}
In particular, $Q_{\al\be}=\emptyset$ implies
$G_{\al\be}(\theta)=0_{n\times n}$. The block sparsity pattern of
the Gramian is therefore given by
\begin{align}\label{eq:gramian-block-pattern-note}
\operatorname{Sp}_{\mathrm{blk}}(G)
=\{(\al,\be):Q_{\al\be}\neq\emptyset\}
=\bigcup_{q=1}^{N_\Omega}T_q\times T_q.
\end{align}
Thus, the Gramian block pattern is constructed by inserting all pairs
of subdomains that occur together in an active set.
Formula \eqref{eq:gramian-block-Q-note} also shows that assembling each
block requires only the rows in $Q_{\al\be}$ and the Jacobian entries
belonging to blocks $\al$ and $\be$. The latter will be utilized in our assembly process reliant on sparse AD.

\subsection{Geometric characterization of the FBPINN block graph}
\label{sub:block-graph}

We represent the Gramian block interactions \eqref{eq:gramian-block-pattern-note} by a graph determined by the overlapping decomposition. This graph provides a compact description of the sparsity structure and will be used for block-count estimates, sparse AD coloring, and the construction of preconditioners for sparse iterative solvers.

\begin{definition}[Block graph]
The block graph is the undirected graph
\begin{align}\label{eq:fbpinn-block-graph}
\mathfrak G_{\mathrm{blk}}
=
\bigl(\mathcal V_{\mathrm{blk}},\mathcal E_{\mathrm{blk}}\bigr),
\qquad
\mathcal V_{\mathrm{blk}}=\{1,\ldots,K\},
\end{align}
where $\{\al,\be\}\in\mathcal E_{\mathrm{blk}}$ if and only if
$\al\neq\be$ and $Q_{\al\be}\neq\emptyset$.
\end{definition}

To characterize the block graph geometrically, we identify each
subdomain label $\al$ with its grid multi-index
$a_\al\in\mathcal I_{K_0}$, using the row-major ordering
\eqref{lex ordering} defined in
\ref{app1: finite basis networks}. For $a,b\in\mathcal I_{K_0}$, define the Chebyshev distance between grid indices
\begin{align}\label{eq:chebyshev-distance}
d_\infty(a,b)
:=
\|a-b\|_\infty
=
\max_{1\leq j\leq d}|a_j-b_j|,
\end{align}
and, with Chebyshev radius
$r_\delta:=\lceil\delta\rceil-1$, the corresponding neighborhood
\begin{align}\label{eq:chebyshev-neighborhood}
\mathcal N_\delta(a_\al)
:=
\left\{
a\in\mathcal I_{K_0}:
\|a-a_\al\|_\infty\leq r_\delta
\right\}.
\end{align}
By the construction of the overlapping subdomains, the supports of
$\omega_\al$ and $\omega_\be$ overlap when
\begin{align}
|a_{\al,j}-a_{\be,j}|<\delta,
\qquad j=1,\ldots,d.
\end{align}
Since the differences $|a_{\al,j}-a_{\be,j}|$ are integers, this is
equivalent to
\begin{align}
\|a_\al-a_\be\|_\infty\leq r_\delta.
\end{align}
Moreover, the zonal sampling resolves every such overlap region.
Consequently, for $\al\neq\be$, we get
\begin{align}\label{eq:chebyshev-cooccurrence-condition}
Q_{\al\be}\neq\emptyset
\quad\Longleftrightarrow\quad
\|a_\al-a_\be\|_\infty\leq r_\delta.
\end{align}
Hence, the neighbors of $\al$ in $\mathfrak G_{\mathrm{blk}}$ are
the subdomains whose grid indices belong to
$\mathcal N_\delta(a_\al)\setminus\{a_\al\}$.
For $1<\delta\leq2$, one has $r_\delta=1$; in two dimensions, the
resulting block graph is therefore the king graph on the subdomain grid,
up to boundary truncation.

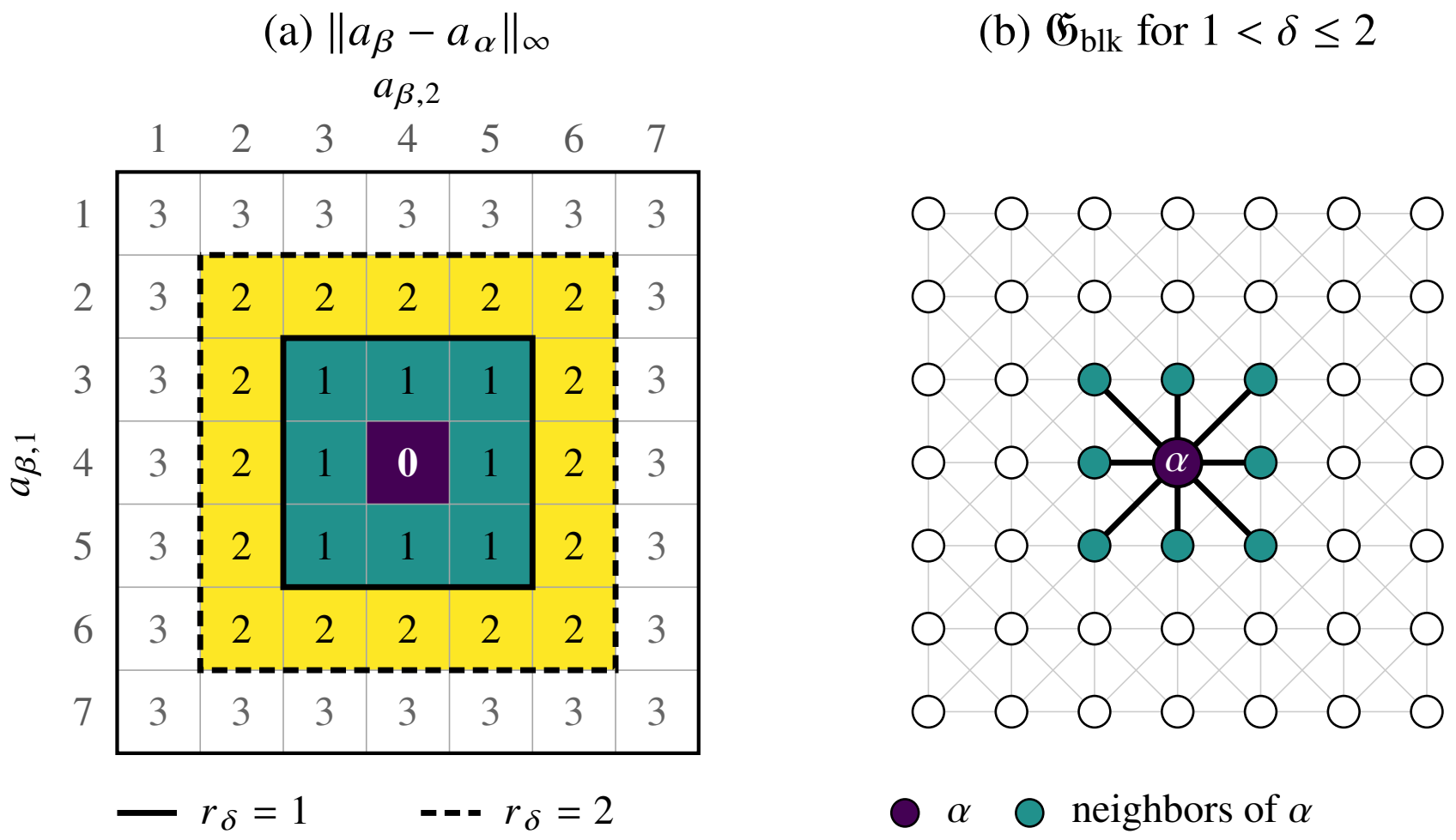
\begin{figure}[t]
\centering
\begin{tikzpicture}[font=\fontsize{9}{11}\selectfont]
\def\K{7}\def\h{0.640}\def\gap{5.930}\def\ci{4}\def\cj{4}
\pgfmathtruncatemacro{\Kp}{\K+1}
\pgfmathsetmacro{\rv}{max(0.19*\h,0.088)}
\pgfmathsetmacro{\ra}{max(0.29*\h,0.158)}
\pgfmathsetmacro{\rl}{max(0.165*\h,0.098)}
\definecolor{dz}{HTML}{440154}  
\definecolor{done}{HTML}{21918C}
\definecolor{dtwo}{HTML}{FDE725}
\tikzset{
  panel title/.style={font=\fontsize{10}{12}\selectfont, anchor=base},
  axis tick/.style={font=\fontsize{9}{11}\selectfont, text=black!65},
  nbhd one/.style={line width=1.3pt},
  nbhd two/.style={line width=1.3pt, dash pattern=on 3.4pt off 1.8pt},
}

\begin{scope}
  \foreach \i in {1,...,\K}{\foreach \j in {1,...,\K}{
    \pgfmathtruncatemacro{\d}{max(abs(\i-\ci),abs(\j-\cj))}
    \ifcase\d \def\fc{dz}\or \def\fc{done}\or \def\fc{dtwo}\else \def\fc{white}\fi
    \fill[\fc] ({(\j-1)*\h},{-(\i-1)*\h}) rectangle ++(\h,-\h);
    \ifcase\d
      \node[text=white] at ({(\j-0.5)*\h},{-(\i-0.5)*\h}) {$\mathbf{0}$};
    \or
      \node[text=black] at ({(\j-0.5)*\h},{-(\i-0.5)*\h}) {$\d$};
    \or
      \node[text=black] at ({(\j-0.5)*\h},{-(\i-0.5)*\h}) {$\d$};
    \else
      \node[text=black!60] at ({(\j-0.5)*\h},{-(\i-0.5)*\h}) {$\d$};
    \fi}}
  \draw[black!35, line width=0.3pt] (0,0) grid[step=\h] ({\K*\h},{-\K*\h});
  \draw[line width=0.8pt] (0,0) rectangle ({\K*\h},{-\K*\h});
  \draw[nbhd one] ({(\cj-2)*\h},{-(\ci-2)*\h}) rectangle ++({3*\h},{-3*\h});
  \draw[nbhd two] ({(\cj-3)*\h},{-(\ci-3)*\h}) rectangle ++({5*\h},{-5*\h});
  \foreach \j in {1,...,\K} \node[axis tick] at ({(\j-0.5)*\h},0.26) {$\j$};
  \foreach \i in {1,...,\K} \node[axis tick] at (-0.26,{-(\i-0.5)*\h}) {$\i$};
  \node at ({\K*\h/2},0.60) {$a_{\beta,2}$};
  \node[rotate=90] at (-0.70,{-\K*\h/2}) {$a_{\beta,1}$};
  \begin{scope}[yshift={-\K*\h cm-0.46cm}]
    \draw[nbhd one] (0,0) -- (0.46,0);
    \node[anchor=west] at (0.52,0) {$r_\delta=1$};
    \draw[nbhd two] ({\K*\h/2+0.10},0) -- ({\K*\h/2+0.56},0);
    \node[anchor=west] at ({\K*\h/2+0.62},0) {$r_\delta=2$};
  \end{scope}
  \node[panel title] at ({\K*\h/2},0.98) {(a) $\|a_\beta-a_\alpha\|_\infty$};
\end{scope}

\begin{scope}[xshift=\gap cm]
  \def\s{\h}
  \foreach \i in {1,...,\K}{\foreach \j in {1,...,\K}{
    \foreach \di/\dj in {0/1,1/0,1/1,1/-1}{
      \pgfmathtruncatemacro{\ii}{\i+\di}\pgfmathtruncatemacro{\jj}{\j+\dj}
      \ifnum\ii<\Kp\ifnum\jj<\Kp\ifnum\jj>0
        \draw[black!20, line width=0.3pt] ({(\j-0.5)*\s},{-(\i-0.5)*\s}) -- ({(\jj-0.5)*\s},{-(\ii-0.5)*\s});
      \fi\fi\fi}}}
  \foreach \di in {-1,0,1}{\foreach \dj in {-1,0,1}{
    \draw[black, line width=1.3pt] ({(\cj-0.5)*\s},{-(\ci-0.5)*\s}) -- ({(\cj+\dj-0.5)*\s},{-(\ci+\di-0.5)*\s});}}
  \foreach \i in {1,...,\K}{\foreach \j in {1,...,\K}{
    \pgfmathtruncatemacro{\d}{max(abs(\i-\ci),abs(\j-\cj))}
    \ifcase\d \def\fc{dz}\or \def\fc{done}\else \def\fc{white}\fi
    \filldraw[fill=\fc, draw=black, line width=0.5pt] ({(\j-0.5)*\s},{-(\i-0.5)*\s}) circle (\rv cm);}}
  \filldraw[fill=dz, draw=black, line width=0.7pt] ({(\cj-0.5)*\s},{-(\ci-0.5)*\s}) circle (\ra cm);
  \node[text=white] at ({(\cj-0.5)*\s},{-(\ci-0.5)*\s}) {$\alpha$};
  \begin{scope}[yshift={-\K*\h cm-0.46cm}]
    \filldraw[fill=dz, draw=black, line width=0.5pt] (0.14,0) circle (\rl cm);
    \node[anchor=west] at (0.34,0) {$\alpha$};
    \filldraw[fill=done, draw=black, line width=0.5pt] (1.10,0) circle (\rl cm);
    \node[anchor=west] at (1.30,0) {neighbors of $\alpha$};
  \end{scope}
  \node[panel title] at ({\K*\h/2},0.98) {(b) $\mathfrak G_{\mathrm{blk}}$ for $1<\delta\le2$};
\end{scope}
\end{tikzpicture}
\caption{Block graph on a $7\times7$ subdomain grid, $a_\alpha=(4,4)$.
(a) Chebyshev distance to $a_\alpha$; solid and dashed frames are
$\mathcal N_\delta(a_\alpha)$ \eqref{eq:chebyshev-neighborhood} for $r_\delta=1,2$.
(b) For $1<\delta\le2$, $\mathfrak G_{\mathrm{blk}}$ is the king graph: $\alpha$ is
adjacent to the $8$ blocks with $\|a_\alpha-a_\beta\|_\infty\le r_\delta=1$.}
\label{fig:chebyshev-block-graph}
\end{figure}

As a first combinatorial application, we bound the number of
structurally nonzero blocks in the upper-triangular part (including the diagonal) of the Gramian by a straightforward counting argument on the graph.
\begin{proposition}[Gramian block count]
\label{thm:gramian-block-overlap-bound}
Every vertex of $\mathfrak G_{\mathrm{blk}}$ satisfies
\begin{align}\label{eq:block-graph-degree-bound}
\deg(\al)
:=\bigl|\{\be\neq\al:\{\al,\be\}\in\mathcal E_{\mathrm{blk}}\}\bigr|
\leq(2r_\delta+1)^d-1.
\end{align}
Let
\begin{align}\label{eq:block storage}
\operatorname{Sp}_{\mathrm{blk}}^+(G)
:=\{(\al,\be)\in\operatorname{Sp}_{\mathrm{blk}}(G):\al\leq\be\}
\end{align}
denote the upper-triangular block sparsity pattern, including the
diagonal. Then
\begin{align}\label{eq:upper-block-count-bound}
\bigl|\operatorname{Sp}_{\mathrm{blk}}^+(G)\bigr|
\leq K+|\mathcal E_{\mathrm{blk}}|
\leq\frac K2\bigl((2r_\delta+1)^d+1\bigr).
\end{align}
\end{proposition}
\begin{proof}
Fix $\alpha\in\{1,\ldots,K\}$. By the Chebyshev characterization of the block pattern \eqref{eq:chebyshev-cooccurrence-condition}, a structurally nonzero block $G_{\alpha\beta}$ can occur only if and only if $\|a_\alpha-a_\beta\|_\infty\leq r_\delta$. Thus, for each $j$, the index $a_{\beta,j}$  has at most $2r_\delta+1$ possible values. Since the subdomain grid is tensor-product, the $\alpha$-th block row then has at most $(2r_\delta+1)^d$ structurally nonzero block columns, including the diagonal block. Hence, the degree of the corresponding vertex in $\mathfrak G_{\mathrm{blk}}$ is at most $(2r_\delta+1)^d-1$.

Since $\mathfrak G_{\mathrm{blk}}$ is undirected, the degree sum counts each edge in the graph twice. Therefore, we have
\begin{align*}
2|\mathcal E_{\mathrm{blk}}|
=
\sum_{\alpha=1}^{K}\deg(\alpha)
\leq
K\left((2r_\delta+1)^d-1\right).
\end{align*}
Counting $K$ diagonal blocks in addition to the above estimate yields \eqref{eq:upper-block-count-bound}.
\end{proof}
\noindent
We note that since $G(\theta)$ is symmetric, there is no need to store its entire block representation. In fact, the numerical solvers used in this work operate directly on
the upper-triangular block part
$\operatorname{Sp}_{\mathrm{blk}}^+(G)$. If each stored block is represented
as a dense $n\times n$ matrix, Proposition~\ref{thm:gramian-block-overlap-bound} yields the following
bound on the number of scalar values required to store the upper-triangular
part $\operatorname{Sp}_{\mathrm{blk}}^+(G)$:
\begin{align}\label{eq:upper-gramian-scalar-count}
\nnz_{\mathrm{stored}}(G)
:=
n^2\bigl|\operatorname{Sp}_{\mathrm{blk}}^+(G)\bigr|
\leq
\frac{n^2K}{2}
\bigl((2r_\delta+1)^d+1\bigr).
\end{align}
However, this estimate is crude, as it treats every subdomain as if its
full Chebyshev neighborhood were present, whereas subdomains near the boundary
of the grid have fewer neighbors. Nevertheless, it exposes the structural dependencies. In particular, \eqref{eq:upper-gramian-scalar-count} is linear in $K$ but
quadratic in the number $n$ of parameters in one local network. Moreover, increasing the overlap enlarges $r_\delta$ and therefore increases
the number of structurally nonzero Gramian blocks. 

To illustrate the resulting scale, we report representative
values in Table~\ref{tab:one-two-layer-memory-comparison}.
\begin{table}[h!]
\centering
\small
\begin{tabular}{ccccccc}
\toprule
Dimension
& Partition
& $K$
& Overlap ratio
& Stored GN blocks
& 1-layer net
& 2-layer net \\
\midrule
$d=2$
& $20^2$
& $400$
& $1<\delta \leq 2$
&  $5K$
&  $0.27$ GB
&  $22.47$ GB
\\
$d=2$
& $20^2$
& $400$
& $2< \delta \leq 3$
& $13K$
& $0.69$ GB
& $58.42$ GB
\\
$d=2$
& $50^2$
& $2500$
& $1<\delta \leq 2$
& $5K$
& $1.66$ GB
& $140.42$ GB
\\
$d=2$
& $50^2$
& $2500$
& $2< \delta \leq 3$
& $13K$
& $4.33$ GB
& $365.10$ GB
\\
\midrule
$d=3$
& $20^3$
& $8000$
& $1<\delta \leq 2$
& $14K$
& $23.23$ GB
& $1.33$ TB
\\
$d=3$
& $20^3$
& $8000$
& $2 < \delta \leq 3$
& $63K$
& $104.51$ GB
& $5.97$ TB
\\
$d=3$
& $50^3$
& $125000$
& $1<\delta \leq 2$
& $14K$
& $362.89$ GB
& $20.74$ TB
\\
$d=3$
& $50^3$
& $125000$
& $2 < \delta \leq 3$
& $63K$
& $1.63$ TB
& $93.31$ TB
\\
\bottomrule
\end{tabular}
\caption{Estimated dense-block storage for the upper-triangular Gramian block
pattern in float64 precision. The one-hidden-layer network has width $32$ and
parameter count $n=32d+65$. The two-hidden-layer network has width $32$ in each
hidden layer and parameter count $n=32d+1121$.}
\label{tab:one-two-layer-memory-comparison}
\end{table}
The estimates illustrate that adding a second hidden layer increases $n$ sufficiently that the quadratic dependence on $n$ in \eqref{eq:upper-gramian-scalar-count} becomes dominant. From the viewpoint of Gramian storage, scaling through a larger number of small
local subnetworks is therefore substantially more favorable than increasing
the size of each subnetwork. In addition, larger overlap reduces
sparsity by increasing the number of neighboring blocks and can lead to a
substantial increase in memory, particularly in three dimensions.

Finally, we note that these estimates should be interpreted with care, since they account only for the numerical values stored in the assembled Gramian and do not represent the full memory footprint of the optimization procedure. In particular, AD may additionally require storing intermediate quantities and parts of the computational graph, whose
memory cost can substantially exceed that of the assembled matrix itself.
Thus, Table~\ref{tab:one-two-layer-memory-comparison} should primarily be
viewed as an indication of the scale problem.

\section{Sparse automatic differentiation for FBPINNs}\label{Sec: sparse AD}

The preceding section identified the sparsity pattern of the residual Jacobian
$J_{\mathrm R}(\theta)$ and the induced block sparsity of the Gauss--Newton
matrix. We now use this structure to reduce the number of JVPs and VJPs required for materialization.
Recall that for a parameter tangent $v\in\R^P$, forward-mode AD computes the
JVP
\begin{align}\label{eq:jvp-definition-revised}
\operatorname{JVP}_\theta(v)
:=
\diff\mathrm R(\theta)[v]
=
J_{\mathrm R}(\theta)v.
\end{align}
Choosing the unit seed vector $v=e_j$ gives the $j$-th column of
$J_{\mathrm R}(\theta)$, so dense column-wise materialization of
$J_{\mathrm R}(\theta)$ requires $P$ JVPs. Reverse-mode AD evaluates the adjoint action of the residual differential.
For a vector $w\in\R^{N_\Omega}$ in the residual space,
\begin{align}\label{eq:vjp-definition-revised}
\operatorname{VJP}_\theta(w)
:=
\diff\mathrm R(\theta)^*[w]
=
J_{\mathrm R}(\theta)^\top w
\in\R^P.
\end{align}
Choosing the unit seed vector $w=e_q$ gives the transpose of the $q$-th row of
$J_{\mathrm R}(\theta)$. Thus, dense row-wise materialization of
$J_{\mathrm R}(\theta)$ requires $N_\Omega$ VJPs. 
A single JVP or VJP is matrix-free and cost a small multiple of one residual evaluation, 
although reverse mode typically stores additional intermediate states and is therefore more memory intensive. We refer to
\cite{griewank2008evaluating} for further details on AD.

\subsection{Forward-mode sparse automatic differentiation for FBPINNs}

For FBPINNs, forward-mode sparse AD colors parameter columns of the Jacobian matrix,
so it is natural to start from the parameter blocks associated with the subdomains.
A block coloring is a map
\begin{align}
\kappa: \{1,\ldots,K\}\to \{1,\ldots,C_{\mathrm{blk}}\}
\end{align}
such that adjacent vertices in the block graph $\mathfrak G_{\mathrm{blk}}$ receive different colors: for any $\alpha\neq\beta$ with $Q_{\alpha\beta}\neq\emptyset$, one has
$\kappa(\alpha)\neq\kappa(\beta)$. Geometrically, it means that overlapping subdomains receive different block colors. As discussed above, for a rectangular computational domain
with a uniform tensor-product domain decomposition, the block graph
$\mathfrak G_{\mathrm{blk}}$ is a Chebyshev-neighborhood graph
on the subdomain grid. This structure admits an optimal periodic
coloring given by an explicit formula, as established next.
\begin{proposition}\label{prop:block-coloring}
Let $r_{\delta}=\lceil\delta\rceil-1$, where $\delta>1$ is the overlap
ratio, and define the block coloring period $p_{\delta}=r_{\delta}+1$.
Assume that $K_0\geq p_{\delta}$, where $K_0$ is the number of
subdomains in each coordinate direction. Using the lexicographic
indexing map $\operatorname{lex}$ defined in \eqref{lex ordering},
the map
\begin{align}\label{block coloring map}
\kappa: \{1,\ldots,K\}\to \{1,\ldots,p_{\delta}^{d}\},
\qquad
\kappa(\alpha)
:=
1+
\sum_{j=1}^{d}
\Bigl(
\bigl((\operatorname{lex}^{-1}(\alpha))_j-1\bigr)
\bmod p_{\delta}
\Bigr)
p_{\delta}^{d-j}
\end{align}
defines a valid coloring of $\mathfrak G_{\mathrm{blk}}$ with
$p_\delta^d$ colors. This number of colors is optimal, i.e.,
\begin{align}\label{Chrom number}
\chi(\mathfrak G_{\mathrm{blk}})=p_\delta^d,
\end{align}
where $\chi(\mathfrak G_{\mathrm{blk}})$ denotes the chromatic number
of $\mathfrak G_{\mathrm{blk}}$.
\end{proposition}
\begin{proof}
See~\ref{app:proof-block-coloring}.
\end{proof}

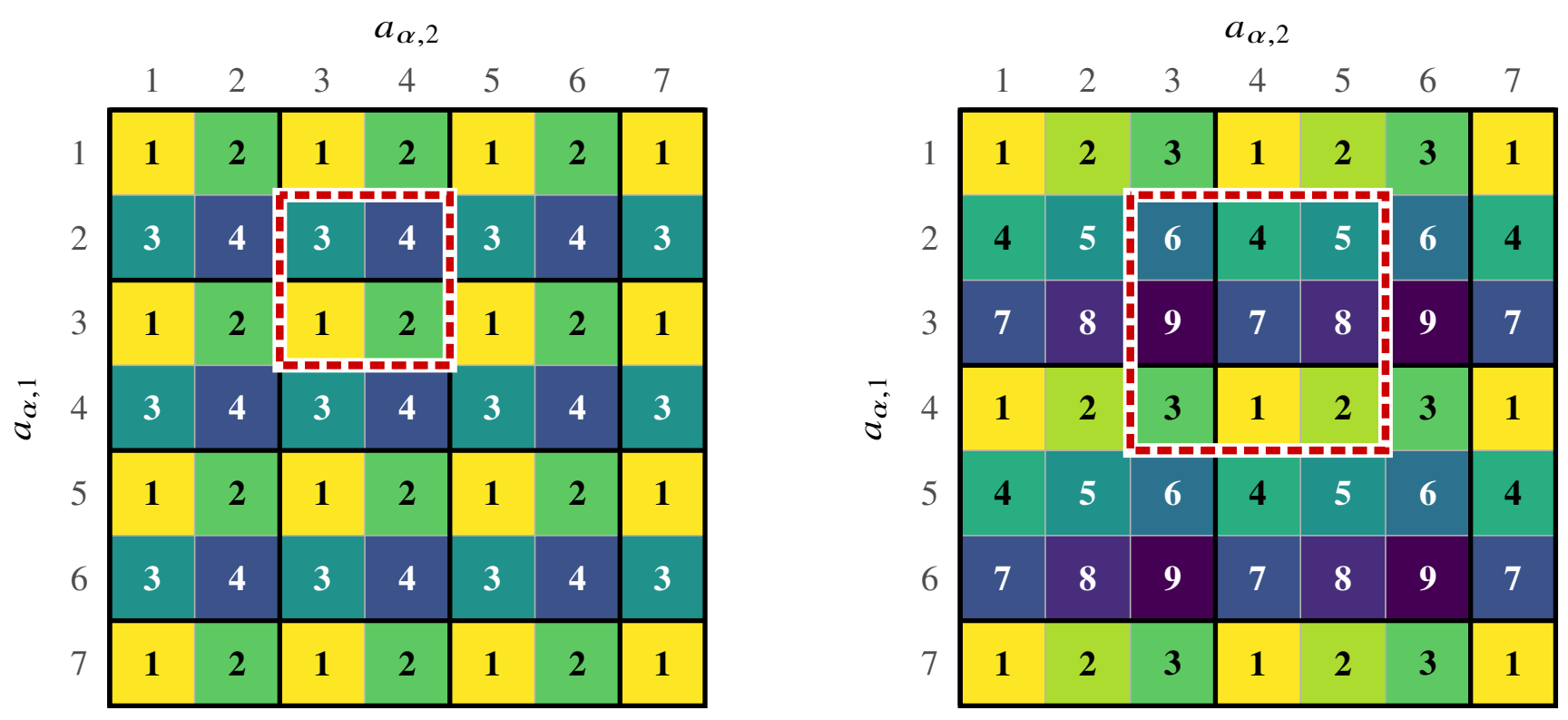
\begin{figure}[t]
\centering
\begin{tikzpicture}[font=\scriptsize]
\def\K{7}
\def\h{0.62}
\definecolor{va1}{HTML}{FDE725}\definecolor{va2}{HTML}{5EC962}
\definecolor{va3}{HTML}{21918C}\definecolor{va4}{HTML}{3B528B}
\definecolor{vb1}{HTML}{FDE725}\definecolor{vb2}{HTML}{ADDC30}\definecolor{vb3}{HTML}{5EC962}
\definecolor{vb4}{HTML}{28AE80}\definecolor{vb5}{HTML}{21918C}\definecolor{vb6}{HTML}{2C728E}
\definecolor{vb7}{HTML}{3B528B}\definecolor{vb8}{HTML}{472D7B}\definecolor{vb9}{HTML}{440154}

\foreach \p/\xo/\pal/\dark/\ttl in {2/0/va/2/{(a) $1<\delta\le 2$: $r_\delta=1,\ p_\delta=2$, $C_{\mathrm{blk}}=4$},
                        3/6.2/vb/4/{(b) $2<\delta\le 3$: $r_\delta=2,\ p_\delta=3$, $C_{\mathrm{blk}}=9$}}{
  \begin{scope}[xshift=\xo cm]
    \foreach \i in {1,...,\K}{
      \foreach \j in {1,...,\K}{
        \pgfmathtruncatemacro{\k}{1+\p*mod(\i-1,\p)+mod(\j-1,\p)}
        \fill[\pal\k] ({(\j-1)*\h},{-(\i-1)*\h}) rectangle ++(\h,-\h);
        \ifnum\k>\dark
          \node[text=white] at ({(\j-0.5)*\h},{-(\i-0.5)*\h}) {$\mathbf{\k}$};
        \else
          \node at ({(\j-0.5)*\h},{-(\i-0.5)*\h}) {$\mathbf{\k}$};
        \fi
      }}
    \draw[white!70!black, line width=0.3pt] (0,0) grid[step=\h] ({\K*\h},{-\K*\h});
    \foreach \t in {0,\p,...,\K}{
      \draw[line width=1.0pt] ({\t*\h},0) -- ({\t*\h},{-\K*\h});
      \draw[line width=1.0pt] (0,{-\t*\h}) -- ({\K*\h},{-\t*\h});}
    \draw[line width=1.0pt] (0,0) rectangle ({\K*\h},{-\K*\h});
    \draw[white, line width=3.0pt] ({2*\h},{-1*\h}) rectangle ++({\p*\h},{-\p*\h});
    \draw[red!80!black, line width=1.6pt, dash pattern=on 3pt off 1.5pt] ({2*\h},{-1*\h}) rectangle ++({\p*\h},{-\p*\h});
    \foreach \j in {1,...,\K} \node[text=black!70] at ({(\j-0.5)*\h},0.22) {$\j$};
    \foreach \i in {1,...,\K} \node[text=black!70] at (-0.22,{-(\i-0.5)*\h}) {$\i$};
    \node at ({\K*\h/2},0.55) {$a_{\alpha,2}$};
    \node[rotate=90] at (-0.6,{-\K*\h/2}) {$a_{\alpha,1}$};
    \node[font=\footnotesize, anchor=north] at ({\K*\h/2},{-\K*\h-0.18}) {\ttl};
  \end{scope}}
\end{tikzpicture}
\caption{Periodic block coloring \eqref{block coloring map} on a $7\times7$
subdomain grid; the number in each cell is the color $\kappa(\alpha)$ of the
subdomain with grid index $a_\alpha=(a_{\alpha,1},a_{\alpha,2})$. Any $p_\delta\times p_\delta$
window (dashed) contains every color exactly once and is a clique of
$\mathfrak G_{\mathrm{blk}}$.}
\label{fig:graph periodic coloring}
\end{figure}

\begin{Example}[Periodic coloring of two-dimensional block graphs]
For \(d=2\), the coloring map \eqref{block coloring map} takes the form
\begin{align}
\kappa(\alpha)
=
1
+
p_\delta
\Bigl((a_{\alpha,1}-1)\bmod p_\delta\Bigr)
+
\Bigl((a_{\alpha,2}-1)\bmod p_\delta\Bigr),
\end{align}
where
$a_\alpha=(a_{\alpha,1},a_{\alpha,2})
=\operatorname{lex}^{-1}(\alpha)\in\mathcal I_{K_0}$
is the grid index of the $\alpha$-th subdomain. If $1<\delta\leq2$, then $r_\delta=1$ and $p_\delta=2$, and
$\mathfrak{G}_{\mathrm{blk}}$ is the two-dimensional king graph with a
$2$-periodic coloring. If $2<\delta\leq3$, then $r_\delta=2$ and
$p_\delta=3$, and $\mathfrak{G}_{\mathrm{blk}}$ can be colored by a
$3$-periodic coloring. We refer to Fig.~\ref{fig:graph periodic coloring},
where these two coloring schemes are illustrated.
\end{Example}

The block coloring groups subdomains, but each subdomain contributes an entire column block of $n$ local parameter columns to $J_{\mathrm R}$. Hence, after a block color has been assigned, one still has to distinguish the $n$ columns inside each block. For this purpose, each block color $\kappa(\al)\in \{1,\ldots,C_{\mathrm{blk}}\}$ is expanded into a palette of $n$ scalar column colors:
\begin{align}\label{eq:scalar-palette-revised}
\mathfrak{P}_{\kappa(\al)}:=\{ n(\kappa(\al)-1)+1,\ldots,n\kappa(\al) \}.
\end{align}
Thus, the $\alpha$-th subdomain with block color $\kappa(\alpha)$ uses the palette $\mathfrak{P}_{\kappa(\alpha)}$. The $p$-th local parameter column in this subdomain is assigned the $p$-th color $\gamma(\al,p)\in \mathfrak{P}_{\kappa(\alpha)}$ in this palette via a simple index linearization:
\begin{align}\label{eq:scalar-column-color-revised}
\gamma(\al,p)
:=
n\bigl(\kappa(\al)-1\bigr)+p.
\end{align}
Subdomains with the same block color use the same palette: the equality of scalar colors $\gamma(\be,s)=\gamma(\al, p)$ implies $s=p$ and $\kappa(\be) = \kappa(\al)$ for $\be \neq \al$. This yields 
\begin{align}
C_{\mathrm{col}}=nC_{\mathrm{blk}}
=n(r_\delta+1)^d
\end{align}
column colors in total. It is favorable to choose a small subnetwork size $n$ and $\delta\leq2$ to keep this number small.

For a subdomain label $\alpha\in\{1,\ldots,K\}$ and a local parameter index $p\in\{1,\ldots,n\}$ in the Jacobian block $\Phi_\alpha$, we define the corresponding global scalar column index by
\begin{align}\label{eq:global-column-index-revised}
\iota(\alpha,p):=(\alpha-1)n+p.
\end{align}
Using the block coloring map \eqref{block coloring map} and the scalar column colors \eqref{eq:scalar-column-color-revised}, 
we construct the column seed matrix $D \in\{0,1\}^{P\times C_{\mathrm{col}}}$
by grouping Jacobian columns with the same scalar color. Its entries are defined by
\begin{align}\label{eq:forward-seed-matrix-revised}
D_{\iota(\alpha,p),\ell}
=
\begin{cases}
1, & \ell=\gamma(\alpha,p),\\
0, & \text{otherwise},
\end{cases}
\qquad
\ell=1,\ldots,C_{\mathrm{col}}.
\end{align}
Here, the $\ell$-th column defines a seed vector whose nonzero entries
correspond to the global parameter columns with scalar color $\ell$. Applying forward-mode AD to these seed directions gives the compressed Jacobian
\begin{align}\label{eq:foward-compressed-jac}
Y:=J_{\mathrm R}D \in \mathbb R^{N_\Omega\times C_{\mathrm{col}}}.
\end{align}
Computing $Y$ requires $C_{\mathrm{col}}$ JVPs. 
Since evaluating all of them
simultaneously can still lead to substantial peak memory usage, 
we allow the block colors to be processed in batches of at most $B_{\mathrm{fwd}}$ colors,
together with their complete palettes.
A batch of $B_{\mathrm{fwd}}$ block
colors therefore accounts for $nB_{\mathrm{fwd}}$ JVPs.
Thus, all block colors are covered in
$\lceil C_{\mathrm{blk}}/B_{\mathrm{fwd}}\rceil$ batches, while the total
number of JVPs remains $C_{\mathrm{col}}$.

\begin{figure}[t]
\centering
\begin{tikzpicture}[font=\scriptsize,>=Latex]
 
\definecolor{scalarA}{HTML}{005BFF} 
\definecolor{scalarC}{HTML}{009E2F} 
\definecolor{scalarB}{HTML}{E00000} 
\definecolor{scalarD}{HTML}{7A00CC} 
 
\tikzset{
  cell/.style   ={draw=black!65, minimum width=7.2mm, minimum height=7.0mm,
                  inner sep=0pt, text height=2.3mm, text depth=1.1mm},
  lab/.style    ={font=\scriptsize, text=black!85, inner sep=0pt},
  hdr/.style    ={font=\tiny, text=black!75, inner sep=0pt},
  shd/.style    ={font=\tiny, text=black!48, inner sep=0pt},
  ptitle/.style ={font=\scriptsize\bfseries, anchor=base},
  note/.style   ={font=\tiny, text=black!65, align=center, inner sep=0pt},
  grp/.style    ={line width=0.5pt, black!50},
}
\newcommand{\Jc}[2]{{\color{#1}\scriptstyle J_{#2}}}
\newcommand{\Yc}[2]{\raisebox{1.1mm}{$\color{#1}\scriptstyle J_{#2}$}}
 
\matrix (J) [matrix of math nodes, nodes={cell}, nodes in empty cells,
             row sep=-\pgflinewidth, column sep=-\pgflinewidth] at (0,0) {
  \Jc{scalarA}{11} & \Jc{scalarC}{12} & \Jc{scalarB}{13} & \Jc{scalarD}{14} &
                   &                  &                  &                  \\
                   &                  & \Jc{scalarB}{23} & \Jc{scalarD}{24} &
                   &                  &                  &                  \\
                   &                  &                  &                  &
  \Jc{scalarA}{35} & \Jc{scalarC}{36} & \Jc{scalarB}{37} & \Jc{scalarD}{38} \\
                   &                  &                  &                  &
                   &                  & \Jc{scalarB}{47} & \Jc{scalarD}{48} \\
};
\begin{scope}[on background layer]
  \foreach \ca/\cb/\c in {1/2/lightA, 3/4/lightB, 5/6/lightA, 7/8/lightB}
    \fill[\c, opacity=0.55] (J-1-\ca.north west) rectangle (J-4-\cb.south east);
\end{scope}
 
\foreach \ca/\cb/\al/\kap in {1/2/1/1, 3/4/2/2, 5/6/3/1, 7/8/4/2}{
  \draw[grp] ([xshift=0.6mm,yshift=4.3mm]J-1-\ca.north west)
          -- ([xshift=-0.6mm,yshift=4.3mm]J-1-\cb.north east);
  \node[hdr, anchor=south, inner sep=0.5mm]
       at ([yshift=4.3mm]$(J-1-\ca.north west)!0.5!(J-1-\cb.north east)$)
       {$\Phi_{\al}$, $\kappa{=}\kap$};}
\foreach \c/\l/\col in {1/1/scalarA, 2/2/scalarC, 3/3/scalarB, 4/4/scalarD,
                        5/1/scalarA, 6/2/scalarC, 7/3/scalarB, 8/4/scalarD}
  \node[hdr, text=\col, anchor=south] at ([yshift=0.7mm]J-1-\c.north) {$\ell{=}\l$};
\foreach \r in {1,2,3,4}
  \node[lab, anchor=east] at ([xshift=-1.3mm]J-\r-1.west) {$q_{\r}$};
 
\matrix (D) [matrix of math nodes, nodes={cell},
             row sep=-\pgflinewidth, column sep=-\pgflinewidth] at (6.24,0) {
  1 & 0 & 0 & 0 \\
  0 & 1 & 0 & 0 \\
  0 & 0 & 1 & 0 \\
  0 & 0 & 0 & 1 \\
  1 & 0 & 0 & 0 \\
  0 & 1 & 0 & 0 \\
  0 & 0 & 1 & 0 \\
  0 & 0 & 0 & 1 \\
};
\begin{scope}[on background layer]
  \foreach \ra/\rb/\c in {1/2/lightA, 3/4/lightB, 5/6/lightA, 7/8/lightB}
    \fill[\c, opacity=0.55] (D-\ra-1.north west) rectangle (D-\rb-4.south east);
\end{scope}
\foreach \c/\col in {1/scalarA, 2/scalarC, 3/scalarB, 4/scalarD}
  \node[hdr, text=\col, anchor=south] at ([yshift=0.7mm]D-1-\c.north) {$\ell{=}\c$};
\foreach \r/\a/\p in {1/1/1, 2/1/2, 3/2/1, 4/2/2, 5/3/1, 6/3/2, 7/4/1, 8/4/2}
  \node[shd, anchor=east] at ([xshift=-1.0mm]D-\r-1.west) {$(\a,\p)$};
\node[note, anchor=north] at ([yshift=-1.4mm]D.south)
     {one column per scalar color};
 
\matrix (Y) [matrix of nodes, nodes={cell},
             row sep=-\pgflinewidth, column sep=-\pgflinewidth] at (10.84,0) {
  {\Yc{scalarA}{11}} & {\Yc{scalarC}{12}} & {\Yc{scalarB}{13}} & {\Yc{scalarD}{14}} \\
  {}                 & {}                 & {\Yc{scalarB}{23}} & {\Yc{scalarD}{24}} \\
  {\Yc{scalarA}{35}} & {\Yc{scalarC}{36}} & {\Yc{scalarB}{37}} & {\Yc{scalarD}{38}} \\
  {}                 & {}                 & {\Yc{scalarB}{47}} & {\Yc{scalarD}{48}} \\
};
\begin{scope}[on background layer]
  \foreach \c/\col in {1/scalarA,2/scalarC,3/scalarB,4/scalarD}
    \fill[\col!10] (Y-1-\c.north west) rectangle (Y-4-\c.south east);
\end{scope}
\foreach \r/\c/\a/\p in {1/1/1/1, 1/2/1/2, 1/3/2/1, 1/4/2/2,
                                           2/3/2/1, 2/4/2/2,
                         3/1/3/1, 3/2/3/2, 3/3/4/1, 3/4/4/2,
                                           4/3/4/1, 4/4/4/2}
  \node[shd, anchor=south] at ([yshift=0.6mm]Y-\r-\c.south) {$(\a,\p)$};
\foreach \c/\col in {1/scalarA,2/scalarC,3/scalarB,4/scalarD}
  \node[hdr, text=\col, anchor=south] at ([yshift=0.7mm]Y-1-\c.north) {$\ell{=}\c$};
\foreach \ca/\cb/\pal/\kap in {1/2/1/1, 3/4/2/2}{
  \draw[grp] ([xshift=0.6mm,yshift=4.3mm]Y-1-\ca.north west)
          -- ([xshift=-0.6mm,yshift=4.3mm]Y-1-\cb.north east);
  \node[hdr, anchor=south, inner sep=0.5mm]
       at ([yshift=4.3mm]$(Y-1-\ca.north west)!0.5!(Y-1-\cb.north east)$)
       {$\mathfrak P_{\pal}$, $\kappa{=}\kap$};}
\foreach \r in {1,2,3,4}
  \node[lab, anchor=east] at ([xshift=-1.3mm]Y-\r-1.west) {$q_{\r}$};
\node[note, anchor=north] at ([yshift=-1.4mm]Y.south)
     {$C_{\mathrm{col}}=4$ JVPs};
 
\node[text=black!70, font=\small] at (3.45,0) {$\times$};
\node[text=black!70, font=\small] at (8.25,0) {$=$};
 
\node[ptitle] at (0,2.45)     {Block-column Jacobian $J_{\mathrm R}$};
\node[ptitle] at (6.24,3.45)  {Seed matrix $D$};
\node[ptitle] at (10.84,2.45) {Compressed Jacobian $Y$};
 
\end{tikzpicture}
\caption{Forward-mode column compression for
Example~\ref{Forward mode example}. The seed
matrix $D$ has one column per scalar color and gathers the structurally
orthogonal columns carrying that color.}
\label{fig:forward-compression-tikz}
\end{figure}
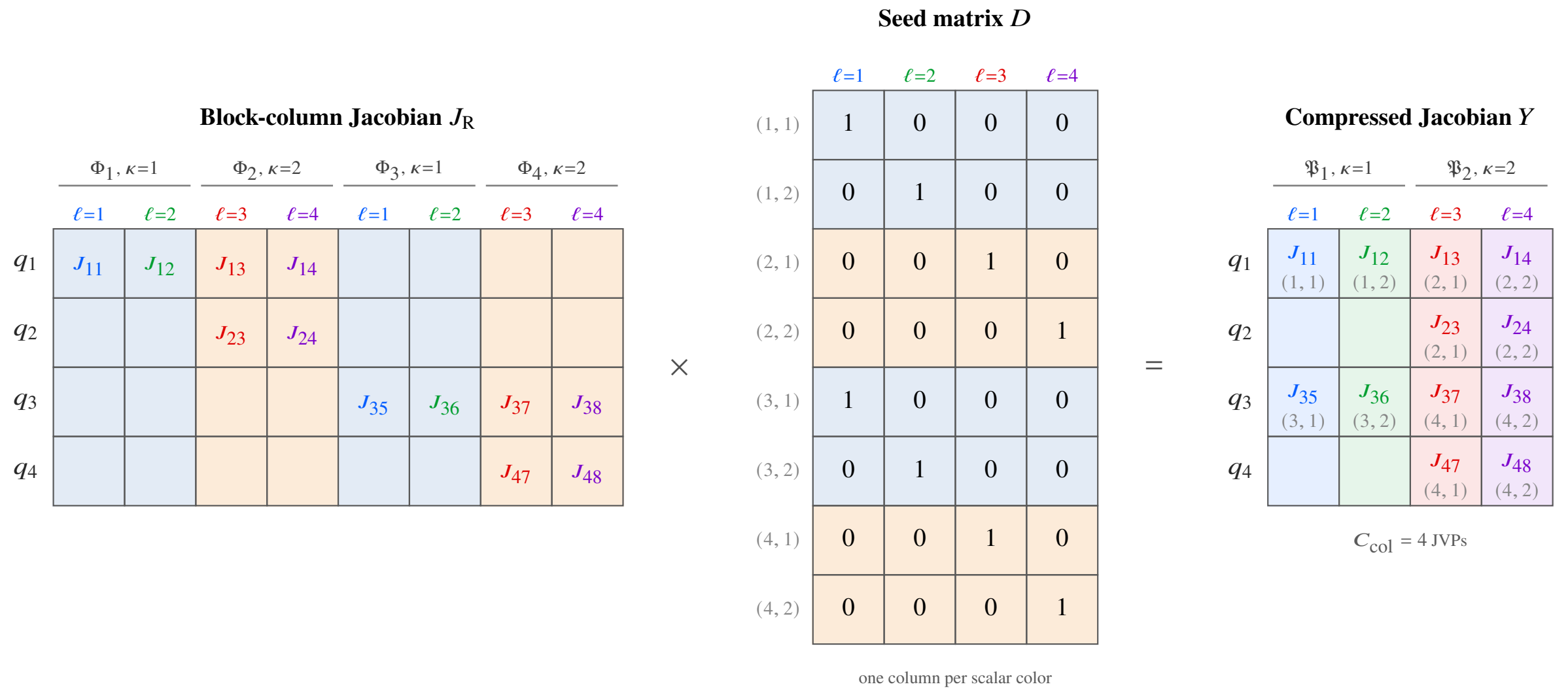

\begin{Example}[Forward-mode Jacobian compression]\label{Forward mode example}
We use this and the following examples to illustrate the sparse AD constructions
on small toy problems. Consider four subdomains, each with \(n=2\) local parameters. Thus, the Jacobian has the block-column form
\begin{align}
J_{\mathrm R}
=
\begin{bmatrix}
\Phi_1 & \Phi_2 & \Phi_3 & \Phi_4
\end{bmatrix},
\qquad
\Phi_\alpha\in\mathbb R^{N_\Omega\times 2}.
\end{align}
Suppose that the block coloring \eqref{block coloring map} assigns $\kappa(1)=\kappa(3)=1$ and $\kappa(2)=\kappa(4)=2$. Thus, subdomains $1$ and $3$ may be grouped together, and subdomains $2$ and $4$ may be grouped together as well. Since $n=2$, each block color is expanded into a palette of two scalar colors:
\begin{align}
\mathfrak P_1
= \{1,2\},
\qquad
\mathfrak P_2 = \{3,4\}.
\end{align}
Therefore, the scalar column colors are $\gamma(1,1)=\gamma(3,1)=1$, $\gamma(1,2)=\gamma(3,2)=2$, and $\gamma(2,1)=\gamma(4,1)=3$, $\gamma(2,2)=\gamma(4,2)=4$. The seed matrix $D$ has one column for each scalar color. 
Hence, computing $Y=J_{\mathrm R}D$ requires 4 JVP directions, one for each scalar color. 
We refer to Fig.~\ref{fig:forward-compression-tikz} for
a visual illustration. 
\end{Example}

Since the seed matrix \eqref{eq:forward-seed-matrix-revised} consists of columns which are structurally orthogonal, JVP contributions in the compressed Jacobian do not collide. By the sparsity pattern characterization \eqref{eq: structural equivalence of active sets},
$T_q$ identifies the subdomain Jacobian blocks that are structurally present in
row $q$. Hence, for every $\alpha\in T_q$ and $p=1,\ldots,n$,
\begin{align}\label{eq:forward-decompression}
(\Phi_\alpha)_{q,p}
=
Y_{q,\gamma(\alpha,p)}.
\end{align}
Thus, each membership $(q,\alpha)$ with $\alpha\in T_q$ determines a local row
$(\Phi_\alpha)_{q,:}$ recovered from $Y$;
we refer to this as \textit{decompression}. 
Using \eqref{eq:gramian-block-Q-note} and \eqref{eq:forward-decompression}, each nonzero Gramian block can then be assembled as follows:
\begin{align}\label{eq:forward-gramian-from-compressed}
\bigl(G_{\al\be}\bigr)_{ps}
=
\sum_{q\in Q_{\al\be}}
Y_{q,\gamma(\al,p)}\,
Y_{q,\gamma(\be,s)},
\qquad
p,s=1,\ldots,n.
\end{align}
%
%
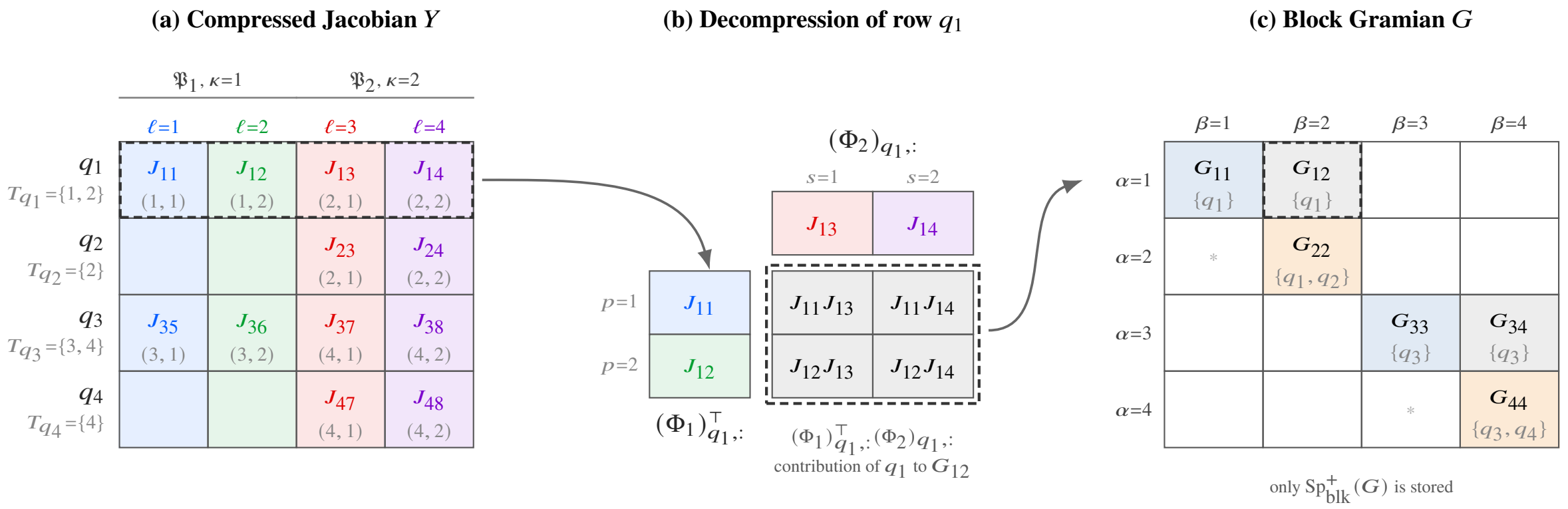
\begin{figure}[H]
\centering
\begin{tikzpicture}[font=\scriptsize,>=Latex]

\definecolor{scalarA}{HTML}{005BFF} 
\definecolor{scalarC}{HTML}{009E2F} 
\definecolor{scalarB}{HTML}{E00000} 
\definecolor{scalarD}{HTML}{7A00CC} 

\tikzset{
  ycell/.style  ={draw=black!65, minimum width=8.8mm, minimum height=7.6mm,
                  inner sep=0pt, text height=2.3mm, text depth=1.1mm},
  ocell/.style  ={draw=black!65, minimum width=10.0mm, minimum height=6.3mm,
                  inner sep=0pt, text height=2.3mm, text depth=1.1mm},
  gcell/.style  ={draw=black!65, minimum width=9.8mm, minimum height=7.6mm,
                  inner sep=0pt, text height=2.3mm, text depth=1.1mm},
  lab/.style    ={font=\scriptsize, text=black!85, inner sep=0pt},
  hdr/.style    ={font=\tiny, text=black!75, inner sep=0pt},
  shd/.style    ={font=\tiny, text=black!48, inner sep=0pt},
  ptitle/.style ={font=\scriptsize\bfseries, anchor=base},
  note/.style   ={font=\tiny, text=black!65, align=center, inner sep=0pt},
  arr/.style    ={-{Latex[length=2.0mm]}, line width=0.7pt, black!60},
  grp/.style    ={line width=0.5pt, black!50},
  acc/.style    ={draw=black!80, line width=0.8pt, dash pattern=on 2.4pt off 1.4pt},
}
\newcommand{\Yc}[2]{\raisebox{1.2mm}{$\color{#1}\scriptstyle J_{#2}$}}
\newcommand{\Gc}[1]{\raisebox{1.2mm}{$\scriptstyle #1$}}

\matrix (Y) [matrix of nodes, nodes={ycell},
             row sep=-\pgflinewidth, column sep=-\pgflinewidth] at (0,0) {
  {\Yc{scalarA}{11}} & {\Yc{scalarC}{12}} & {\Yc{scalarB}{13}} & {\Yc{scalarD}{14}} \\
  {}                 & {}                 & {\Yc{scalarB}{23}} & {\Yc{scalarD}{24}} \\
  {\Yc{scalarA}{35}} & {\Yc{scalarC}{36}} & {\Yc{scalarB}{37}} & {\Yc{scalarD}{38}} \\
  {}                 & {}                 & {\Yc{scalarB}{47}} & {\Yc{scalarD}{48}} \\
};
\foreach \r/\c/\a/\p in {1/1/1/1, 1/2/1/2, 1/3/2/1, 1/4/2/2,
                                           2/3/2/1, 2/4/2/2,
                         3/1/3/1, 3/2/3/2, 3/3/4/1, 3/4/4/2,
                                           4/3/4/1, 4/4/4/2}
  \node[shd, anchor=south] at ([yshift=0.65mm]Y-\r-\c.south) {$(\a,\p)$};

\begin{scope}[on background layer]
  \foreach \c/\col in {1/scalarA,2/scalarC,3/scalarB,4/scalarD}
    \fill[\col!10] (Y-1-\c.north west) rectangle (Y-4-\c.south east);
\end{scope}

\foreach \c/\col in {1/scalarA,2/scalarC,3/scalarB,4/scalarD}
  \node[hdr, text=\col, anchor=south] at ([yshift=0.7mm]Y-1-\c.north) {$\ell{=}\c$};
\foreach \ca/\cb/\pal/\kap in {1/2/1/1, 3/4/2/2}{
  \draw[grp] ([yshift=4.3mm]Y-1-\ca.north west) -- ([yshift=4.3mm]Y-1-\cb.north east);
  \node[hdr, anchor=south, inner sep=0.5mm]
       at ([yshift=4.3mm]$(Y-1-\ca.north west)!0.5!(Y-1-\cb.north east)$)
       {$\mathfrak P_{\pal}$, $\kappa{=}\kap$};}

\foreach \r/\T in {1/{\{1,2\}},2/{\{2\}},3/{\{3,4\}},4/{\{4\}}}{
  \node[lab, anchor=east] at ([xshift=-1.3mm,yshift=1.3mm]Y-\r-1.west) {$q_{\r}$};
  \node[shd, anchor=east] at ([xshift=-1.3mm,yshift=-1.7mm]Y-\r-1.west) {$T_{q_\r}{=}\T$};}
\node[acc, fit=(Y-1-1)(Y-1-4), inner sep=-0.15mm] (qrow) {};

\matrix (O) [matrix of nodes, nodes={ocell},
             row sep=-\pgflinewidth, column sep=-\pgflinewidth] at (5.73,-0.395) {
  {$\scriptstyle J_{11}J_{13}$} & {$\scriptstyle J_{11}J_{14}$} \\
  {$\scriptstyle J_{12}J_{13}$} & {$\scriptstyle J_{12}J_{14}$} \\
};
\begin{scope}[on background layer]
  \fill[black!7] (O-1-1.north west) rectangle (O-2-2.south east);
\end{scope}
\node[acc, fit=(O-1-1)(O-2-2), inner sep=0.5mm] (obox) {};

\node[ocell, fill=scalarB!10] (v1) at ([yshift=7.9mm]O-1-1.center) {$\color{scalarB}\scriptstyle J_{13}$};
\node[ocell, fill=scalarD!10] (v2) at ([yshift=7.9mm]O-1-2.center) {$\color{scalarD}\scriptstyle J_{14}$};
\node[ocell, fill=scalarA!10] (u1) at ([xshift=-12.2mm]O-1-1.center) {$\color{scalarA}\scriptstyle J_{11}$};
\node[ocell, fill=scalarC!10] (u2) at ([xshift=-12.2mm]O-2-1.center) {$\color{scalarC}\scriptstyle J_{12}$};

\node[lab, anchor=south] at ([yshift=3.2mm]$(v1.north)!0.5!(v2.north)$)
     {$(\Phi_2)_{q_1,:}$};
\node[shd, anchor=south] at ([yshift=0.65mm]v1.north) {$s{=}1$};
\node[shd, anchor=south] at ([yshift=0.65mm]v2.north) {$s{=}2$};
\node[shd, anchor=east]  at ([xshift=-1.0mm]u1.west) {$p{=}1$};
\node[shd, anchor=east]  at ([xshift=-1.0mm]u2.west) {$p{=}2$};
\node[lab, anchor=north] at ([yshift=-1.3mm]u2.south) {$(\Phi_1)_{q_1,:}^{\top}$};
\node[note, anchor=north] at ([yshift=-1.3mm]O.south)
     {$(\Phi_1)_{q_1,:}^{\top}(\Phi_2)_{q_1,:}$\\contribution of $q_1$ to $G_{12}$};

\matrix (G) [matrix of nodes, nodes={gcell},
             row sep=-\pgflinewidth, column sep=-\pgflinewidth] at (10.59,0) {
  {\Gc{G_{11}}}                         & {\Gc{G_{12}}} & {}                                    & {}            \\
  {$\color{black!32}\scriptstyle \ast$} & {\Gc{G_{22}}} & {}                                    & {}            \\
  {}                                    & {}            & {\Gc{G_{33}}}                         & {\Gc{G_{34}}} \\
  {}                                    & {}            & {$\color{black!32}\scriptstyle \ast$} & {\Gc{G_{44}}} \\
};
\foreach \r/\c/\Q in {1/1/{\{q_1\}}, 1/2/{\{q_1\}}, 2/2/{\{q_1,q_2\}},
                      3/3/{\{q_3\}}, 3/4/{\{q_3\}}, 4/4/{\{q_3,q_4\}}}
  \node[shd, anchor=south] at ([yshift=0.65mm]G-\r-\c.south) {$\Q$};

\begin{scope}[on background layer]
  \foreach \r/\c in {1/lightA,2/lightB,3/lightA,4/lightB}
    \fill[\c, opacity=0.6] (G-\r-\r.north west) rectangle (G-\r-\r.south east);
  \fill[black!7] (G-1-2.north west) rectangle (G-1-2.south east);
  \fill[black!7] (G-3-4.north west) rectangle (G-3-4.south east);
\end{scope}
\node[acc, fit=(G-1-2), inner sep=-0.15mm] (gbox) {};

\foreach \c in {1,2,3,4}
  \node[hdr, anchor=south] at ([yshift=0.7mm]G-1-\c.north) {$\beta{=}\c$};
\foreach \r in {1,2,3,4}
  \node[hdr, anchor=east]  at ([xshift=-1.1mm]G-\r-1.west) {$\alpha{=}\r$};
\node[note, anchor=north] at ([yshift=-1.4mm]G.south)
     {only $\operatorname{Sp}^{+}_{\mathrm{blk}}(G)$ is stored};

\draw[arr] ([xshift=0.8mm]qrow.east) to[out=0,in=105] ([xshift=1.0mm]u1.north);

\draw[arr] ([xshift=0.8mm,yshift=0.4mm]obox.east) to[out=0,in=200] ([xshift=-8.0mm]G-1-1.west);

\coordinate (tline) at (0,2.66);
\node[ptitle] at (Y.north |- tline) {(a) Compressed Jacobian $Y$};
\node[ptitle] at (5.145,2.66)       {(b) Decompression of row $q_1$};
\node[ptitle] at (G.north |- tline) {(c) Block Gramian $G$};

\end{tikzpicture}
\caption{Assembly of the Gauss--Newton Gramian for
Example~\ref{Forward mode example}.
(a) Compressed Jacobian $Y=J_{\mathrm R}\Dmat$, with gray pairs
identifying the corresponding original Jacobian columns $(\alpha,p)$.
(b) Decompression of row $q_1$ and its contribution to $G_{12}$.
(c) Blockwise assembly over $Q_{\alpha\beta}$, indicated
in gray. Only upper-triangular blocks are stored, with
$\ast$ denoting their symmetric counterparts.}
\label{fig:forward-assembly-tikz}
\end{figure}

The same compressed data also gives the loss gradient by
\begin{align}\label{eq:forward-jtr-clean}
(J_{\mathrm R}^\top \mathrm{R})_{\iota(\alpha,p)}
= \sum_{q\in Q_\alpha} Y_{q,\gamma(\alpha,p)}\,\mathrm{R}_q .
\end{align}
With all the ingredients above, the Gauss-Newton system \eqref{eq:GN_system} can now be assembled.

\subsection{Reverse-mode sparse automatic differentiation for FBPINNs}
Reverse-mode sparse AD groups structurally orthogonal rows of the Jacobian such that they can be evaluated with a single VJP \cite{curtis1974estimation}.
We define a row-coloring map
\begin{align}\label{Sec 3: coloring map}
\chi : \{1, \dots , N_{\Omega} \} \rightarrow \{1, \dots, C_{\mathrm{row}} \}
\end{align}
such that, for any distinct $i,j\in\{1,\ldots,N_{\Omega}\}$, $\chi(i) = \chi(j)$ implies that $T_{i} \cap T_{j} = \emptyset$. The corresponding color classes are
\begin{align}\label{eq:row-color-classes}
\mathcal C_c^\chi
:=
\{\, i\in\{1,\ldots,N_{\Omega}\} : \chi(i)=c \,\},
\qquad
c\in\{1,\ldots,C_{\mathrm{row}}\}.
\end{align}
By the characterization \eqref{eq: structural equivalence of active sets}, we obtain by construction that rows
assigned to the same color do not share structurally nonzero columns in the
Jacobian matrix. 
The reverse-mode seed matrix $S\in\{0,1\}^{N_\Omega\times C_{\mathrm{row}}}$ is defined by
\begin{align}\label{eq:reverse-seed-matrix-clean}
S_{i,c}
=
\begin{cases}
1,&i\in\mathcal C_c^{\chi},\\
0,&\text{otherwise.}
\end{cases}
\end{align}
Stacking the corresponding $C_{\mathrm{row}}$ VJP evaluations gives the compressed Jacobian
\begin{align}\label{eq:reverse-compressed-jacobian-clean}
Z
:=
S^\top J_{\mathrm R}(\theta)
\in
\R^{C_{\mathrm{row}}\times P}.
\end{align}
To reduce peak memory usage, which is particularly important for the
more memory-intensive reverse mode, we may partition the row colors into
$N_{\mathrm{bat}}=\lceil C_{\mathrm{row}}/B_{\mathrm{rev}}\rceil$
disjoint batches $\mathcal I_b$ of at most $B_{\mathrm{rev}}$ colors.
For each batch, we compute
\begin{align}\label{eq:reverse-color-batching}
Z^{(b)}:=(S_{:,\mathcal I_b})^\top J_{\mathrm R}
\in\R^{|\mathcal I_b|\times P},
\end{align}
and assemble \eqref{eq:reverse-compressed-jacobian-clean} batch by batch.
\begin{Example}[Reverse-mode Jacobian compression]\label{ex:reverse-row-compression}
For this example, we consider the same setting as in Example \ref{Forward mode example}.  The four residual rows have the following active subdomain sets
\begin{align*}
T_1&=\{1,2\},
&
T_2&=\{2\},
&
T_3&=\{3,4\},
&
T_4&=\{4\}.
\end{align*}
Rows conflict if their active subdomain sets intersect. Thus, rows $1$ and $2$
conflict because $T_1\cap T_2=\{2\}$, while rows $1$ and $3$ do not conflict
because $T_1\cap T_3=\emptyset$. A valid row coloring is then
\begin{align*}
\chi(1)=\chi(3)=1,
\qquad
\chi(2)=\chi(4)=2.
\end{align*}
Hence, the row color classes are $\mathcal C_1^\chi=\{1,3\}$ and $\mathcal C_2^\chi=\{2,4\}$. Each color class contains structurally orthogonal rows. The transpose of the
seed matrix is given by
\begin{align}
S^{\top}
=
\begin{bmatrix}
1 & 0 & 1 & 0 \\
0 & 1 & 0 & 1
\end{bmatrix}.
\end{align}
Since the rows in each color class have disjoint column supports, no nonzero entries collide in the compressed Jacobian $Z$. Observe that only two row colors are required in this setting, while four column colors are needed for the sparse forward mode in Example~\ref{Forward mode example}.
\end{Example}

%
\begin{figure}[t]
\centering
\begin{tikzpicture}[font=\scriptsize,>=Latex]

\colorlet{rowA}{blue!10}    
\colorlet{rowB}{orange!16}  

\tikzset{
  cell/.style   ={draw=black!65, minimum width=7.2mm, minimum height=7.0mm,
                  inner sep=0pt, text height=2.3mm, text depth=1.1mm},
  lab/.style    ={font=\scriptsize, text=black!85, inner sep=0pt},
  hdr/.style    ={font=\tiny, text=black!75, inner sep=0pt},
  shd/.style    ={font=\tiny, text=black!48, inner sep=0pt},
  ptitle/.style ={font=\scriptsize\bfseries, anchor=base},
  note/.style   ={font=\tiny, text=black!65, align=center, inner sep=0pt},
  grp/.style    ={line width=0.5pt, black!50},
}
\newcommand{\Jc}[1]{{\scriptstyle J_{#1}}}

\matrix (S) [matrix of math nodes, nodes={cell}, nodes in empty cells,
             row sep=-\pgflinewidth, column sep=-\pgflinewidth] at (0,0) {
  1 & 0 & 1 & 0 \\
  0 & 1 & 0 & 1 \\
};
\begin{scope}[on background layer]
  \fill[rowA] (S-1-1.north west) rectangle (S-1-4.south east);
  \fill[rowB] (S-2-1.north west) rectangle (S-2-4.south east);
\end{scope}
\foreach \c in {1,2,3,4}
  \node[hdr, anchor=south] at ([yshift=0.7mm]S-1-\c.north) {$q_{\c}$};
\foreach \r/\C in {1/{\{q_1,q_3\}}, 2/{\{q_2,q_4\}}}{
  \node[lab, anchor=east] at ([xshift=-1.3mm,yshift=1.3mm]S-\r-1.west) {$c_{\r}$};
  \node[shd, anchor=east] at ([xshift=-1.3mm,yshift=-1.7mm]S-\r-1.west) {$\C$};}

\matrix (J) [matrix of math nodes, nodes={cell}, nodes in empty cells,
             row sep=-\pgflinewidth, column sep=-\pgflinewidth] at (7.08,0) {
  \Jc{11} & \Jc{12} & \Jc{13} & \Jc{14} &         &         &         &         \\
          &         & \Jc{23} & \Jc{24} &         &         &         &         \\
          &         &         &         & \Jc{35} & \Jc{36} & \Jc{37} & \Jc{38} \\
          &         &         &         &         &         & \Jc{47} & \Jc{48} \\
};
\begin{scope}[on background layer]
  \foreach \r/\c in {1/rowA, 2/rowB, 3/rowA, 4/rowB}
    \fill[\c] (J-\r-1.north west) rectangle (J-\r-8.south east);
\end{scope}
\foreach \ca/\cb/\al in {1/2/1, 3/4/2, 5/6/3, 7/8/4}{
  \draw[grp] ([xshift=0.6mm,yshift=4.3mm]J-1-\ca.north west)
          -- ([xshift=-0.6mm,yshift=4.3mm]J-1-\cb.north east);
  \node[hdr, anchor=south, inner sep=0.5mm]
       at ([yshift=4.3mm]$(J-1-\ca.north west)!0.5!(J-1-\cb.north east)$)
       {$\Phi_{\al}$};}
\foreach \r/\T in {1/{\{1,2\}}, 2/{\{2\}}, 3/{\{3,4\}}, 4/{\{4\}}}{
  \node[lab, anchor=east] at ([xshift=-1.3mm,yshift=1.3mm]J-\r-1.west) {$q_{\r}$};
  \node[shd, anchor=east] at ([xshift=-1.3mm,yshift=-1.7mm]J-\r-1.west) {$T_{q_\r}{=}\T$};}

\matrix (Z) [matrix of math nodes, nodes={cell}, nodes in empty cells,
             row sep=-\pgflinewidth, column sep=-\pgflinewidth] at (7.08,-3.26) {
  \Jc{11} & \Jc{12} & \Jc{13} & \Jc{14} & \Jc{35} & \Jc{36} & \Jc{37} & \Jc{38} \\
          &         & \Jc{23} & \Jc{24} &         &         & \Jc{47} & \Jc{48} \\
};
\begin{scope}[on background layer]
  \fill[rowA] (Z-1-1.north west) rectangle (Z-1-8.south east);
  \fill[rowB] (Z-2-1.north west) rectangle (Z-2-8.south east);
\end{scope}
\foreach \r in {1,2}
  \node[lab, anchor=east] at ([xshift=-1.3mm]Z-\r-1.west) {$c_{\r}$};
\node[note, anchor=north] at ([yshift=-1.4mm]Z.south)
     {$C_{\mathrm{row}}=2$ VJPs};

\node[text=black!70, font=\small] at (2.10,0) {$\times$};

\node[ptitle] at (0,1.35)     {Row seed matrix $S^{\top}$};
\node[ptitle] at (7.08,2.45)  {Block Jacobian $J_{\mathrm R}$};
\node[ptitle] at (7.08,-2.16) {Compressed Jacobian $Z=S^{\top}J_{\mathrm R}$};

\end{tikzpicture}
\caption{Reverse-mode compression in Example~\ref{ex:reverse-row-compression}.
Rows of one color have disjoint active sets $T_q$, so no two entries of
$Z=S^{\top}J_{\mathrm R}$ collide and the block row $(\Phi_\alpha)_{q,:}$ is the
slice $Z_{\chi(q),\mathcal C_\alpha}$, at a cost of $C_{\mathrm{row}}=2$ VJPs.}
\label{fig:reverse-compression-tikz-clean}
\end{figure}
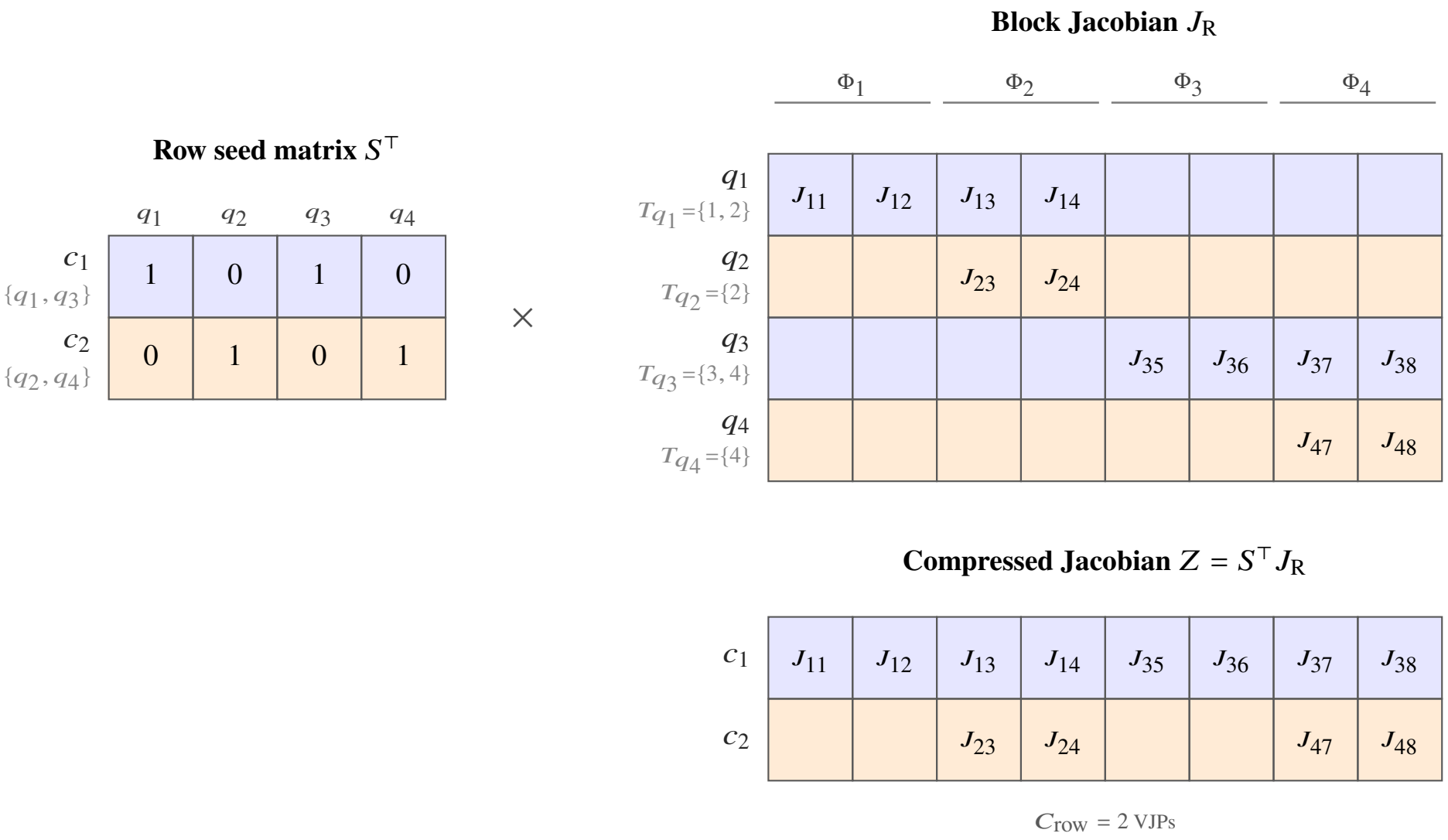

The compressed Jacobian matrix $Z=S^\top J_{\mathrm R}$ can be used directly to assemble the Gramian. In particular, for each residual row $q$ and active subdomain $\alpha\in T_q$,
\begin{align}
Z_{\chi(q),\mathcal C_\alpha}
=
(\Phi_\alpha)_{q,:},
\end{align}
where $\mathcal C_\alpha$ is given by \eqref{eq:jacobian-column-blocks}. Therefore, the Gramian blocks can be assembled as
\begin{align}\label{eq:reverse-gramian-clean}
G_{\alpha\beta}
=
\sum_{q\in Q_{\alpha\beta}}
Z_{\chi(q),\mathcal C_\alpha}^{\top}
Z_{\chi(q),\mathcal C_\beta}.
\end{align}
Similarly, the right-hand side is given by
\begin{align}\label{eq:reverse-jtr-clean}
(J_{\mathrm R}^\top \mathrm R)_{\mathcal C_\alpha}
=
\sum_{q\in Q_\alpha}
Z_{\chi(q),\mathcal C_\alpha}^{\top}\,\mathrm R_q .
\end{align}
Therefore, the Gauss--Newton system \eqref{eq:GN_system} can now be assembled.

\begin{figure}[t]
\centering
\begin{tikzpicture}[
row/.style={circle,draw,minimum size=7mm,inner sep=0pt,font=\small},
sub/.style={rectangle,draw,minimum width=8mm,minimum height=6mm,inner sep=1pt,font=\small},
edge/.style={gray!70,line width=0.55pt},
>=Latex
]
\node[row,fill=blue!10]   (q1) at (0,0)    {$q_1$};
\node[row,fill=orange!16] (q2) at (0,-0.9) {$q_2$};
\node[row,fill=blue!10]   (q3) at (0,-1.8) {$q_3$};
\node[row,fill=orange!16] (q4) at (0,-2.7) {$q_4$};
\node[sub] (a1) at (3.0,0)    {$\alpha=1$};
\node[sub] (a2) at (3.0,-0.9) {$\alpha=2$};
\node[sub] (a3) at (3.0,-1.8) {$\alpha=3$};
\node[sub] (a4) at (3.0,-2.7) {$\alpha=4$};
\draw[edge] (q1) -- (a1);
\draw[edge] (q1) -- (a2);
\draw[edge] (q2) -- (a2);
\draw[edge] (q3) -- (a3);
\draw[edge] (q3) -- (a4);
\draw[edge] (q4) -- (a4);
\node[font=\small] at (0,0.55) {residual rows};
\node[font=\small] at (3.0,0.55) {subdomain blocks};
\end{tikzpicture}
\caption{Row-subdomain incidence graph for the active sets $T_1=\{1,2\}$,
$T_2=\{2\}$, $T_3=\{3,4\}$, $T_4=\{4\}$: an edge joins $q_i$ to $\alpha$
whenever $\alpha\in T_i$. The coloring shown is valid because rows of one
color share no subdomain.}
\label{fig:row-coloring-bipartite}
\end{figure}
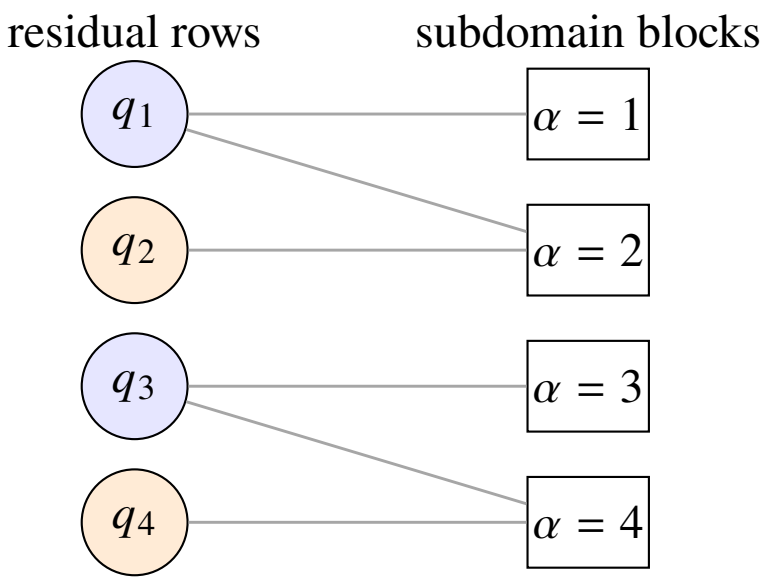

It remains to specify the row coloring used to construct the seed matrix $S$. In the FBPINN setting, this coloring problem is naturally formulated using the row-subdomain incidence graph
\begin{align}
\mathfrak{G}_{b}
=
(\mathcal{V}_{1}, \mathcal{V}_{2}, E),
\qquad
\mathcal{V}_{1} = \{1,\ldots,N_\Omega\},
\qquad
\mathcal{V}_{2} = \{1,\ldots,K\},
\end{align}
where $(q,\alpha)\in E$ if and only if $\alpha\in T_q$. Two residual rows
$i$ and $j$ cannot be compressed into the same reverse-mode seed whenever they
share an active subdomain, that is, whenever $T_i\cap T_j\neq\emptyset$. In $\mathfrak{G}_{b}$, this is precisely the case when there exists
$\alpha\in\mathcal{V}_{2}$ such that $i \;-\; \alpha \;-\; j$ is a path of length two, see Fig. \ref{fig:row-coloring-bipartite}. Hence, the row coloring problem is a partial distance-two coloring of $\mathfrak{G}_{b}$ on the residual-row vertex set $\mathcal{V}_{1}$: rows
connected by such a length-two path must receive different colors, while the
subdomain vertices are left uncolored. 
 
To solve the above problem, we use the greedy coloring
algorithm \cite[Algorithm~3.2]{gebremedhin2005color}. We define
\begin{align}\label{eq:used-color-set}
\mathcal{M}_\al
:=
\left\{
 c:\text{a previously colored row }i\text{ satisfies }\al \in T_i\text{ and }\chi(i)=c
\right\},
\end{align}
which records the row colors already present at subdomain $\al$. When row $i$ is
processed, all colors that have appeared at one of its active subdomains are forbidden, defining the set of forbidden colors
\begin{align}
F_{i} = \bigcup_{\al \in T_{i}} \mathcal{M}_{\al}.
\end{align}
The greedy rule assigns the smallest admissible color
\begin{align}
\chi(i): = \min \ \{ c \in \{1,2,.. \}  \, : \, c \notin F_{i} \}.
\end{align}
Afterwards, we update \eqref{eq:used-color-set} for all subdomains $\al \in T_{i}$ active at row $i$ as $\mathcal{M}_{\al} \leftarrow \mathcal{M}_{\al} \cup \{\chi(i) \}$. The number of colors produced by the greedy procedure depends on
the row ordering. We visit rows in decreasing order of $|T_i|$,
then lexicographically by the increasingly ordered elements of $T_i$,
and finally by increasing row index. Sorting adds
$\mathcal O(N_\Omega\log N_\Omega)$ preprocessing cost.
The row coloring procedure is summarized in
Algorithm~\ref{alg:active-set-coloring-simple}.

\begin{algorithm}[H]
\caption{Greedy FBPINN row coloring}
\label{alg:active-set-coloring-simple}
\begin{algorithmic}[1]
\Require Active-subdomain sets
$T_1,\ldots,T_{N_\Omega}\subseteq\{1,\ldots,K\}$
\Ensure An admissible coloring
$\chi:\{1,\ldots,N_\Omega\}\to\{1,\ldots,C_{\mathrm{row}}\}$

\State Order the rows by decreasing $|T_i|$, then
lexicographically by the ordered elements of $T_i$,
and finally by increasing row index $i$
\State Initialize $\mathcal M_\al\gets\emptyset$ for
$\al=1,\ldots,K$

\For{each row $i$ in this order}
    \State $F_i\gets\displaystyle\bigcup_{\al\in T_i}\mathcal M_\al$
    \State $\chi(i)\gets\min\{c\in\{1,2,\ldots\}:c\notin F_i\}$
    \State $\mathcal M_\al\gets\mathcal M_\al\cup\{\chi(i)\}$
    for all $\al\in T_i$
\EndFor

\State $C_{\mathrm{row}}\gets\max_i\chi(i)$
\State $\mathcal C_c^\chi\gets\{i:\chi(i)=c\}$,
$c=1,\ldots,C_{\mathrm{row}}$
\State \Return $\chi$ and
$\{\mathcal C_c^\chi\}_{c=1}^{C_{\mathrm{row}}}$
\end{algorithmic}
\end{algorithm}
 
For a vertex set $\mathcal V$, let $\Delta(\mathcal V):=\max_{v\in\mathcal V}\deg_{\mathfrak G_b}(v)$ denote its maximum vertex degree, i.e., the largest number of edges incident to any vertex in $\mathcal V$. For
$\mathfrak G_b$, we readily obtain
\begin{align}\label{eq:incidence-degrees}
\Delta(\mathcal V_1)
:=\max_{q=1,\ldots,N_\Omega}|T_q|, \qquad
\Delta(\mathcal V_2)
:=\max_{\al=1,\ldots,K}|Q_\al|,
\end{align}
where $Q_\al$ is defined in \eqref{eq: active-row set}. Since the sets $Q_\al$
are determined by the active-subdomain sets $T_q$, $\Delta(\mathcal V_2)$ is known before
training. Appealing to \cite[Lemma~3.13]{gebremedhin2005color}, the cost of
Algorithm~\ref{alg:active-set-coloring-simple} is then
$\mathcal O\bigl(|E|\,\Delta(\mathcal V_2)\bigr)$. Lemma~\ref{lem:overlap-multiplicity} yields
\begin{align}\label{eq:edge-count}
|E|
=\sum_{q=1}^{N_\Omega}|T_q|
\leq N_\Omega p_\delta^d.
\end{align}
Following Remark~\ref{rem:zonal-rounding} and using
\eqref{eq:no-clipping-spacing}, we obtain
\begin{align}
|Q_\al|
&=\frac{N_\Omega|\Omega_\al|}{|\Omega|}
+\mathcal O_{\delta,d}(1) =N_\Omega\left(\frac{\delta}{K_0-1+\delta}\right)^d
+\mathcal O_{\delta,d}(1) \leq \frac{N_\Omega}{K}\delta^d
+\mathcal O_{\delta,d}(1),
\end{align}
where $K=K_0^d$ and $\delta>1$.
Thus, for fixed $d$ and $\delta$, choosing $N_\Omega$ proportional
to $K$ ensures that the number of collocation points per subdomain
remains bounded under zonal sampling. Consequently,
$\Delta(\mathcal V_2)$ remains bounded, and the greedy coloring has linear complexity in $N_\Omega$.

Furthermore, let $\chi_{\mathrm{row}}$ denote the minimal number of colors of a map satisfying
\eqref{Sec 3: coloring map}, and let $C_{\mathrm{row}}^{\mathrm{greedy}}$ be the number
of colors returned by Algorithm~\ref{alg:active-set-coloring-simple}. For a bipartite
graph, \cite[Lemma~3.14 and Sec.~3.5.2]{gebremedhin2005color} give the following bounds
on the number of colors of a partial distance-two coloring:
\begin{align}\label{eq:gmp-bounds}
\Delta(\mathcal V_2)
\ \leq\
\chi_{\mathrm{row}}
\ \leq\
C_{\mathrm{row}}^{\mathrm{greedy}}
\ \leq\
\min\bigl\{\,\Delta(\mathcal V_1)\bigl(\Delta(\mathcal V_2)-1\bigr)+1,\ N_\Omega\,\bigr\}.
\end{align}
Lemma~\ref{lem:overlap-multiplicity} implies that
$\Delta(\mathcal V_1)=\max_q|T_q|\leq p_\delta^d$. Therefore, the upper bound on $C_{\mathrm{row}}^{\mathrm{greedy}}$ is quite pessimistic. In practice, however, we observe that the number of colors returned by Algorithm~\ref{alg:active-set-coloring-simple} is very close to $\Delta(\mathcal V_2)$.
 
Summarizing the above discussion, sparse forward and reverse mode require
$C_{\mathrm{col}}=n\,p_\delta^{d}$ JVPs and
$C_{\mathrm{row}}^{\mathrm{greedy}}$ VJPs, respectively. The forward mode count depends on the overlap ratio and the local
network size $n$, but not on $N_\Omega$, whereas the reverse mode count is
independent of $n$ and grows with the number of collocation points per
subdomain. Since both counts are known before training, the mode can be selected
a priori, taking into account the relative cost of a JVP and a VJP. In general, sparse forward mode is favored by
small local networks and fine sampling, whereas reverse mode is favored by
larger local networks and coarser sampling. However, as discussed in Section~\ref{sub:block-graph}, practically relevant configurations favor small $n$ and $\delta\leq 2$, while multiscale PDEs often require dense zonal sampling, rendering sparse forward-mode AD the preferable choice in many cases.

\section{Numerical solvers for the block-sparse Gauss--Newton system}
\label{sec:solvers}

In this section, we describe several approaches for solving the block-sparse Gauss--Newton system. Recall that, at each Gauss--Newton step, we solve the damped normal equations
\begin{align}\label{eq:solver-gn-system}
G_{\mu}(\theta)\delta\theta = -\nabla L(\theta),
\qquad
G_{\mu}(\theta):=G(\theta)+\mu I_P,
\qquad
\nabla  L(\theta):=J_{\mathrm R}(\theta)^\top \mathrm R(\theta),
\end{align}
where the damped Gauss--Newton blocks are given by
\begin{align*}
G_{\mu,\al\be}
=
\begin{cases}
G_{\al\al}+\mu I_n, & \al=\be,\\
G_{\al\be}, & \al\neq\be .
\end{cases}
\end{align*}
Since $G(\theta)$ is block symmetric, we store only $\operatorname{Sp}_{\mathrm{blk}}^+(G)$. For later use, we separate its
off-diagonal part by defining
\begin{align*}
\Eplus
:=
\{(\al,\be)\in\operatorname{Sp}_{\mathrm{blk}}^+(G):\al<\be\}.
\end{align*} 
The damping makes $G_\mu(\theta)$ positive definite for $\mu>0$, even when $G(\theta)$ is only positive semidefinite in general. Hence, direct Cholesky-type solvers and
preconditioned conjugate gradient (PCG) methods are natural choices to solve \eqref{eq:solver-gn-system}.

We begin with the direct solver. Direct factorization of a sparse matrix may introduce many additional nonzero entries in the factors, a phenomenon known as fill-in \cite{duff2017direct, saad2003iterative}. This fill-in may impose substantial additional storage requirements. Direct solvers are therefore attractive when high accuracy is required and the problem is still small enough that this fill-in remains manageable. In this case, the
upper-triangular block representation
$\operatorname{Sp}_{\mathrm{blk}}^+(G)$ is sufficient for solvers that accept
triangular input for symmetric positive definite matrices. For example, one may
apply a sparse Cholesky-type factorization
\begin{align}\label{eq:direct-cholesky}
P_{\mathrm{ch}}G_\mu P_{\mathrm{ch}}^\top
=
F_\mu F_\mu^\top,
\end{align}
where $P_{\mathrm{ch}}\in\R^{P\times P}$ is a fill-reducing
permutation matrix and $F_\mu\in\R^{P\times P}$ is a lower
triangular Cholesky factor \cite{duff2017direct}.

Since $G_\mu$ is symmetric positive definite, the system
\eqref{eq:solver-gn-system} can also be solved by conjugate gradients. Efficient use of CG, however, requires specifying two ingredients: a cheap realization of the Gramian-vector
product $v\mapsto G_\mu v$ and an effective preconditioner. Without
preconditioning, the convergence may be very slow for ill-conditioned
Gauss--Newton systems; cf. \cite{de2023operator}. We exploit the block structure of $G_\mu$ to construct both operations. We start with the Gramian-vector product. For $v=(v^1,\ldots,v^K)$ with
$v^\al\in\R^n$, the $\al$-th block of the product is given by
\begin{align}\label{eq:block-matvec-component}
(G_\mu v)^\al
=
\sum_{\be:(\al,\be)\in\operatorname{Sp}_{\mathrm{blk}}(G)}
G_{\al\be}\,v^\be
+
\mu v^\al .
\end{align}
This product is realized directly from the upper-triangular block storage
$\operatorname{Sp}_{\mathrm{blk}}^+(G)$. We initialize
\begin{align}
y^\al
\gets
(G_{\al\al}+\mu I_n)v^\al,
\qquad
\al=1,\ldots,K.
\end{align}
For each stored off-diagonal block $(\al,\be)\in\Eplus$,
symmetry gives the two updates
\begin{align}\label{eq:block-matvec-update}
y^\al
\leftarrow
y^\al+G_{\al\be}v^\be, \qquad
y^\be
\leftarrow
y^\be+G_{\al\be}^\top v^\al .
\end{align}
After all stored off-diagonal blocks have been processed, we obtain
$y=G_\mu v$. This operation is highly parallelizable and consists only of
dense matrix-vector products on small blocks, which are quite efficient on GPUs.

We proceed by constructing preconditioners from small patches of parameter blocks on the
Cartesian subdomain graph introduced in Section~\ref{sub:block-graph}; cf. also \cite{shang2025overlapping, chen2022overlapping}. Let $a_\al\in\mathcal I_{K_0}$ be the grid index of
subdomain $\al$. For each patch center $\al$, we consider
\begin{align}\label{eq:patch-definitions-solver-section}
\Patch_\al^{\mathrm{block}}
&:=\{\al\},
&
\Patch_\al^{\mathrm{cross}}
&:=\{\be:\|a_\be-a_\al\|_1\leq1\},
&
\Patch_\al^{\mathrm{king}}
&:=\{\be:\|a_\be-a_\al\|_\infty\leq1\}.
\end{align}
The block patch contains only the central parameter block and therefore leads to
a block Jacobi preconditioner. The cross and king patches additionally include,
respectively, the nearest coordinate-wise and diagonal neighbors of the central
subdomain. Near the boundary, only indices belonging to
$\mathcal I_{K_0}$ are retained.  

In the following, we write $\Patch_\al$ for any one of the three patch choices in
\eqref{eq:patch-definitions-solver-section}. Let 
\begin{align}
\Patch_\alpha=\{\beta_{\alpha,1},\ldots,\beta_{\alpha,s_\alpha}\}
\end{align}
be the patch centered at $\alpha$ with $s_\alpha:=|\Patch_\alpha|$. We define the restriction operator
\begin{align}\label{eq:patch-restriction-solver-section}
E_\alpha:\R^{Kn}\to\R^{s_\alpha n},
\qquad
E_\alpha v
=
\begin{bmatrix}
v^{\beta_{\alpha,1}}\\
\vdots\\
v^{\beta_{\alpha,s_\alpha}}
\end{bmatrix}.
\end{align}
Its transpose $E_\alpha^\top$ extends a local patch vector to the
global parameter space by setting all blocks outside the patch to zero.
The local patch matrix is defined by
$A_\al:=E_\al G_\mu E_\al^\top
\in\R^{s_\al n\times s_\al n}$ and has the block form
\begin{align}
A_\alpha
=
\begin{bmatrix}
G_{\beta_{\alpha,1}\beta_{\alpha,1}}+\mu I_n
&
\cdots
&
G_{\beta_{\alpha,1}\beta_{\alpha,s_\alpha}}
\\
\vdots
&
\ddots
&
\vdots
\\
G_{\beta_{\alpha,s_\alpha}\beta_{\alpha,1}}
&
\cdots
&
G_{\beta_{\alpha,s_\alpha}\beta_{\alpha,s_\alpha}}+\mu I_n
\end{bmatrix}.
\end{align}
At the $m$-th PCG iteration, let
$r_{\mathrm{CG}}^{(m)}
:=\nabla L(\theta) + G_\mu(\theta)\delta\theta^{(m)}$
denote the residual of the Gauss--Newton system
\eqref{eq:solver-gn-system}. The additive Schwarz preconditioner
is then given by
\begin{align}\label{eq:additive-schwarz-preconditioner-solver-section}
z^{(m)}
=M_{\mathrm{AS}}^{-1}r_{\mathrm{CG}}^{(m)}
:=
\sum_{\al=1}^{K}
E_\al^\top A_\al^{-1}E_\al r_{\mathrm{CG}}^{(m)}.
\end{align}
Each $A_\al$ is symmetric positive definite and admits a dense
Cholesky factorization. The local solves are independent and can
be performed in parallel. The factors can be computed once per
Gauss--Newton step and reused throughout the PCG iterations.
The preconditioner application is summarized in
Algorithm~\ref{alg:additive-schwarz}.

\begin{algorithm}[H]
\caption{Additive Schwarz application}
\label{alg:additive-schwarz}
\begin{algorithmic}[1]
\Require PCG residual $r_{\mathrm{CG}}^{(m)}$
and local Cholesky factors $A_\al=F_\al F_\al^\top$,
$\al=1,\ldots,K$
\Ensure $z^{(m)}=M_{\mathrm{AS}}^{-1}r_{\mathrm{CG}}^{(m)}$
\State Initialize $z^{(m)}\gets0$
\For{$\al=1,\ldots,K$}
    \State Solve $A_\al z_\al=E_\al r_{\mathrm{CG}}^{(m)}$
    using the Cholesky factors
    \State $z^{(m)}\gets z^{(m)}+E_\al^\top z_\al$
\EndFor
\State \Return $z^{(m)}$
\end{algorithmic}
\end{algorithm}

The Gauss--Newton solver is summarized in
Algorithm~\ref{alg:block-gn-linear-solver}. Both direct and iterative branches use the same stored
upper-triangular blocks. A direct solve is preferable when sparse factorization
is affordable and high accuracy is required; PCG avoids a global factorization
and is more suitable for larger systems when the local preconditioner keeps the
iteration count moderate. However, we remark that patch-based preconditioners may themselves
require substantial additional memory, since the Cholesky factors of the local
patch matrices must be stored and reused throughout the PCG iterations. This
cost increases with the patch size and can become significant for the cross and,
in particular, king Schwarz preconditioners.

\begin{algorithm}[H]
\caption{Linear solve for the damped block Gauss--Newton system}
\label{alg:block-gn-linear-solver}
\begin{algorithmic}[1]
\Require Stored Gramian blocks
$\{G_{\al\be}:(\al,\be)\in\operatorname{Sp}_{\mathrm{blk}}^+(G)\}$,
gradient $\nabla L(\theta)$, and damping $\mu>0$
\Require Solver mode
$\mathsf{mode}\in\{\textsc{direct},\textsc{pcg}\}$;
for PCG: patch type, tolerance $\varepsilon_{\mathrm{CG}}$,
and maximum iteration count $m_{\max}$
\Ensure $\delta\theta\approx G_\mu^{-1}\nabla L(\theta)$

\If{$\mathsf{mode}=\textsc{direct}$}
    \State Assemble $G_\mu$ by adding $\mu I_n$ to the diagonal blocks
    \State Compute the sparse Cholesky factorization
    $P_{\mathrm{ch}}G_\mu P_{\mathrm{ch}}^\top=F_\mu F_\mu^\top$
    \State Solve $F_\mu w=P_{\mathrm{ch}}\nabla L(\theta)$ and
    $F_\mu^\top\widehat{\delta\theta}=w$
    \State $\delta\theta\gets
    P_{\mathrm{ch}}^\top\widehat{\delta\theta}$
\Else
    \State Define $v\mapsto G_\mu v$ using
    \eqref{eq:block-matvec-component}--\eqref{eq:block-matvec-update}
    \State Compute the local Cholesky factorizations
    $A_\al=F_\al F_\al^\top$, $\al=1,\ldots,K$
    \State Apply PCG with the additive Schwarz preconditioner
    from Algorithm~\ref{alg:additive-schwarz}
    \State Stop when
    $\|r_{\mathrm{CG}}^{(m)}\|_2
    \leq\varepsilon_{\mathrm{CG}}\|\nabla L(\theta)\|_2$
    or $m=m_{\max}$
\EndIf
\State \Return $\delta\theta$
\end{algorithmic}
\end{algorithm}

\section{Numerical examples}

We evaluate the proposed methodology on three PDE problems:
a multiscale Poisson equation with several spatial frequencies,
an elliptic equation with rapidly oscillating coefficients in
a mixed first-order formulation, and Stokes flow in perforated
domains. We assess the convergence, accuracy, and computational
cost of the proposed block-sparse Gauss--Newton method, including
comparisons with Adam \cite{kingma2014adam}. In our numerical experiments, sparse direct
solvers generally remain the preferable option. Nevertheless,
we also investigate preconditioned conjugate gradients with block
Jacobi and additive Schwarz preconditioners on the multiscale
Poisson problem to assess the capabilities of iterative solvers. We implement our methods using the JAX programming framework \cite{jax2018github}.
For sparse direct solves, we employ the NVIDIA cuDSS library%
\footnote{\url{https://docs.nvidia.com/cuda/cudss/}},
accessed from JAX through the Spineax interface \cite{viljoen2026scaling}.
Our implementation combines ASDEX \cite{hill2026asdex} for sparse reverse-mode AD with dedicated implementations of sparse forward-mode AD and iterative solvers. All numerical experiments are performed on an NVIDIA A100 GPU
with 80 GB of memory.

\subsection{Multiscale Poisson problem}

Our first benchmark is the Poisson problem
\begin{equation}\label{eq:poisson_problem}
    -\Delta u = f
    \quad \text{in } \Omega=(0,1)^2,
    \qquad
    u=0
    \quad \text{on } \partial\Omega .
\end{equation}
Following \cite{dolean2024multilevel}, we choose the multiscale right-hand side
\begin{equation*}
    f(x_1,x_2)
    =
    \frac{1}{N_{\mathrm f}}
    \sum_{j=1}^{N_{\mathrm f}}
    2(\pi\omega_j)^2
    \sin(\pi\omega_j x_1)
    \sin(\pi\omega_j x_2),
    \qquad \omega_j=2^j,
\end{equation*}
for which the exact solution is given by
\begin{equation}\label{eq:poisson_exact}
    u_{\mathrm{ex}}(x_1,x_2)
    =
    \frac{1}{N_{\mathrm f}}
    \sum_{j=1}^{N_{\mathrm f}}
    \sin(\pi\omega_j x_1)
    \sin(\pi\omega_j x_2).
\end{equation}
Multiplying the
FBPINN ansatz \eqref{eq:local-fbpinn-contributions} by the factor
\begin{equation}
    l(x)
    =
    \prod_{\ell=1}^{2}
    \tanh\!\left(\frac{x_\ell}{\sigma_c}\right)
    \tanh\!\left(\frac{1-x_\ell}{\sigma_c}\right),
    \qquad
    \sigma_c=\omega_{\max}^{-1},
\end{equation}
we replace $u_\theta$ in
\eqref{eq:local-fbpinn-contributions} by $l\,u_\theta$. Since $l=0$ on
$\partial\Omega$, the homogeneous Dirichlet boundary condition is satisfied exactly by the neural network ansatz, and training
involves only the interior Poisson residual.

We consider two configurations. For $N_{\mathrm f}=6$, we use a
$64\times64$ subdomain decomposition. Each subdomain carries a fully
connected network $2\to12\to1$ with $\tanh$ activation, totaling
$n=49$ parameters per subnetwork. The collocation points are drawn
once using zonal sampling with $100$ points per subdomain. In total,
we have $P=200{,}704$ trainable parameters and $N_\Omega=409{,}600$
collocation points. For $N_{\mathrm f}=5$, we use a coarser
$32\times32$ decomposition with correspondingly wider local networks
$2\to50\to1$, that is, $n=201$ parameters per subnetwork, and again
$100$ zonally sampled points per subdomain, giving $P=205{,}824$
trainable parameters and $N_\Omega=102{,}400$ collocation points.
We note that the configuration for $N_{\mathrm f}=5$ was deliberately
chosen to make reverse-mode AD preferable to
forward-mode AD, while the configuation with smaller network sizes can be more favourable in practice in terms of compute time. 

For $1<\delta\le2$, Proposition~\ref{prop:block-coloring} yields
$p_\delta^2=4$ block colors, hence $C_{\mathrm{col}}=4n$
JVPs, in both configurations. For $N_{\mathrm f}=6$ this gives
$C_{\mathrm{col}}=196$ JVPs, whereas
Algorithm~\ref{alg:active-set-coloring-simple} yields
$252$ and $392$ row colors for $\delta=1.6$ and $\delta=2.0$.
For $N_{\mathrm f}=5$ the wider subnetworks raise the forward
count to $C_{\mathrm{col}}=804$ JVPs, while the row coloring, which is
governed by the overlap geometry and not by the network width, yields only
$C_{\mathrm{row}}=252$ and $C_{\mathrm{row}}=380$ for $\delta=1.6$ and $\delta = 2.0$, respectively. In each
configuration we use the mode with fewer colors. For reverse-mode AD, we process the
colors in batches of size $32$. Observe that
$C_{\mathrm{col}} \ll P$ and $C_{\mathrm{row}} \ll N_\Omega$. The resulting
Gramian has a sparsity of $99.78\%$ for $N_{\mathrm f}=6$ and
$99.16\%$ for $N_{\mathrm f}=5$. 

\begin{figure}[H]
    \centering
    \includegraphics[width=0.92\linewidth]{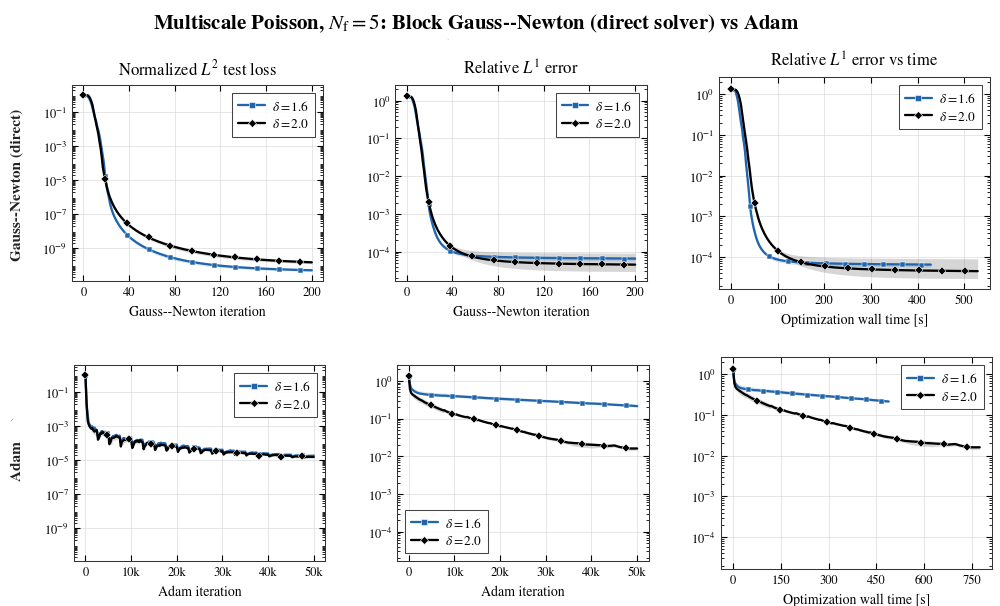}
    \caption{Multiscale Poisson problem with $N_{\mathrm f} = 5$: convergence of direct block Gauss--Newton and Adam for different overlap ratios $\delta$, measured by normalized (by its initial value) test loss, relative $L^1$ error, and optimization wall time. Curves are means over 5
initialization seeds; shaded bands indicate the minimum--maximum range.}
    \label{fig:direct_solver_poisson_N5}
\end{figure}

\begin{figure}[H]
    \centering
    \includegraphics[width=0.92\linewidth]{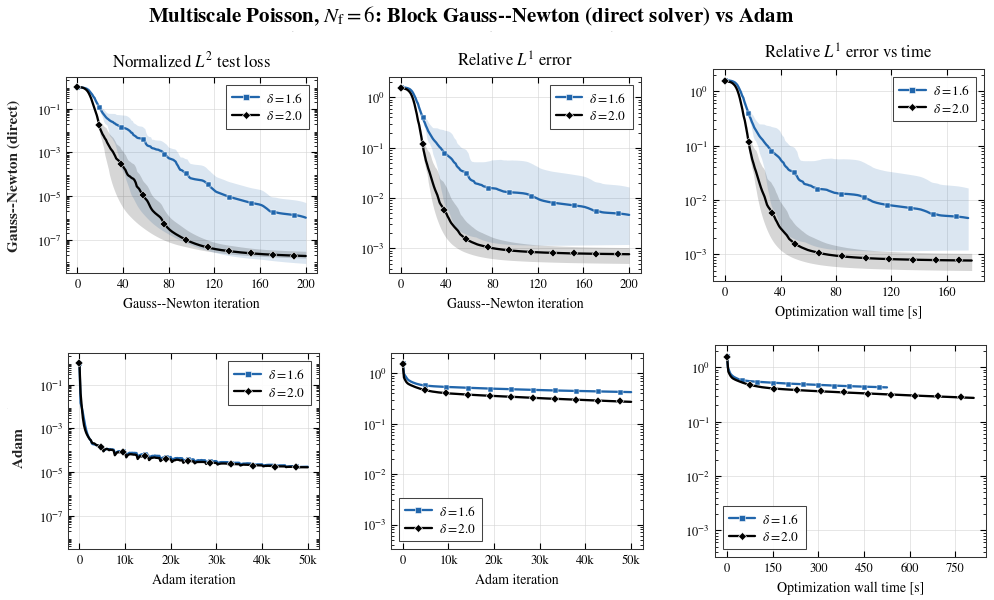}
    \caption{Multiscale Poisson problem with $N_{\mathrm f} = 6$: convergence of direct block Gauss--Newton and Adam for different overlap ratios $\delta$, measured by normalized (by its initial value) test loss, relative $L^1$ error, and optimization wall time. Curves are means over 5
initialization seeds; shaded bands indicate the minimum--maximum range.}
    \label{fig:direct_solver_poisson_N6}
\end{figure}

\begin{figure}[H]
    \centering
    \includegraphics[width=0.9\linewidth]{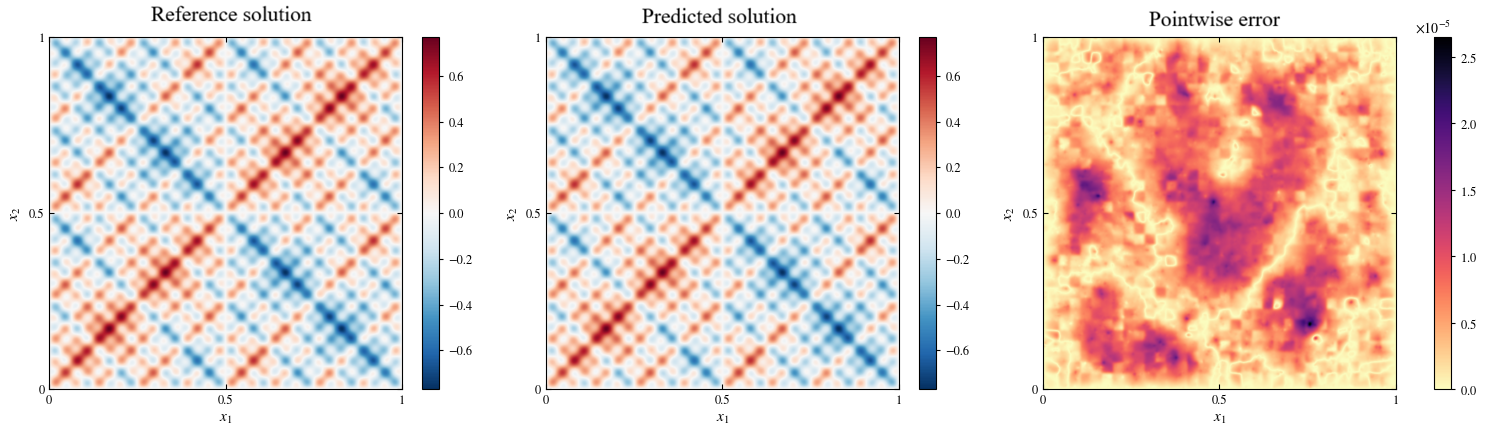}
    \caption{Multiscale Poisson problem with $N_{\mathrm f} = 5$: reference solution, FBPINN prediction obtained with the Gauss--Newton, and the corresponding pointwise absolute error.}
    \label{fig:poisson_solution_32}
\end{figure}

\begin{figure}[H]
    \centering
    \includegraphics[width=0.92\linewidth]{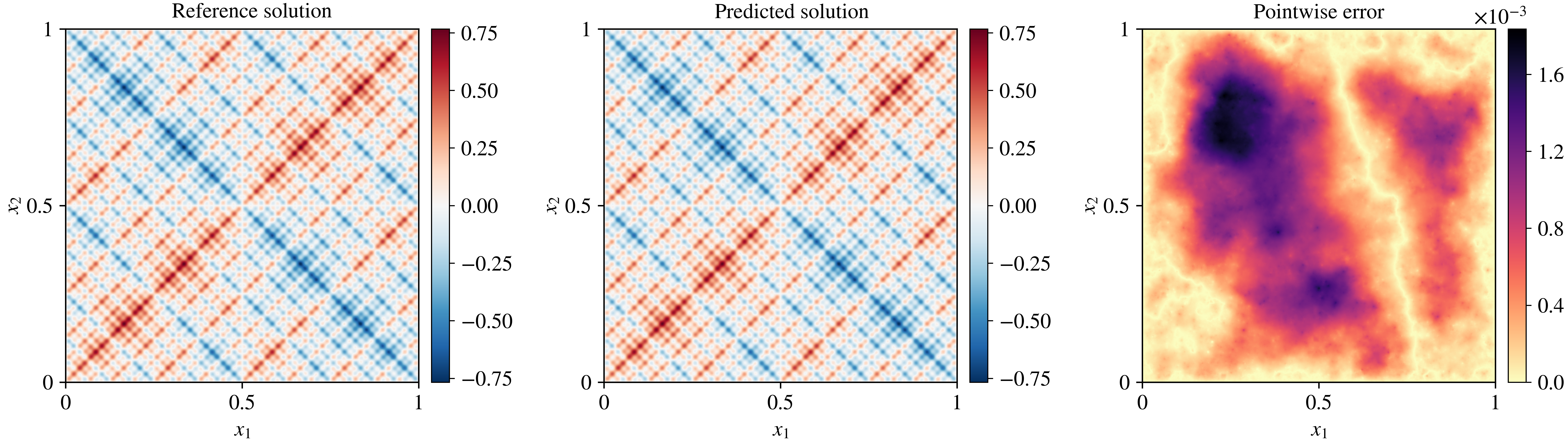}
    \caption{Multiscale Poisson problem with $N_{\mathrm f} = 6$: reference solution, FBPINN prediction obtained with the Gauss--Newton, and the corresponding pointwise absolute error.}
    \label{fig:poisson_solution_64}
\end{figure}

We first apply the block Gauss--Newton method with the sparse direct solver of
Section~\ref{sec:solvers}, using damping $\mu=10^{-3}$
for $N_{\mathrm f}=6$ and $\mu=10^{-6}$ for
$N_{\mathrm f}=5$, and selecting the step size at each iteration using the
line search over the geometric grid $\eta \in \{2^{-s} : s=0,\ldots,10\}$.  We compare it with Adam~\cite{kingma2014adam} using gradients
computed by standard reverse-mode AD in double precision,
the standard Adam hyperparameters
$\beta_1^{\mathrm{Ad}}=0.9$,
$\beta_2^{\mathrm{Ad}}=0.999$,
$\varepsilon^{\mathrm{Ad}}=10^{-8}$,
a learning rate of $\alpha_k=10^{-2}\cdot 0.9^{k/2000}$,
and resampling of the collocation points every $2500$ iterations. Fig.~\ref{fig:direct_solver_poisson_N5} and Fig.~\ref{fig:direct_solver_poisson_N6} show that, although Adam makes some learning progress for larger overlap ratio, this progress is very slow, both in terms of iteration count and optimization wall time, and ultimately yields low accuracy, see also Table \ref{tab:poisson-overlap-accuracy}. In contrast, the block Gauss--Newton approach achieves high accuracy, especially for the larger overlap ratio $\delta=2$, while retaining a very competitive run time given the number of subdomains, parameters, and collocation points. Fig.~\ref{fig:poisson_solution_32} and Fig.~\ref{fig:poisson_solution_64} show representative comparisons between the Gauss--Newton solutions
and the reference solutions for $N_{\mathrm f}=5$ and
$N_{\mathrm f}=6$, respectively. We remark that for fine subdomain decompositions, such as the one
considered for $N_{\mathrm f}=6$, a larger overlap is generally preferable. Indeed, the derivatives of the strongly localized windows enter the PDE residual and introduce additional stiffness into the optimization problem, which is reduced with larger $\delta$. This may explain the substantially larger variability in the errors across independent runs observed for $N_{\mathrm f}=6$ and $\delta=1.6$, as reflected by the wider error bands in Fig.~\ref{fig:direct_solver_poisson_N6}.

\begin{table}
  \centering
  \footnotesize
  \setlength{\tabcolsep}{4pt}
  \renewcommand{\arraystretch}{1.15}
  \begin{tabular}{|l|l|r|r|r|r|}
    \hline
    & & \multicolumn{2}{c|}{$N_f=5$} & \multicolumn{2}{c|}{$N_f=6$} \\
    \cline{3-6}
    & & \multicolumn{4}{c|}{Overlap ratio $\delta$} \\
    \cline{3-6}
    Method & Quantity
    & $1.6$ & $2.0$ & $1.6$ & $2.0$ \\
    \hline
    Gauss--Newton (direct) & Rel. $L^1$ error
    & $(6.55\pm0.55)\times10^{-5}$ & $(4.52\pm2.61)\times10^{-5}$ & $(4.65\pm6.72)\times10^{-3}$ & $(7.69\pm2.40)\times10^{-4}$ \\
    \hline
     & Norm. $L^2$ test loss
    & $(4.84\pm0.64)\times10^{-11}$ & $(1.42\pm0.15)\times10^{-10}$ & $(1.02\pm2.23)\times10^{-6}$ & $(1.78\pm0.64)\times10^{-8}$ \\
    \hline
     & Time [s]
    & $428.3$ & $529.0$ & $175.4$ & $178.0$ \\
    \hline\hline
    PCG (no preconditioner) & Rel. $L^1$ error
    & $(5.28\pm1.44)\times10^{-2}$ & $(1.73\pm0.19)\times10^{-2}$ & $(2.69\pm0.08)\times10^{-1}$ & $(4.45\pm0.45)\times10^{-2}$ \\
    \hline
     & Norm. $L^2$ test loss
    & $(3.49\pm0.18)\times10^{-5}$ & $(3.35\pm0.10)\times10^{-5}$ & $(2.77\pm0.03)\times10^{-5}$ & $(2.43\pm0.04)\times10^{-5}$ \\
    \hline
     & Time [s]
    & $392.6$ & $493.4$ & $156.2$ & $158.8$ \\
    \hline\hline
    PCG (block Jacobi) & Rel. $L^1$ error
    & $(4.83\pm0.61)\times10^{-1}$ & $(1.66\pm1.00)\times10^{-2}$ & $(5.38\pm0.26)\times10^{-1}$ & $(1.26\pm0.15)\times10^{-1}$ \\
    \hline
     & Norm. $L^2$ test loss
    & $(8.03\pm1.89)\times10^{-4}$ & $(4.36\pm0.15)\times10^{-8}$ & $(1.62\pm0.46)\times10^{-4}$ & $(2.92\pm2.16)\times10^{-7}$ \\
    \hline
     & Time [s]
    & $425.0$ & $525.8$ & $169.6$ & $172.3$ \\
    \hline\hline
    PCG (cross Schwarz) & Rel. $L^1$ error
    & $(4.05\pm0.48)\times10^{-3}$ & $(1.06\pm0.25)\times10^{-3}$ & $(3.41\pm0.61)\times10^{-2}$ & $(6.71\pm0.83)\times10^{-3}$ \\
    \hline
     & Norm. $L^2$ test loss
    & $(1.75\pm0.10)\times10^{-8}$ & $(4.16\pm0.30)\times10^{-8}$ & $(2.64\pm0.20)\times10^{-8}$ & $(5.67\pm0.07)\times10^{-8}$ \\
    \hline
     & Time [s]
    & $737.2$ & $835.5$ & $268.9$ & $273.5$ \\
    \hline\hline
    PCG (king Schwarz) & Rel. $L^1$ error
    & $(1.93\pm0.51)\times10^{-3}$ & $(9.98\pm2.67)\times10^{-4}$ & $(1.34\pm0.11)\times10^{-2}$ & $(3.57\pm1.09)\times10^{-3}$ \\
    \hline
     & Norm. $L^2$ test loss
    & $(1.62\pm0.07)\times10^{-8}$ & $(4.17\pm0.27)\times10^{-8}$ & $(1.97\pm0.23)\times10^{-8}$ & $(5.62\pm0.15)\times10^{-8}$ \\
    \hline
     & Time [s]
    & $1327.6$ & $1346.3$ & $400.5$ & $390.2$ \\
    \hline\hline
    Adam & Rel. $L^1$ error
    & $(2.14\pm0.13)\times10^{-1}$ & $(1.61\pm0.12)\times10^{-2}$ & $(4.26\pm0.06)\times10^{-1}$ & $(2.73\pm0.08)\times10^{-1}$ \\
    \hline
     & Norm. $L^2$ test loss
    & $(1.79\pm0.06)\times10^{-5}$ & $(1.50\pm0.02)\times10^{-5}$ & $(1.76\pm0.02)\times10^{-5}$ & $(1.66\pm0.02)\times10^{-5}$ \\
    \hline
     & Time [s]
    & $488.5$ & $774.5$ & $525.8$ & $810.9$ \\
    \hline
 \end{tabular}
\caption{Multiscale Poisson problem: final relative $L^1$ errors,
normalized $L^2$ test losses, and optimization times (s) for
different overlap ratios $\delta$. Results are means $\pm$
standard deviations over five seeds, using $200$ Gauss--Newton
or $50{,}000$ Adam iterations.}
  \label{tab:poisson-overlap-accuracy}
\end{table}

\begin{figure}
\centering
\includegraphics[width=0.77\linewidth]{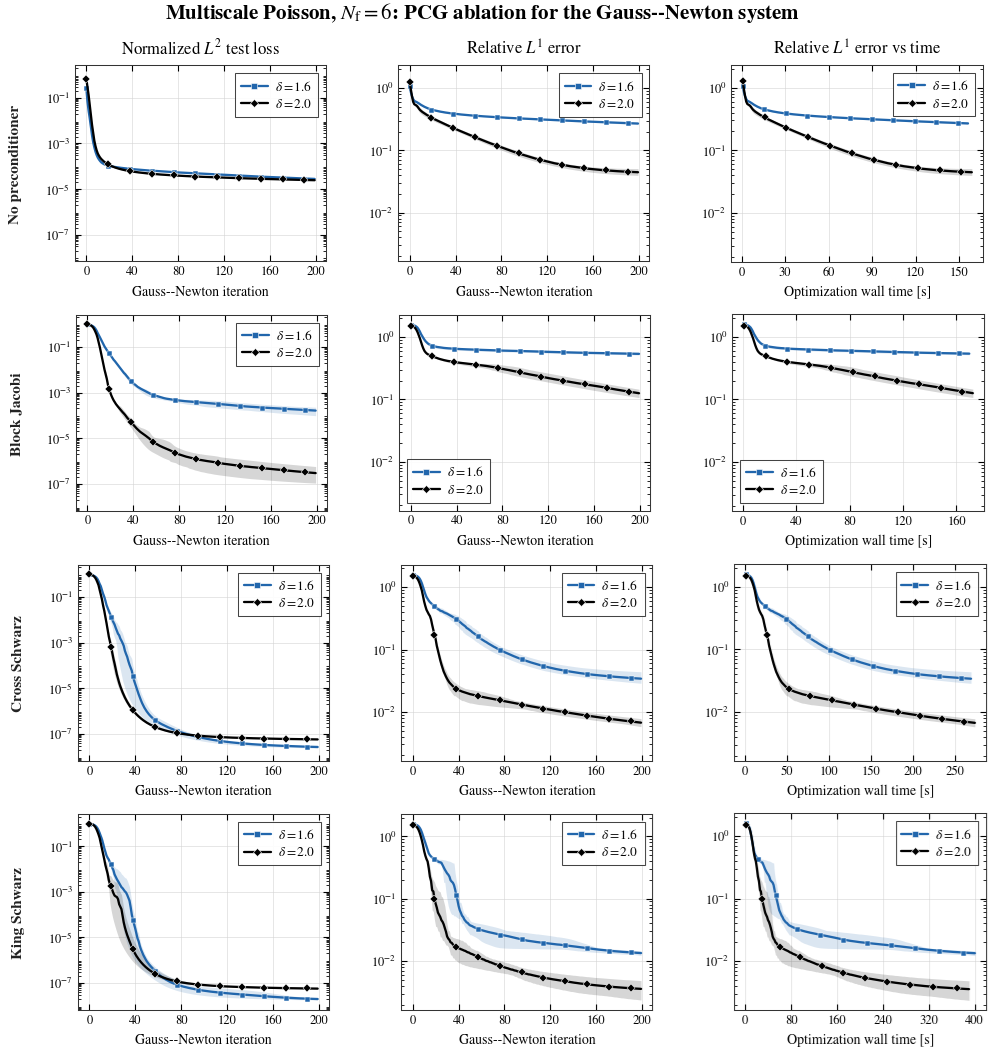}
\caption{Multiscale Poisson problem, $N_{\mathrm f}=6$: PCG with different preconditioners and overlap ratios $\delta$, showing the $L^2$ test loss normalized by its initial value, relative $L^1$ error, and optimization wall time. Curves show means over 5 initialization seeds; shaded bands indicate the minimum--maximum range.}
\label{fig:iterative_solver_poisson}
\end{figure}

We further investigate the iterative solvers from Section~\ref{sec:solvers},
comparing PCG without preconditioning, block Jacobi, and additive Schwarz
with cross and king patches. We use $\mu=10^{-2}$, $\varepsilon_{\mathrm{CG}}=10^{-5}$, and $m_{\max}=200$
in Algorithm~\ref{alg:block-gn-linear-solver}, with the same line search as for the direct solver. The results are presented in Table~\ref{tab:poisson-overlap-accuracy} for both $N_{\mathrm f}=5$ and $N_{\mathrm f}=6$, and in Fig.~\ref{fig:iterative_solver_poisson} for $N_{\mathrm f}=6$. Preconditioning substantially improves PCG accuracy,
with king Schwarz performing best among the iterative methods, albeit at a higher cost and lower accuracy than the sparse direct solver. Further improvements may require more effective regularization to improve the conditioning of the Gauss--Newton system; see, e.g., \cite{van2026local}.

Table~\ref{tab:poisson-stage-timing} shows that the sparse-AD stage remains a significant part of the Gauss--Newton iteration in both configurations, so that the differentiation, and not the linear algebra, is the main cost of the method. The sparse direct solve is comparatively inexpensive, which is what makes the block-sparse Gramian attractive at these problem sizes. The assembly of the right-hand side and of the Gramian from the compressed Jacobian is negligible throughout. The preconditioners built on the richer patches, cross and king Schwarz, are by contrast expensive, costing more per iteration than the differentiation itself. 

\begin{table}[!htbp]
  \centering
  \vspace{0.9em}
  \small
  \setlength{\tabcolsep}{4pt}
  \begin{tabular}{@{}lccccccc@{}}
  \toprule
  & \multicolumn{2}{c}{Common stages}
  & \multicolumn{5}{c}{Linear solve} \\
  \cmidrule(lr){2-3}\cmidrule(l){4-8}
  $N_{\mathrm f}$
  & Sparse AD & Assembly
  & Direct & None
  & \shortstack{Block}
  & \shortstack{Cross}
  & \shortstack{King} \\
  \midrule
  $5$ & 1.860 & 0.018 & 0.711 & 0.522 & 0.685 & 2.233 & 4.787 \\
  $6$ & 0.583 & 0.012 & 0.218 & 0.114 & 0.181 & 0.687 & 1.271 \\
  \bottomrule
\end{tabular}
\caption{Multiscale Poisson: Mean stage times (s) per Gauss--Newton
iteration for $\delta=2.0$. Linear solvers are cuDSS (Direct)
and PCG without (None) or with the indicated preconditioners.
Assembly includes the gradient and the Gramian.}
 \label{tab:poisson-stage-timing}
\end{table}

\subsection{Mixed formulation for elliptic problems with oscillatory coefficients}
\label{subsec:oscillatory-mixed}

\begin{figure}[H]
\centering
\includegraphics[width=0.75\linewidth]{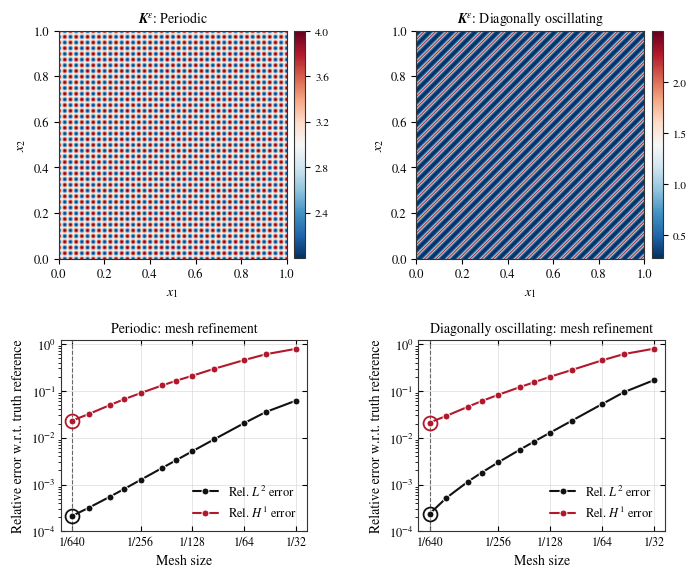}
\caption{Oscillatory coefficients (top) and FEM mesh-refinement validation
(bottom): periodic coefficient (left) and diagonally oscillating
coefficient (right), with $\varepsilon=0.05$. The lower panels show relative
$L^2$ errors  and $H^1$ errors of the $P_1$ solutions
with respect to the finest computed FEM solution, $h=1/812$.
Dashed lines and open circles mark the reference mesh $h=1/640$ used
in the FBPINN experiments.}
\label{fig:coeff_and_fem}
\end{figure}

We consider elliptic PDEs with rapidly oscillating coefficients,
which arise in the modeling of heterogeneous
media; cf.~\cite{wu2002analysis,low2020influence,hintermuller2026hybrid}.
Specifically, we investigate the following boundary-value problem:
\begin{align}\label{eq:mixed-oscillatory-problem}
-\operatorname{div}\big(
\boldsymbol K^{\varepsilon}(x)
(\nabla u(x)+\boldsymbol e_1)
\big)
=f(x) \quad\text{in }\Omega=(0,1)^2, \qquad
u=0 \quad\text{on }\partial\Omega,
\end{align}
where $\boldsymbol e_1=(1,0)^\top$ and the right-hand side is given by
\begin{align*}
f(x)
:=10\exp\bigg(-100\sum_{j=1}^2
\big(x_j-\frac12\big)^2\bigg).
\end{align*}
We consider the following coefficient fields:
\begin{align}\label{eq:mixed-oscillatory-coefficients}
\boldsymbol K_{\mathrm{diag}}^{\varepsilon}(x)
:=\frac{1}{2+1.6\sin\bigl(2\pi(x_1-x_2)/\varepsilon\bigr)}, \qquad
\boldsymbol K_{\mathrm{per}}^{\varepsilon}(x)
:=3+\prod_{j=1}^2
\cos\big(2\pi x_j/\varepsilon\big).
\end{align}
In both coefficients, we set $\varepsilon=0.05$. The first coefficient varies along the diagonal direction,
whereas the second is periodic in both coordinate directions, see Fig.\ref{fig:coeff_and_fem}.

Instead of minimizing the second-order residual of
\eqref{eq:mixed-oscillatory-problem}, we introduce the flux
\begin{equation}\label{eq:mixed-oscillatory-flux}
\boldsymbol\phi:=\boldsymbol K^{\varepsilon}
(\nabla u+\boldsymbol e_1)
\end{equation}
and obtain the equivalent first-order system
\begin{equation}\label{eq:mixed-oscillatory-system}
\begin{aligned}
-\operatorname{div}\boldsymbol\phi&=f &&\text{in }\Omega,\\
\boldsymbol\phi-\boldsymbol K^{\varepsilon}
(\nabla u+\boldsymbol e_1)&=0 &&\text{in }\Omega,\\
u&=0 &&\text{on }\partial\Omega.
\end{aligned}
\end{equation}
The first equation expresses conservation and the second imposes the
constitutive relation. The corresponding weighted first-order system
least-squares (FOSLS) formulation~\cite{cai1994first} is given by
\begin{equation}\label{eq:mixed-oscillatory-fosls}
\min_{\substack{v\in H_0^1(\Omega)\\
\boldsymbol\psi\in H(\operatorname{div};\Omega)}}
\left\{
\frac12\|{-\operatorname{div}\boldsymbol\psi-f}\|_{L^2(\Omega)}^2
+\frac{\lambda_{\mathrm{flux}}}{2}
\|\boldsymbol\psi-\boldsymbol K^{\varepsilon}
(\nabla v+\boldsymbol e_1)\|_{L^2(\Omega)}^2
\right\}.
\end{equation}
This formulation requires only first-order
spatial derivatives of the network outputs and avoids spatial differentiation of
the oscillatory coefficient, which can substantially worsen the
conditioning of the optimization problem; cf.~\cite{li2024solving}.
It also approximates the flux \eqref{eq:mixed-oscillatory-flux} directly, which is a quantity of interest in heat
transfer problems.

To assess the FEM reference, we perform a mesh-refinement study with
$P_1$ finite elements, comparing coarser solutions with the finest computed
solution at $h=1/812$; see Fig.~\ref{fig:coeff_and_fem}. The reference used
in the FBPINN experiments is $u_h$ with $h=1/640$. We report the $L^2$ and $H^1$ relative errors
\begin{equation}\label{eq:mixed-relative-errors}
E_{L^2}:=\frac{\|u_{\theta_u}-u_h\|_{L^2(\Omega)}}
{\|u_h\|_{L^2(\Omega)}},
\qquad
E_{H^1}:=\frac{\|\nabla u_{\theta_u}-\nabla u_h\|_{L^2(\Omega)}}
{\|\nabla u_h\|_{L^2(\Omega)}}.
\end{equation}
Here, the $H^1$ error is a relative seminorm error. However, under homogeneous Dirichlet boundary conditions, it defines a norm on $H_0^1(\Omega)$.
For validation, we use the FEM triangulation and approximate
the integrals in the above error measures using a six-point
Dunavant quadrature rule \cite{dunavant1985high} on the FEM triangles.

We approximate $u$ and $\boldsymbol\phi$ by separate FBPINN networks
$u_{\theta_u}$ and $\boldsymbol\phi_{\theta_{\boldsymbol\phi}}$,
with parameters $\theta_u\in\mathbb R^{P_u}$ and
$\theta_{\boldsymbol\phi}\in\mathbb R^{P_{\boldsymbol\phi}}$.
The overall parameter vector is
$\theta=(\theta_u,\theta_{\boldsymbol\phi})\in\mathbb R^P$,
where $P=P_u+P_{\boldsymbol\phi}$.
The homogeneous Dirichlet condition is imposed exactly by multiplying
the raw field output by $\ell(x):=16x_1x_2(1-x_1)(1-x_2)$.
The pointwise conservation and constitutive residuals are given by
\begin{align*}
r^{\mathrm{cons}}(x;\theta)
:=-\operatorname{div}\boldsymbol\phi_{\theta_{\boldsymbol\phi}}(x)-f(x), \qquad
\boldsymbol r^{\mathrm{const}}(x;\theta)
:=\boldsymbol\phi_{\theta_{\boldsymbol\phi}}(x)
-\boldsymbol K^{\varepsilon}(x)
\bigl(\nabla u_{\theta_u}(x)+\boldsymbol e_1\bigr).
\end{align*}
Thus, $N_{\mathrm{eq}}=3$, and we obtain
\begin{align*}
\mathcal R(\theta)(x)=
\begin{bmatrix}
r^{\mathrm{cons}}(x;\theta)\\
\sqrt{\lambda_{\mathrm{flux}}}\,\boldsymbol r^{\mathrm{const}}(x;\theta)
\end{bmatrix},
\qquad
\mathrm R(\theta)\in\mathbb R^{3N_\Omega},
\qquad
J_{\mathrm R}(\theta)\in\mathbb R^{3N_\Omega\times P},
\end{align*}
where $\mathrm R$ is the quadrature-weighted discrete residual
and $\lambda_{\mathrm{flux}}=\varepsilon^{-2}=400$. The extension of the preceding analysis from $N_{\mathrm{eq}}=1$
to $N_{\mathrm{eq}}=3$ is straightforward. The residual Jacobian
is formed by stacking the Jacobians of the three scalar residual
components, while the field and flux parameters are grouped by
subdomain. This yields blocks
$\Phi_\alpha\in\mathbb R^{3N_\Omega\times n}$,
where $n$ counts all parameters per subdomain.
The sparse AD constructions \eqref{eq:foward-compressed-jac} and
\eqref{eq:reverse-compressed-jacobian-clean} apply to the stacked
Jacobian, with reverse-mode coloring performed on individual
scalar residual rows.

Both benchmarks use a $40\times40$ decomposition with overlap ratio
$\delta=2$. Each subdomain carries two fully connected networks with one
$\tanh$ hidden layer: a field network $2\to15\to1$ with $61$ parameters
and a flux network $2\to20\to2$ with $102$ parameters. Thus, $n=163$ and $P=260{,}800$.
Zonal sampling uses a budget of $300$ points per subdomain, giving
$N_\Omega=480{,}000$. The collocation points are fixed for Gauss--Newton
and resampled every $2{,}500$ iterations for Adam. We compare $100$ direct block Gauss--Newton iterations with
$50{,}000$ Adam iterations over the same five initialization seeds.
Gauss--Newton uses fixed damping $\mu=10^{-6}$ and geometric grid line
search. For the Adam algorithm, we use the same setup as in the first Poisson example. We use sparse forward-mode AD with
$C_{\mathrm{blk}}=4$ block colors and
$C_{\mathrm{col}}=nC_{\mathrm{blk}}=652$ seed directions, whereas
row coloring gives $C_{\mathrm{row}}=3456$ colors.
To reduce peak memory usage, we split the block colors into two
batches of $B_{\mathrm{fwd}}=2$ colors, at the expense of increased
computation time on a single GPU. 

\begin{figure}[H]
\centering
\includegraphics[width=1.00\linewidth]{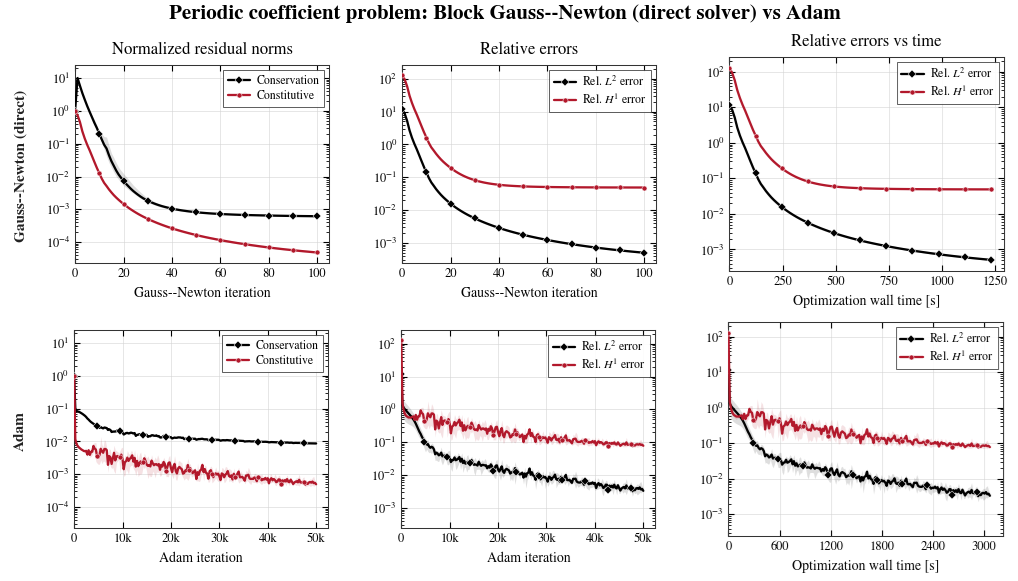}
\caption{Periodic coefficient: direct block Gauss--Newton (top) versus
Adam (bottom), with $\delta=2$. Left: conservation and constitutive test
residual norms, each normalized by its own initial value. Middle:
relative $L^2$ and $H^1$-norm errors versus iteration. Right: the same
errors versus optimization wall time. Curves show means over 5 initialization seeds; shaded bands indicate the minimum--maximum range.}
\label{fig:mixed-periodic-convergence}
\end{figure}

\begin{figure}[H]
\centering
\includegraphics[width=0.9\linewidth]{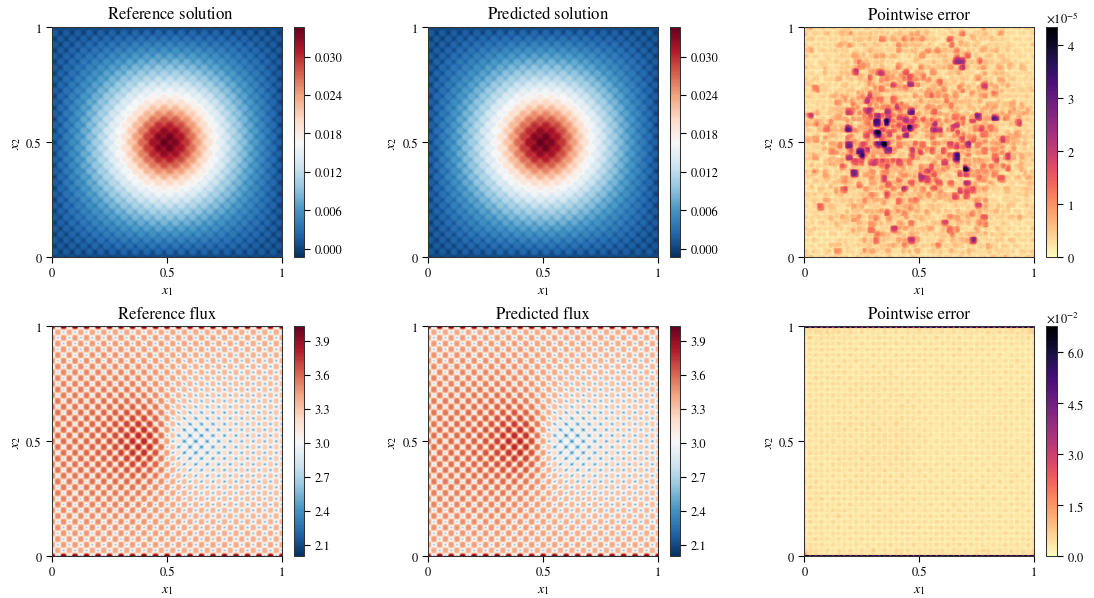}
\caption{Periodic coefficient: Gauss--Newton results.
Top row: FEM reference solution (left), FBPINN prediction
(middle), and pointwise absolute error (right).
Bottom row: FEM reference flux (left), predicted flux
(middle), and pointwise Euclidean norm of the flux error
$\|\boldsymbol{\phi}_{\theta_{\boldsymbol{\phi}}}
-\boldsymbol{\phi}_h\|_2$ (right).}
\label{fig:mixed-periodic-fields}
\end{figure}

\begin{figure}[H]
\centering
\includegraphics[width=1.00\linewidth]{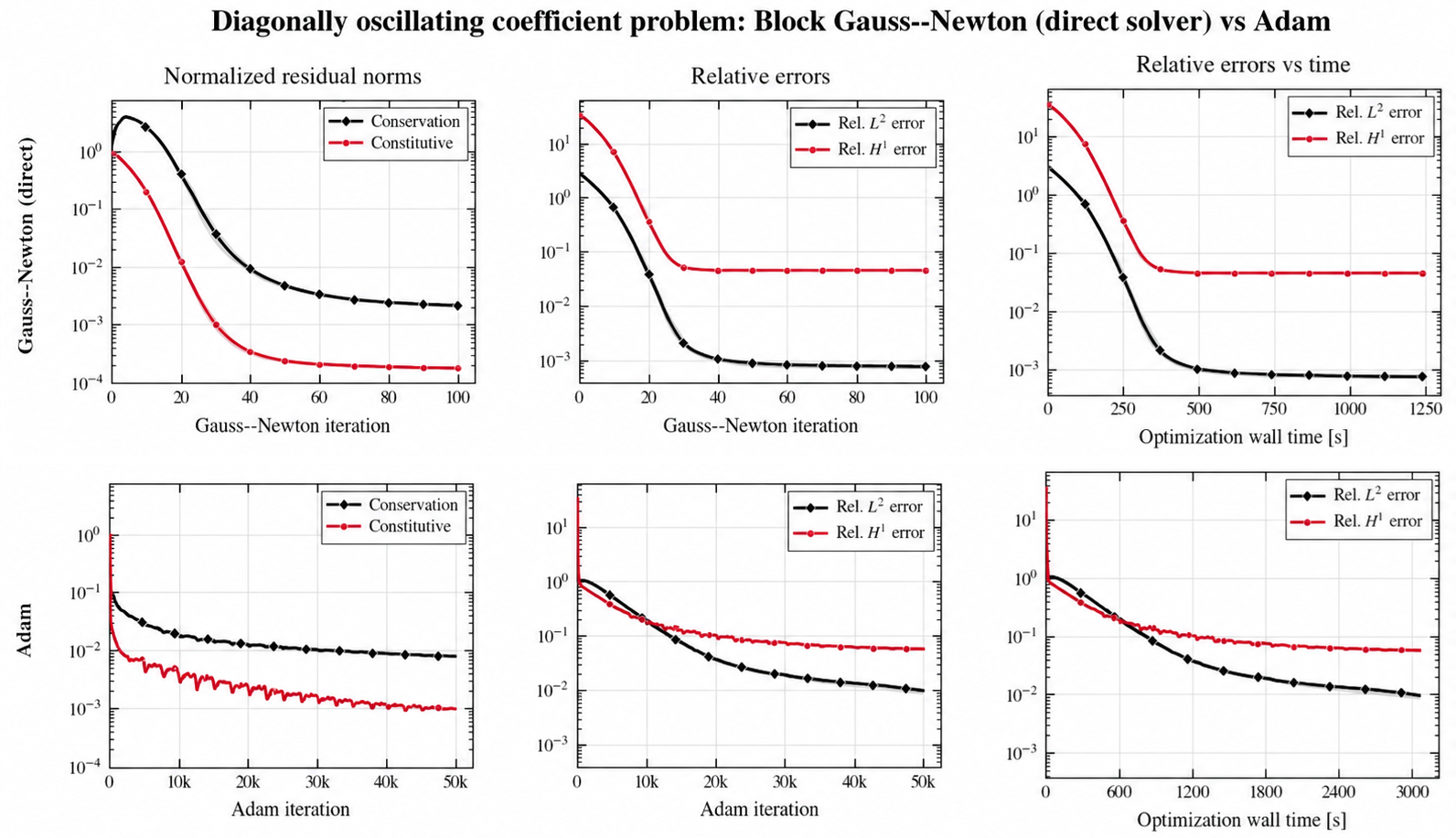}
\caption{Diagonally oscillating coefficient: direct block Gauss--Newton (top) versus
Adam (bottom), with $\delta=2$. Left: conservation and constitutive test
residual norms, each normalized by its own initial value. Middle:
relative $L^2$ and $H^1$-norm errors versus iteration. Right: the same
errors versus optimization wall time. Curves show means over 5 initialization seeds; shaded bands indicate the minimum--maximum range.}
\label{fig:mixed-diagonal-convergence}
\end{figure}

\begin{figure}[H]
\centering
\includegraphics[width=0.9\linewidth]{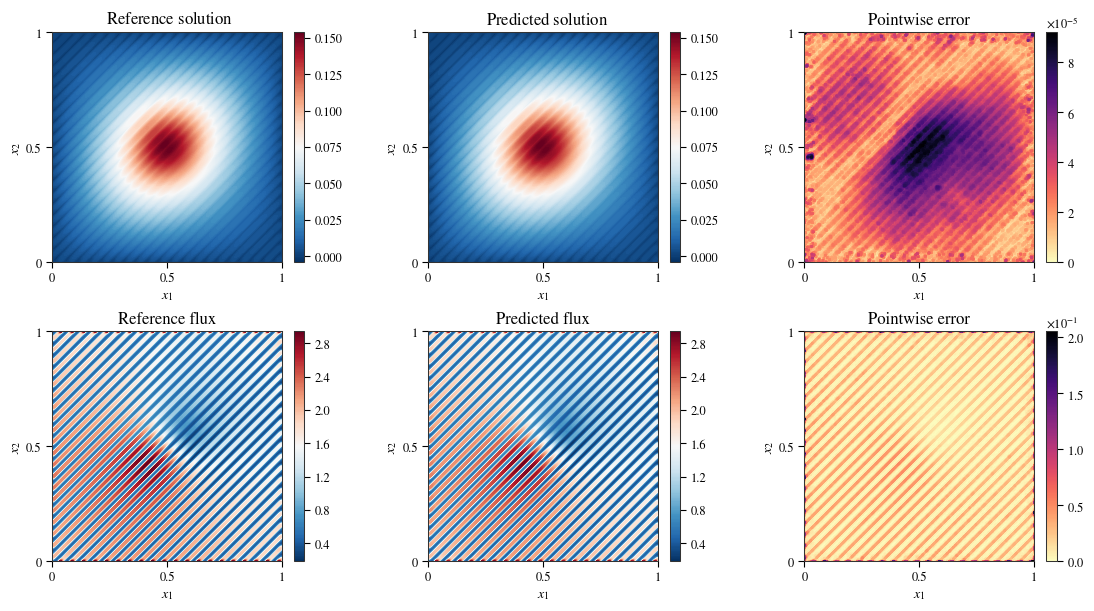}
\caption{Diagonally oscillating coefficient: Gauss--Newton results.
Top row: FEM reference solution (left), FBPINN prediction
(middle), and pointwise absolute error (right).
Bottom row: FEM reference flux (left), predicted flux
(middle), and the pointwise flux error
$\|\boldsymbol{\phi}_{\theta_{\boldsymbol{\phi}}}
-\boldsymbol{\phi}_h\|_2$ (right).}
\label{fig:mixed-diagonal-fields}
\end{figure}

For the periodic coefficient, Fig.~\ref{fig:mixed-periodic-convergence}
shows faster residual and error reduction with Gauss--Newton, both per
iteration and in wall time. Fig.~\ref{fig:mixed-periodic-fields}
compares the resulting solution and flux with the FEM reference (for one representative realization),
including their pointwise errors, showing good agreement. For the diagonally oscillating
coefficient, Fig.~\ref{fig:mixed-diagonal-convergence} shows the same
advantage of Gauss--Newton over Adam. The corresponding solution and flux,
together with their pointwise errors, are shown in
Fig.~\ref{fig:mixed-diagonal-fields}, again showing good agreement with the reference. In both cases, Gauss--Newton
consistently achieves mean relative errors below $10^{-3}$ in the
$L^2$ norm and around $5\%$ in the $H^1$ norm, and outperforms the Adam algorithm. Observe further that FBPINNs exhibit negligible variation in accuracy across random seeds for both Gauss--Newton and Adam.

For both problems, the final Gauss--Newton errors in $L^2$ and the
$H^1$ norm are of the same order as the FEM discrepancies under
further refinement of the reference mesh in Fig.~\ref{fig:coeff_and_fem}.
This indicates that the FBPINN approximations approach the accuracy
level resolved by the available FEM reference.  Table~\ref{tab:mixed-elliptic-accuracy-timing} summarizes the final
accuracy and cost. Gauss--Newton obtains smaller errors than Adam
in about $20.5$ rather than $51.1$ minutes. Sparse AD dominates the
reported mean per-iteration cost; system assembly and the sparse direct
solve are substantially cheaper.

\begin{table}[H]
\centering
\small
\setlength{\tabcolsep}{4pt}

\textbf{(A) Accuracy and total cost}\\[0.35em]
\begin{tabular}{@{}llccc@{}}
\toprule
Coefficient & Optimizer & Relative $L^2$ error & Relative $H^1$ error
  & Total time (min) \\
\midrule
$\boldsymbol K_{\mathrm{diag}}^{\varepsilon}$
  & Sparse GN & $(7.47\pm0.57)\times10^{-4}$
  & $(4.5088\pm0.0003)\times10^{-2}$ & $20.55$ \\
  & Adam & $(9.43\pm1.11)\times10^{-3}$
  & $(5.644\pm0.053)\times10^{-2}$ & $51.08$ \\
\addlinespace[0.4em]
$\boldsymbol K_{\mathrm{per}}^{\varepsilon}$
  & Sparse GN & $(5.09\pm0.26)\times10^{-4}$
  & $(4.893\pm0.004)\times10^{-2}$ & $20.48$ \\
  & Adam & $(3.35\pm0.28)\times10^{-3}$
  & $(7.82\pm0.32)\times10^{-2}$ & $51.13$ \\
\bottomrule
\end{tabular}

\vspace{0.9em}
\textbf{(B) Sparse Gauss--Newton stage timing}\\[0.35em]
\begin{tabular}{@{}lccc@{}}
\toprule
 & Sparse AD & Assembly & Sparse direct solve \\
\midrule
Mean time (s/iteration) & $10.869$ & $0.160$ & $0.924$ \\
\bottomrule
\end{tabular}

\caption{Mixed elliptic benchmarks.
\textbf{(A)} Accuracy and total optimization cost; errors are means $\pm$ sample
standard deviations over 5 seeds.
\textbf{(B)} Mean Gauss--Newton stage timings per iteration, averaged over all
iterations and seeds. Assembly comprises the gradient and the Gramian.}
\label{tab:mixed-elliptic-accuracy-timing}
\end{table}

\section{Conclusion and outlook}

In this work, we have developed a sparse automatic differentiation
framework for efficient Gauss--Newton optimization of FBPINNs.
The approach exploits the block-sparse structure of the residual
Jacobian and the associated Gauss--Newton Gramian induced by the
localization of the subnetworks. For tensor-product subdomain
decompositions, we derived an explicit optimal block coloring
for forward-mode AD and developed a greedy row coloring strategy
for reverse-mode AD based on active-subdomain sets. The resulting
compressed Jacobian representations allow for the assembly of the
Gauss--Newton system directly in block-sparse form, without forming
the full Jacobian. This enables the use of sparse direct and
preconditioned iterative solvers. The efficiency of the proposed
approach is demonstrated by numerical experiments involving
multiscale PDEs. A natural extension concerns more general overlapping subdomain
decompositions. In this setting, closed-form coloring formulas may
no longer be available, requiring efficient greedy coloring strategies
based on active-subdomain sets. Such an extension also necessitates
the construction of suitable window functions adapted to more
general subdomain partitions.

Although sparse direct solvers were more efficient than the
iterative methods in our numerical experiments, the latter are
of particular interest for parallel implementations. In fact,
sparse AD can be performed in batches, individual Gramian blocks
can be assembled independently, and both block matrix--vector
products and local preconditioner solves admit parallelization.
This motivates the development of distributed implementations
on GPU clusters. A further question concerns efficient regularization strategies
based on rank-revealing numerical linear algebra techniques,
which have not been investigated here. By identifying and
eliminating parameter redundancies, such strategies may reduce
ill-conditioning of the Gauss--Newton matrix and its individual
blocks and improve the performance of both sparse direct and
iterative solvers. The development of such regularization strategies,
together with distributed implementations and the investigation
of their scalability, constitutes a subject of our future work in this area.

\section*{Data availability}
The code required to reproduce the results will be available upon the publication of the manuscript.

\section*{Acknowledgements} D. Korolev acknowledges the support of the German Federal Ministry
for Research, Technology and Space  under project HybridSolver (funding reference: 01IS24081).
A. Hill gratefully acknowledges funding from the German Federal Ministry
for Research, Technology and Space under the grant BIFOLD26B.
\appendix

\section{Overlap geometry and window construction}
\label{app1: finite basis networks}
\renewcommand{\thetheorem}{\Alph{section}.\arabic{theorem}}

We specify the grid indexing, uniform cover, and windows used in the main text.

\paragraph{Geometric subdomain grid} Let $K_0$ denote the number of subdomains along each coordinate direction
of the $d$-dimensional uniform tensor-product cover. The total number of
subdomains is therefore $K=K_0^d$. We write
$\mathcal I_{K_0}:=\{1,\ldots,K_0\}^d$ and introduce the row-major bijection
\begin{align}\label{lex ordering}
\operatorname{lex}:\mathcal I_{K_0}\to\{1,\ldots,K\},
\qquad
\operatorname{lex}(a)
:=
1+\sum_{r=1}^{d}(a_r-1)K_0^{d-r}.
\end{align}
For each $\al\in\{1,\ldots,K\}$, let
\begin{align*}
a_\al
=
(a_{\al,1},\ldots,a_{\al,d})
:=
\operatorname{lex}^{-1}(\al)
\end{align*}
denote the position of $\Omega_\al$ on the grid. Thus, $\al$ identifies
the corresponding subdomain and parameter block, whereas $a_\al$
specifies its geometric location.

\begin{figure}[H]
\centering
\includegraphics[width=1.0\linewidth]{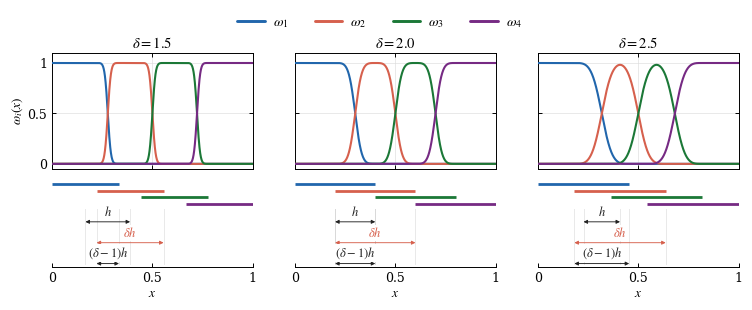}
\caption{Normalized one-dimensional windows \eqref{eq:one-dimensional-pou} on
$[0,1]$ with $K_0=4$, for $\delta=1.5,2.0,2.5$. Bars show the supports $I_i$,
annotated with the spacing $h$, width $\delta h$ and overlap $(\delta-1)h$. At
most $p_\delta=\lceil\delta\rceil$ windows are nonzero at a point: two for
$\delta\le2$, three for $2<\delta\le3$.}
\label{fig:window-functions}
\end{figure}

\paragraph{Overlap ratio}
Fix $\delta>1$. In coordinate direction $j$, define the center spacing and
interval centers by
\begin{align}\label{eq:no-clipping-spacing}
h_j
&:=\frac{x_j^{\max}-x_j^{\min}}{K_0-1+\delta},
&
c_i^{(j)}
&:=x_j^{\min}+\left(i-1+\frac{\delta}{2}\right)h_j,
\qquad i=1,\ldots,K_0.
\end{align}
The one-dimensional intervals and tensor-product boxes are given by
\begin{align}\label{eq:no-clipping-1d-interval}
I_i^{(j)}
:=\left[c_i^{(j)}-\frac{\delta h_j}{2},\;
         c_i^{(j)}+\frac{\delta h_j}{2}\right], \qquad \Omega_\alpha
:=\prod_{j=1}^d I_{a_{\alpha,j}}^{(j)}.
\end{align}
For adjacent intervals,
\begin{align}\label{eq:overlap-ratio-definition}
|I_i^{(j)}|=\delta h_j,
\qquad
c_{i+1}^{(j)}-c_i^{(j)}=h_j,
\qquad
|I_i^{(j)}\cap I_{i+1}^{(j)}|=(\delta-1)h_j.
\end{align}
Thus, $\delta$ is the ratio of subdomain width to center spacing. Fig.~\ref{fig:window-functions} illustrates the
present placement.

\paragraph{Window construction}
Following \cite[Section~2.2.1, Eq.~(10)]{dolean2024multilevel}, we use
the compactly supported squared raised-cosine profile
\begin{align}\label{eq:raw-window-profile}
W(t):=
\begin{cases}
\displaystyle\left(\frac{1+\cos(\pi t)}{2}\right)^2,
& |t|\leq1,\\[1ex]
0,& |t|>1.
\end{cases}
\end{align}
Its zero extension belongs to $C^3(\R)$, which suffices for the first- and
second-order PDE residuals considered here. For
$x_j^{\min}<t<x_j^{\max}$, set
\begin{align}\label{eq:one-dimensional-pou}
\widehat\omega_i^{(j)}(t)
&:=W\!\left(\frac{2(t-c_i^{(j)})}{\delta h_j}\right),
&
\omega_i^{(j)}(t)
&:=\frac{\widehat\omega_i^{(j)}(t)}
{\sum_{k=1}^{K_0}\widehat\omega_k^{(j)}(t)}.
\end{align}
The denominator is positive in the open interval because $\delta>1$.
At the physical endpoints, the normalized windows are defined by continuity:
the first window equals one at $x_j^{\min}$, the last equals one at
$x_j^{\max}$, and all other endpoint values are zero.
The $d$-dimensional windows are given by
\begin{align}\label{eq:tensor-product-window}
\omega_\alpha(x)
:=\prod_{j=1}^d\omega_{a_{\alpha,j}}^{(j)}(x_j).
\end{align}
Since the one-dimensional weights sum to one, this construction satisfies
\eqref{eq: partition of unity}.

\begin{lemma}\label{lem:overlap-multiplicity}
The tensor-product windows \eqref{eq:tensor-product-window} satisfy
\begin{align}
|T(x)|\leq p_\delta^d,
\qquad x\in\Omega.
\end{align}
\end{lemma}

\begin{proof}
Fix $j\in\{1,\ldots,d\}$ and suppose that
$\omega_i^{(j)}(x_j)\neq0$ and $\omega_k^{(j)}(x_j)\neq0$
for some $i<k$. By \eqref{eq:no-clipping-1d-interval},
\begin{align*}
x_j^{\min}+(k-1)h_j
<x_j<
x_j^{\min}+(i-1+\delta)h_j,
\end{align*}
which implies $k-i<\delta$, and hence
$k-i\leq r_\delta$. Therefore, at most
$p_\delta=r_\delta+1$ one-dimensional windows are active
in each coordinate direction. The tensor-product structure
\eqref{eq:tensor-product-window} gives the assertion.
\end{proof}

\section{Proofs of the coloring results}
\label{app:coloring-proofs}

\subsection{Proof of Proposition~\ref{prop:block-coloring}}
\label{app:proof-block-coloring}

\begin{proof}
For each subdomain index $\al$, let $a_{\alpha} = (a_{\al,1},\dots, a_{\al,d}) \in \mathcal I_{K_0}$ be its grid index given by \eqref{lex ordering}. For $j=1,\dots, d$, define the coordinate residue
\begin{align}
\rho_{j}(a_{\al}):=(a_{\al, j}-1) \bmod p_{\delta} \in \{0,\dots, p_{\delta} - 1 \}
\end{align}
Thus, each grid index $a_{\al}$ determines the residue vector
\begin{align*}
\rho(a_{\al}) = (\rho_{1}(a_{\al}), \dots, \rho_{d}(a_{\al})) \in \{0, \dots, p_{\delta} - 1\}^{d}.
\end{align*}
The coloring is obtained by scalarizing this residue vector. Indeed, \eqref{block coloring map} is the bijective base-$p_{\delta}$ linearization of $\rho(a_{\al})$. Therefore, two subdomain labels $\al$ and $\be$ have the same color if and only if  $\rho(a_{\al}) = \rho(a_{\be})$.

We now prove that $\kappa$ is a valid coloring. Suppose that two vertices $\alpha$ and $\beta$  in $\mathfrak G_{\mathrm{blk}}$ have the same color. Then their residue vectors agree, which implies that $a_{\alpha,j}-a_{\beta,j}$ is a multiple of $p_\delta$ for every coordinate direction $j$. Assume, for contradiction, that $\alpha\neq\beta$ and that $\alpha$ and $\beta$ are adjacent in $\mathfrak G_{\mathrm{blk}}$. In the block graph, this implies
\begin{align}
|a_{\alpha,j}-a_{\beta,j}|
\leq
r_\delta = p_\delta-1,
\qquad
j=1,\ldots,d.
\end{align}
The only multiple of $p_\delta$ whose absolute value is at most $p_\delta-1$ is zero. Therefore,
\begin{align}
a_{\alpha,j}= a_{\beta,j}, \qquad j=1,\ldots,d,
\end{align}
which implies that $a_\alpha=a_\beta$. Since the lexicographic map \eqref{lex ordering} is bijective, we get $\alpha=\beta$, contradicting $\alpha\neq\beta$. Therefore, no two distinct adjacent vertices receive the same color, and $\kappa$ is a valid coloring.

It remains to prove optimality. Since $K_0\geq p_\delta$, the tensor-product grid $I_{K_0}$ contains the subgrid $\{1,\ldots,p_\delta\}^d \subset \mathcal I_{K_0}$. Any two distinct multi-indices in this subgrid differ by at most $p_\delta-1 = 
r_\delta$ in every coordinate direction. Hence, the corresponding vertices are pairwise adjacent in $\mathfrak G_{\mathrm{blk}}$. Therefore, this subgrid induces a clique of size $p_\delta^d$  in $\mathfrak G_{\mathrm{blk}}$, which implies that every valid coloring of $\mathfrak G_{\mathrm{blk}}$ requires at least $p_\delta^d$ colors. Since the coloring $\kappa$ uses exactly $p_\delta^d$ colors, it is optimal, and \eqref{Chrom number} holds.
\end{proof}



\bibliographystyle{elsarticle-num} 
\bibliography{bib.bib}
\end{document}